\documentclass[reqno,twoside,11pt,english]{amsart}
\usepackage{amsmath,amsfonts,amssymb,amsthm,epsfig}
\newcommand{\RNum}[1]{\uppercase\expandafter{\romannumeral #1\relax}}
\usepackage{comment}
\usepackage{color}
\usepackage[final,linkcolor = blue,citecolor = blue,colorlinks=true]{hyperref}
\usepackage[leqno]{amsmath}

\usepackage{amscd}
\usepackage{mathrsfs}
\usepackage{}
\usepackage[all]{xy}
\usepackage[OT1]{fontenc}
\usepackage[latin1]{inputenc}
\usepackage{amssymb,latexsym}
\usepackage{type1cm,ae}
\usepackage{graphicx}
\usepackage{appendix}
\usepackage{xspace}
\usepackage{setspace}
\usepackage[english]{babel}

\theoremstyle{definition}
		\newtheorem{theorem}{Theorem}[section]
				\newtheorem{proposition}[theorem]{Proposition}
				\newtheorem{lemma}[theorem]{Lemma}
				\newtheorem{thmA}{Theorem}[section]
				
                \newtheorem{example}[theorem]{Example}
                \newtheorem{question}[theorem]{Question}

     	        \newtheorem{definition}[theorem]{Definition}
	            \newtheorem{remark}[theorem]{Remark}
\numberwithin{equation}{section}
\newcommand*{\bR}{\ensuremath{\mathbb{R}}}

\newcommand*{\loc}{\mathrm{loc}}
\newcommand*{\closure}[1]{\overline{#1}}
\newcommand*{\bdary}[1]{\partial #1}

\newcommand*{\Wert}{\mathord{\mbox{|\kern-1.5pt|\kern-1.5pt|}}}
\newcommand*{\ie}{\mbox{i.e.}\xspace}

\newcommand\numberthis{\addtocounter{equation}{1}\tag{\theequation}}

\newcommand{\N}{\mathbb{N}}	% Numeri naturali
\newcommand{\R}{\mathbb{R}}	% Numeri reali

\newcommand{\G}{\mathbb{G}}	% Grassmanniane
\newcommand{\Co}{\mathscr{C}}	% Funzioni continue
\newcommand{\X}{\mathfrak{X}}	% Spazio dei campi vettore
\renewcommand{\Vec}{\mathrm{Vec}}	% Spazio dei campi vettore
\newcommand{\Span}{\mathrm{span}}	% Span
\newcommand{\dd}{\,\mathrm{d}}	% 'd' di derivata
\newcommand{\de}{\partial}		% Derivata parziale
\newcommand{\THEN}{\Rightarrow}	% =>
\newcommand{\IFF}{\Leftrightarrow}	% <=>
\DeclareMathOperator{\diam}{diam}

\DeclareMathOperator{\lip}{lip}
\DeclareMathOperator{\Lip}{Lip}
\DeclareMathOperator{\Boxs}{Box}

\DeclareMathOperator{\Vol}{Vol}
\DeclareMathOperator{\card}{card}

\newcommand{\scr}[1]{\mathscr{#1}}
\newcommand{\frk}[1]{\mathfrak{#1}}

\newcommand{\g}{\frk g}

\def\Xint#1{\mathchoice
  {\XXint\displaystyle\textstyle{#1}}%
  {\XXint\textstyle\scriptstyle{#1}}%
  {\XXint\scriptstyle\scriptscriptstyle{#1}}%
  {\XXint\scriptscriptstyle\scriptscriptstyle{#1}}%
  \!\int}
\def\XXint#1#2#3{{\setbox0=\hbox{$#1{#2#3}{\int}$}
  \vcenter{\hbox{$#2#3$}}\kern-.5\wd0}}

\def\dashint{\Xint-}
\def\blue#1{\textcolor{blue}{#1}}

\makeatother

\usepackage{babel}
\begin{document}

\title[Quasiregular mappings between equiregular SubRiemannian manifolds]{Quasiregular mappings between equiregular SubRiemannian manifolds}
\author{Chang-Yu Guo}
\date{\today}

\address[Chang-Yu Guo]{Research Center for Mathematics, Shandong University, Qingdao, P. R. China and Frontiers Science Center for Nonlinear Expectations, Ministry of Education, P. R. China and Department of Physics and Mathematics, University of Eastern Finland, 80101, Joensuu, Finland}
\email{guocybnu@gmail.com}

\author{Sebastiano Nicolussi Golo}
\address[Sebastiano Nicolussi Golo]{Department of Mathematics, University of Fribourg, Ch. du Mus\'ee 23, 1700 Fribourg, Switzerland}
\email{sebastiano2.72@gmail.com}

\author{Marshall Williams}

\author{Yi Xuan}
\address[Yi Xuan]{Research Center for Mathematics, Shandong University, Qingdao, P. R. China and Frontiers Science Center for Nonlinear Expectations, Ministry of Education, P. R. China}
\email{yixuan@sdu.edu.cn}

\thanks{C.-Y. Guo and Y. Xuan are supported by the Young Scientist Program of the Ministry of Science and Technology of China (No.~2021YFA1002200), the National Natural Science Foundation of China (No.~12101362), the Taishan Scholar Project and the Natural Science Foundation of Shandong Province (No.~ZR2022YQ01).}

\begin{abstract}
In this paper, we provide an alternative appraoch to an  expectaion of F\"assler et al [J. Geom. Anal. 2016] by showing that a metrically quasiregular mapping between two equiregular subRiemannian manifolds of homogeneous dimension $Q\geq 2$ has a negligible branch set. One main new ingredient is to develop a suitable extension of the generalized Pansu differentiability theory, in spirit of earlier works by Margulis-Mostow, Karmanova and Vodopyanov. Another new ingredient is to apply the theory of Sobolev spaces based on upper gradients developed by Heinonen, Koskela, Shanmugalingam and Tyson to establish the necessary analytic foundations.
\end{abstract}

\maketitle
\tableofcontents{}

\section{Introduction}\label{sec:introduction}
Let $M$ and $N$ be two equiregular subRiemannian manifolds of homogeneous dimension $Q\geq 2$. Recall that an equiregular subRiemannian manifold is a triple $(M,H,g)$ where $M$ is a smooth and connected $n$-manifold, $H\subset TM$ is a bracket generating and equiregular subbundle, and $g:H\times H\to[0,+\infty)$ is a smooth function whose restriction to each fiber $H_p$ is a scalar product (see Section \ref{subsec:geometry of equiregular SRM} below for precise definition). The homogeneous dimension $Q=\sum_{i=1}^sin_i$ is typically bigger than the topological dimension. Moreover, the Hausdorff measure $\mathcal{H}^Q$ (with respect to the subRiemannian distance) is equivalent to the Lebesgue measure on $M$ (see \cite{cambridge}). 

A map $f\colon M\to N$ between two equiregular subRiemannian manifolds of homogeneous dimension $Q\geq 2$ is said to be a \textit{branched cover} if it is continuous, discrete, open and sense-preserving. Here a mapping $f\colon X\to Y$ between (oriented) topological manifolds is \textit{discrete} if each fiber is a discrete set in $X$, \ie for all $y\in Y$, $f^{-1}(y)$ is a discrete set in $X$, and is \textit{open} if it maps open set in $X$ onto open set in $Y$, and is \textit{sense-preserving} if the local index (or degree) at each image point is positive.  For a branched cover $f:M\to N$, the \textit{linear dilatation} $H_f(x)$ of $f$ at $x\in M$ is defined to be
\begin{equation}\label{eq:def for linear dilatation}
	H_f(x)=\limsup_{r\to 0} H_f(x,r)=\limsup_{r\to 0}\frac{L_f(x,r)}{l_f(x,r)},
\end{equation}
where
\begin{align*}
	L_f(x,r)=\sup\big\{d(f(x),f(y)):y\in \closure{B}(x,r)\big\}
\end{align*}
and
\begin{align*}
	l_f(x,r)=\inf\big\{d(f(x),f(y)):y\in\bdary B(x,r)\big\}.
\end{align*}
If $f\colon M\to N$ is a homeomorphism such that $H_f(x)\leq H$ for all $x\in M$, then $f$ is said to be a \textit{metrically $H$-quasiconformal mapping}.  

In this paper, we are mainly interested in the non-injective class of \emph{metrically quasiregular mappings}, that is, a branched cover $f\colon M\to N$ such that $H_f(x)<\infty$ for all $x\in M$ and $H_f(x)\leq H$ for almost every $x\in M$. The difference between a metrically quasiregular mapping and a metrically quasiconformal mapping lies in the mysterious topological branch set $\mathcal{B}_f$, i.e. the set of all $x\in M$ such that $f$ fails to be a local homeomorphism at $x$.

Quasiconformal/quasiregular mappings between Euclidean domains form a natural extension of the class of planar conformal/analytic mappings and they have wide connections to many other areas such as partial differential equations, geometric topology, geometric group theory and so on; see \cite{v71,re89,r93,im01,m14-survey,gmp17} for a comprehensive introduction of the theory and their profound applications in other branches of mathematics.

The study of quasiconformal mappings beyond Riemannian spaces has first appeared in the celebrated work of Mostow~\cite{m73} on strong rigidity of locally symmetric spaces. The boundary of rank-one symmetric spaces can be identified as certain Carnot groups of step two, and Mostow has developed the basic (metrically) quasiconformal mapping theory in these groups. Inspired by Mostow's work, Pansu~\cite{p89} used the theory of quasiconformal mappings to study quasi-isometries of rank-one symmetric spaces. The systematic study of (metrically) quasiconformal mappings on the Heisenberg group was later done by Kor\'anyi and Reimann~\cite{kr85,kr95}.  Margulis and Mostow~\cite{mm95} studied the absolute continuity of quasiconformal mappings along horizontal lines in the equiregular subRiemannian case. By the break-through work of Heinonen and Koskela~\cite{hk98} (see also \cite{hk95} for the case of Carnot groups), a full-fledged (metrically) quasiconformal mapping theory exists in rather general metric measure spaces.  This theory has subsequently been applied to new rigidity studies in geometric group theory, geometric topology and geometric parametrization of metric spaces; see for instance~\cite{bk02,bk05,hr02,r14} and the references therein. This theory also initiated a new way of looking at weakly differentiable maps between non-smooth spaces. In~\cite{hkst01},  the Sobolev class of Banach space valued mappings was studied and several characterizations of quasiconformal mappings between metric spaces of locally bounded geometry were established.

Along with the development of quasiconformal mappings, the theory of quasiregular mappings also received great interest in nonsmooth spaces. The study of quasiregular mappings between equiregular subRiemannian manifolds was initiated by Heinonen and Holopainen~\cite{hh97}, where the authors developed the theory of quasiregular mappings between Carnot groups from an analytic point of view. In the special setting of Heisenberg groups, Dairbekov~\cite{d99,d00} enhanced the results of Heinonen and Holopainen by considering quasiregular mappings under natural Sobolev regularity assumptions. There are also successful development of the theory in the setting of other Carnot groups by Vodopyanov~\cite{v07}. In the setting of generalized manifolds with restricted topology and controlled geometry, Heinonen and Rickman~\cite{hr02} studies the so-called mappings of bounded length distortion (BLD), which form a proper subclass of quasiregular mappings. This theory has been successfully applied in the geometric parametrization problem by Heinonen and Sullivan \cite{hs02}; see also~\cite{t98,bkr07,bfp12,flp14,g14,v07,w12proc,w14,gw14,gw16-paris,gw16,kmx21} for the latest development of the theory in various non-smooth settings. In particular, the quantitative equivalences of all the three definitions of quasiregularity in the setting of metric spaces with locally bounded geometry has been established in~\cite{gw16}. Note that Riemannian manifolds and Carnot groups are typical examples of metric spaces with locally bounded geometry.

There are also many other nice spaces that do not necessary have bounded geometry, but are of independent interest to develop interesting analysis and geometry. Among these nice spaces, we are particularly interested in equiregular subRiemannian manifolds, which serve as a natural class of singular spaces that lie between Riemannian manifolds and general metric measure spaces. In particular, the geometry of an equiregular subRiemannian manifold is meaningful: the tangent cone at each point admits a natural group structure that makes it into the so-called Carnot group and a well-known differentiation theory of mappings exists for mappings between such groups. Very recently, the study of uniformly quasiregular mappings in the setting of Heisenberg groups, generates new interest of a theory of quasiregular mappings between equiregular subRiemannian manifolds from a metric point of view; see~\cite{bfp12,flp14,flt19}. In particular, in \cite{flp14}, a stronger definition of metrically quasiregular mappings was introduced and in \cite[Remark 1.2]{flp14}, the authors made the following natural expectation:
\medskip 

\textbf{Expectation:} \emph{Both the branch set $\mathcal{B}_f$ and its image $f(\mathcal{B}_f)$ of a metrically quasiregular mapping $f:M\to N$ between two equiregular subRiemannian manifolds have measure zero.} 
\medskip 

The expectation turned out to be affirmative by the recent work of Guo and Williams \cite{gw14,gw16-paris}. However, the approach of Guo and Williams uses advanced tools from quantitative topology and the proof is much more complicated than the Euclidean case. 
The aim of this paper is to provide an alternative approach to the above expectation and develop analytic foundations of the theory of metrically quasiregular mappings between equiregular subRiemannian manifolds. More precisely, we shall confirm the above expectation by showing the following theorem. 

\begin{thmA}\label{thm:branch set}
	Let $f\colon M\to N$ be a weakly metrically quasiregular mapping between two equiregular subRiemannian manifolds of homogeneous dimension $Q\geq 2$. Then $$\Vol_M(\mathcal{B}_f)=\Vol_M(f(\mathcal{B}_f))=0.$$
\end{thmA}

The idea behind the proof is in spirit similar to that used in the Euclidean setting and relies on the almost everywhere \emph{geometric differentiability} of the map $f$ and the almost everywhere positiveness of the Jacobian determinant $J_f$. At a point $x$ of differentiability, we may use the differential $Df(x)$ to approximate $f$ in a neighborhood of $x$. If $J_f(x)>0$, then $Df(x)$ will be a Carnot group isomorphism and so we may construct a suitable homotopy between $f$ and the identity map. Then standard topological degree theory would imply that $f$ is a local homeomorphism around $x$. Since such points have full measure, we conclude $\Vol_M(\mathcal{B}_f)=0$ as desired.   

To make the above idea work, we shall develop a theory of geometric differentiation for mappings between equiregular subRiemannian manifolds and prove that (weakly) metrically quasiregular mappings are $P$-differentiable almost everywhere in $M$. 

\begin{thmA}\label{thm:Differentiability Stepanov}
	Let $f\colon (M,d)\to (\bar M,\bar d)$ be a Borel mapping between two equiregular subRiemannian manifolds.
	Then $f$ is $P$-differentiable for almost every~$o$ in the set
	\[
	L(f) := \left\{o\in M: \limsup_{p\to o} \frac{\bar d(f(o),f(p))}{d(o,p)} < \infty \right\}.
	\]
\end{thmA}

The definition of $P$-differential is given in Section~\ref{sec:Pansu diff} and it is a natural extension of the notion introduced by Pansu~\cite{p89} and Margulis-Mostow~\cite{mm95}. It should be noticed that one cannot use the standard techniques as in the Euclidean setting to prove Theorem~\ref{thm:Differentiability Stepanov}, since it is not always possible to extend a Lipschitz mapping $f:A\to N$ from a closed subset $A$ of $M$ as a global Lipschitz mapping $\hat{f}:M\to N$. Our proof of Theorem~\ref{thm:Differentiability Stepanov} relies on a careful blow-up argument, which seems to be new even in the setting of Carnot groups. As the proof is rather technical, we present a simplified proof in the setting of Carnot groups in the appendix and point out the principal differences from general equiregular subRiemannian case. 

\begin{remark}\label{rmk:on thm stepanov}
There are many other extensions of the definition of Pansu differentiability and Stepanov's theorem.  For instance, Vodopyanov \cite{v07diff} announced a proof of Stepanov's theorem for mappings between equiregular subRiemannian manifolds. However, we have difficulty in understanding the proofs. For example, Property 1.3 in \cite{v07diff} states that the exponential mappings at $g$ of the linear combination of the vector fields $X_i$ and the linear combination of $\hat{X}_i^g$ with the same coefficient are equal for the coefficients which are sufficiently small. But, this seems not correct for us. Although $X_i$ and $\hat{X}_i^g$ are the same at $g$, there is no hope that they are the same everywhere. Furthermore, the exponential mapping for a vector field depends not only on the value of vector field at $g$ but also on the value at other points. Thus, there is no evident to show the above two exponential mappings are equal even for sufficiently small coefficients. 

In another paper \cite{kv09}, Karmanova and Vodopyanov announced a similar proof of Stepanov's theorem as \cite{v07diff}. The arguments in both papers are rather similar and thus we cannot understand it. For instance, in \cite[Lemma 2.1.26]{kv09}, they used a similar result as \cite[Property 1.3]{v07diff}, which seems again incorrect for us. Thus we decided to include a self-contained and different proof. 
\end{remark}

Combining Theorem~\ref{thm:Differentiability Stepanov} with some ideas from the recent development of quasiregular mappings on metric measure spaces~\cite{hkst01,w14,gw16}, we are able to establish the necessary analytic foundations for theory of metrically quasiregular mappings between equiregular subRiemannian manifolds (see Theorem \ref{thm:properties} below).

As the previous overview indicates, to establish the theory of metrically quasiregular mappings between general metric measure spaces, the correct conditions imposing on the underlying spaces seems to be that the spaces have locally bounded geometry. Our approach seems to suggest that when the underlying metric spaces have nice geometry so that a (geometric) differentiability theory for mappings exists, then the basic properties of metrically quasiregular mappings remain valid for spaces without bounded geometry. It is then an interesting problem to investigate to what extent the differentiability theory helps in establishing other useful properties of quasiregular mappings; see Section~\ref{subsec:open questions} for those natural open problems.

This paper is organized as follows. Section~\ref{sec:Sobolev spaces on mms} contains preliminaries on Sobolev spaces on metric measure spaces based on upper gradients, Carnot groups, subRiemannian manifolds and tangent cones. Section~\ref{sec:QR mapping} contains a detailed study of metrically quasiregular mappings between equiregular subRiemannian manifolds. In particular, we investigate all the basic analytic properties of metrically quasiregular mappings from the differentiable point of view and prove Theorem~\ref{thm:properties}. In Section~\ref{sec:Pansu diff}, we generalize the notion of $P$-differentiability of Pansu~\cite{p89} and Margulis-Mostow~\cite{mm95} to the setting of mappings between general subRiemannian manifolds and prove Theorem~\ref{thm:Differentiability Stepanov}. In the appendix, we include a proof of Theorem \ref{thm:Differentiability Stepanov} in the setting of Carnot groups, to help the readers to understand the idea of our approach.

\section{Preliminaries}\label{sec:Sobolev spaces on mms}

\subsection{Sobolev spaces}\label{subsec:metric measure spaces}

In this subsection, we will briefly introduce the Sobolev spaces on metric measure spaces based on an upper gradient approach. For detailed description of this approach, see the monograph~\cite{hkst12}.

\begin{definition}\label{def:metric measure space}
A \textit{metric measure space} is defined to be a triple $(X,d,\mu)$, where $(X,d)$ is a separable metric space and $\mu$ is a nontrivial locally finite Borel regular measure on $X$.
\end{definition}

\begin{definition}\label{def:doubling metric measure space}
A Borel regular measure $\mu$ on a metric space $(X,d)$ is called a \textit{doubling measure} if every ball in $X$ has positive and finite measure and there exists a constant $C_\mu\geq 1$ such that
\begin{equation}\label{eq:doubling measure}
\mu(B(x,2r))\leq C_\mu \mu(B(x,r))
\end{equation}
for all balls $B(x,r)\subset X$ with radius $r<\diam X$. We call the triple $(X,d,\mu)$ a doubling metric measure space if $\mu$ is a doubling measure on $X$. We call $(X,d,\mu)$ an \textit{Ahlfors $Q$-regular} space if there exists a constant $C\geq 1$ such that
\begin{equation}\label{eq:Ahlfors regular measure}
C^{-1}r^Q\leq \mu(B(x,r))\leq Cr^Q
\end{equation}
for all balls $B(x,r)\subset X$ of radius $r<\diam X$.
\end{definition}

Let $X=(X,d,\mu)$ be a metric measure space and $Z=(Z,d_Z)$ be a metric space.

\begin{definition}\label{def:upper gradient}
A Borel function $g\colon X\rightarrow [0,\infty]$ is called an \textit{upper gradient} for a map $f:X\to Z$ if for every rectifiable curve $\gamma:[a,b]\to X$, we have the inequality
\begin{equation}\label{ugdefeq}
\int_\gamma g\,ds\geq d_Z(f(\gamma(b)),f(\gamma(a)))\text{.}
\end{equation}
If inequality \eqref{ugdefeq} merely holds for $p$-almost every compact curve\footnote{See \cite{hkst12} for definition of a property that holds on $p$-almost every curve}, then $g$ is called a \textit{$p$-weak upper gradient} for $f$.  When the exponent $p$ is clear, we omit it.
\end{definition}
The concept of upper gradient was introduced in~\cite{hk98}. It was initially called ``very weak gradient", but the befitting term ``upper gradient" was soon suggested. Functions with $p$-integrable $p$-weak upper gradients were subsequently studied in~\cite{km98}, while the theory of Sobolev spaces based on upper gradient was systematically developed in~\cite{s00} and~\cite{c99}.

A $p$-weak upper gradient $g$ of $f$ is \textit{minimal} if for every $p$-weak upper gradient $\tilde{g}$ of $f$, $\tilde{g}\geq g$ $\mu$-almost everywhere.  If $f$ has an upper gradient in $L^p_\loc(X)$, then $f$ has a unique (up to sets of $\mu$-measure zero) minimal $p$-weak upper gradient by~\cite[Theorem \blue{6.3.20}]{hkst12}.  In this situation, we denote the minimal upper gradient of $f$ by $g_{f}$. The minimal $p$-weak upper gradient $g_f$ should be thought of as a substitute for $|\nabla f|$, or the length of a gradient, for functions defined in metric measure spaces.

Let $\mathbb{V}$ be a Banach space and $\tilde{N}^{1,p}(X,\mathbb{V})$ denote the collection of all maps $u\in L^p(X,\mathbb{V})$ that have an upper gradient in $L^p(X)$. We equip it with seminorm
\begin{equation}
\|u\|_{\tilde{N}^{1,p}(X,\mathbb{V})}=\|u\|_{L^p(X,\mathbb{V})}+\|g_u\|_{L^p(X)},
\end{equation}
where $g_u$ is the minimal $p$-weak upper gradient of $u$.

We obtain a normed space $N^{1,p}(X,\mathbb{V})$ by passing to equivalence classes of functions in $\tilde{N}^{1,p}(X,\mathbb{V})$, where $u_1\sim u_2$ if and only if $\|u_1-u_2\|_{\tilde{N}^{1,p}(X,\mathbb{V})}=0$. Thus
\begin{equation}\label{eq:definition of Nowton-Sobolev space}
N^{1,p}(X,\mathbb{V}):=\tilde{N}^{1,p}(X,\mathbb{V})/\{u\in \tilde{N}^{1,p}(X,\mathbb{V}): \|u\|_{\tilde{N}^{1,p}(X,\mathbb{V})}=0\}.
\end{equation}

Let $\tilde{N}_{\loc}^{1,p}(X,\mathbb{V})$ be the vector space of (Banach-space valued) functions $u:X\to \mathbb{V}$ with the property that every point $x\in X$ has a neighborhood $U_x$ in $X$ such that $u\in \tilde{N}^{1,p}(U_x,\mathbb{V})$. Two functions $u_1$ and $u_2$ in $\tilde{N}_{\loc}^{1,p}(X,\mathbb{V})$ are said to be equivalent if every point $x\in X$ has a neighborhood $U_x$ in $X$ such that the restrictions $u_1|_{U_x}$ and $u_2|_{U_x}$ determine the same element in $\tilde{N}^{1,p}(U_x,\mathbb{V})$. The local Sobolev space $N_{\loc}^{1,p}(X,\mathbb{V})$ is the vector space of equivalent classes of functions in $\tilde{N}_{\loc}^{1,p}(X,\mathbb{V})$ under the preceding equivalence relation.

To define the Sobolev space $N^{1,p}(M, N)$ of mappings $f: M\to N$ (between two subRiemannian manifolds $M$ and $N$), we first fix an isometric embedding $\varphi$ of $N$ into the Banach space $\mathbb{V}=\ell^\infty(N)$. Then the Sobolev space $N^{1,p}(M,N)$ consists of all mappings $f:M\to N$ with $\varphi\circ f\in N^{1,p}(M,\mathbb{V})$ and $\varphi\circ f\in N$ almost everywhere.

\subsection{Carnot groups}\label{subsec:Carnot group}
A Carnot group is a \emph{simply connected nilpotent Lie group $G$} with a graded Lie algebra $$ \mathfrak{g}=V_1 \oplus \cdots \oplus V_s $$ equipped with an inner product on $V_1$. Note that $V_1$ generates $\mathfrak{g}$ as a Lie algebra. As a Lie group, there is a canonical left multiple action on $G$. The left multiple by the element $o$ is denoted by $L_o$. Under that action, we obtain a graded Lie algebra structure on the tangent space $\g_o$ of all points $o\in G$. In the canonical way, we can associate $\g_o$ with a Carnot group $\G_o$. What is more, we obtain an inner produce on the corresponding vector subspace at $o$, which is called $V_1^{o}$. 

We call a curve $\gamma:[0,1] \rightarrow G$ \emph{horizontal} if the tangent vector $\gamma^{\prime}(t)$ lies in $V_1^{\gamma(t)}$ at the point $\gamma(t)$ for almost every $t$. Then, for $x, y \in G$, we define 
$$ d(x, y)=\inf \int_0^1\left\|\gamma^{\prime}(t)\right\| d s $$ where the infimum is taken over all horizontal curves $\gamma$ such that $\gamma(0)=x$ and $\gamma(1)=y$, and where $\|\cdot\|$ denotes the inner product of $V_1^{\gamma(t)}$ at the point $\gamma(t)$. 

Every two points on $G$ can be joined by a horizontal curve, such that $d$ defines a distance on $G$. Every two points on $G$ can be joined by a minimal horizontal curve. Furthermore, $d$ determines a left-invariant metric on $G$. Moreover, there are dilations associated to all the points of $G$ as follows. Carnot groups are viewed as $\R^n$ by the exponential map at the unit of the Lie group. For the point $0$, we define the traditional dilation $\delta_{\epsilon}^0(x_1,x_2,...,x_n)=(\epsilon^{\tau(1)}x_1,\epsilon^{\tau(2)}x_2,...,\epsilon^{\tau(n)}x_n),$ 
where $\tau(i)=t$ is the index of the corresponding vector space $V_t$. Furthermore, since the exponential map at the unit is a diffeomorphism, $\delta_{\epsilon}^0$ is smooth. Moreover, note that the group multiplication in the Carnot group is compatible with the dilation in the sense that for any $a,b\in G$ 
\begin{equation}\label{eeeer}
	\delta_\epsilon^0(ab)=\delta_\epsilon^0(a)\delta_\epsilon^0(b).
\end{equation}
Then, for any other point $g$ of $G$, define \[\delta_{\epsilon}^g(x)=L_g\delta^0_\epsilon L_{g^{-1}}(x).\]
We restrict the dilation to the neighborhood of $g$ by $\delta_\epsilon^g\colon U_1^g\to U_{\epsilon}^g$ when $\epsilon\geq 1$ and by $\delta_\epsilon^g\colon U_{\epsilon}^g\to U_1^g$ when $\epsilon\leq 1$. When there is no confliction, we abbreviate the above symbol as $\delta_\epsilon$. It follows easily that  
$U_{\epsilon}^g=U_{\frac{1}{\epsilon}}^g$.
Furthermore, as $L_g$, $L_{g^{-1}}$ and $\delta_{\epsilon}^0$ are smooth, $\delta_{\epsilon}^g$ is smooth and thus the differential of $\delta_{\epsilon}^g$, denoted by $d\delta_{\epsilon}^g$,  exists.

Choose an unit orthogonal basis for the vector space $V_1$ at the point $0$ and denote the left-invariant vector fields generated by the vectors in this basis by $X_1$,$X_2$,...,$X_r$. Then, for $\epsilon\in(0,1]$, $j\in\{1,\dots,r\}$ and $g\in G$, define $X_j^{g,\epsilon}:=\epsilon\cdot d\delta_{\frac1\epsilon}^g\circ X_j\circ\delta_\epsilon^g\in\Gamma(TU^g_1)$.
\begin{lemma}\label{ttt}
	For $\epsilon\in(0,1]$ and $j\in\{1,\dots,r\}$, 
	\[X_j^{g,\epsilon}=X_j.\]
\end{lemma}
\begin{proof}
	For any $s,t\in G$, 
	\[\delta^g_{\epsilon}(gsgt)=g\delta_\epsilon^0(sgt)=g\delta_\epsilon^0(s)\delta^0_{\epsilon}(gt)=\delta_\epsilon^g(gs)\delta^0_\epsilon(gt).\]
	Then, for any $v\in U^g_1$,
	\[\delta^g_{\frac1\epsilon}\circ L_{\delta^g_\epsilon v}(q)=\delta^g_{\frac1\epsilon}(\delta_\epsilon^g(v)q)=L_v\delta^o_{\frac1{\epsilon}}q,\]
	for any $q$ holds.
	Thus,
	\[X_j^{g,\epsilon}(v)=\epsilon d\delta_{\frac1\epsilon}^g(X_j(\delta^g_\epsilon(v)))=\epsilon d\delta_{\frac1\epsilon}^g( {(L_{\delta_\epsilon^gv})}_*X_j(0))=\epsilon{(L_{v})}_*(d\delta_{\frac1\epsilon}^0(X_j(0)))={(L_{v})}_*((X_j(0))),\] 
	which concludes the proof.
\end{proof}

\begin{remark}\label{rmk:on Carnot group}
	Based on Lemma~\ref{ttt}, we shall give a simplified proof of Theorem \ref{thm:Differentiability Stepanov} in the setting of Carnot groups in the appendix. In the case of subRiemannian manifolds, we have to use blow-up of vector fields and the metrics associated to them, while in the case of Carnot groups, we only need to use the horizontal vector fields and the standard metric. 
\end{remark}

\subsection{Equiregular subRiemannian manifolds}\label{subsec:geometry of equiregular SRM}

\newcommand{\Dst}{\scr{H}}

Let $M$ be a differentiable manifold of topological dimension $n$ and fix a subbundle $H\subset TM$ of rank $r$.
Define the following \emph{flag of distributions} inductively for $k\in\N$:
\[
\begin{cases}
 	\Dst^{(0)} &:= \{0\} \\
	\Dst^{(1)} &:= \Gamma(H) \\
	\Dst^{(k+1)} &:= \Dst^{(k)} + \Co^\infty(M)\text{-}\Span\left\{[X,Z] : X\in\Dst^{(1)},\ Z\in\Dst^{(k)} \right\} ,
\end{cases}
\]
where $\Gamma(H)$ is the set of all smooth sections of $H$
and for any set $S$ of vector field, $\Co^\infty(M)\text{-}\Span(S)$ is the set of linear combinations of elements of $S$ with coefficients in $\Co^\infty(M)$, which is the ring of smooth functions $M\to\R$.

By definition we have
\[
\{0\}\subset \Dst^{(1)} \subset\dots\subset\Dst^{(k)}\subset\Dst^{(k+1)}\subset\dots\subset\Vec(M) .
\]
For any point $p\in M$ we have a \emph{pointwise flag}
\[
\Dst^{(k)}_p := \{ Z(p):Z\in\Dst^{(k)}\} \subset T_pM .
\]
To such a flag we associate some functions $M\to\N\cup\{+\infty\}$:
\begin{description}
\item[ranks] 	
	$r_k(p) := \dim(\Dst_p^{(k)})$. Notice that $r=r_1\le r_2\le\dots\le n$.
\item[growth vector]
	$n_k(p) := r_k(p)-r_{k-1}(p) = \dim( \Dst_p^{(k)}/\scr \Dst_p^{(k-1)} )$.
	Notice that $\sum_{i=1}^k n_i = r_k$.
	The function $p\mapsto (n_1(p),n_2(p),\dots)\in\N^\N$ is usually called \emph{growth vector}.
\item[step]	
	$s(p) := \inf\{ k: \Dst_p^{(k)} = T_pM,\text{ i.e., }r_k(p)=n\}$.
	Notice that if $s(p)<\infty$, then
	\[
	\{0\}\subset\Dst_p^{(1)}\subset\dots\subset\Dst_p^{(s(p))}=T_pM .
	\]
\item[weight]	
	for $i\in\{1,\dots,n\}$: $w_i:=k\IFF i\in\{r_{k-1}+1,\dots,r_k\}$.
\end{description}
The subbundle $H$ is said to be \emph{equiregular} if $r_k$ (hence $n_k$ and $s$) are constant.
It is said to be \emph{bracket generating} if $s<\infty$. I

\begin{definition}[subRiemannian manifold]
 	An \emph{equiregular subRiemannian manifold} is a triple $(M,H,g)$ where $M$ is a smooth and connected manifold, $H\subset TM$ is a bracket generating and equiregular subbundle, and $g:H\times H\to[0,+\infty)$ is a smooth function whose restriction to each fiber $H_p$ is a scalar product. %]
\end{definition}

\begin{definition}[subRiemannian distance]
	An absolutely continuous curve
	$\gamma:[0,1]\to M$ is called \emph{horizontal curve (joining $\gamma(0)$ to $\gamma(1)$)} if $\gamma'(t)\in H$ for almost every $t\in[0,1]$.
	
	The \emph{length} of an horizontal curve $\gamma$ is
	\[
	l(\gamma) := \int_0^1\|\gamma'(t)\|\dd t .
	\]
 We finally define the \emph{subRiemannian distance} as
	\[
	d(p,q) := \inf\left\{
	l(\gamma):
	\text{ $\gamma$ is a horizontal curve joining $p$ to $q$}
	\right\} .
	\]
\end{definition}

An equiregular subRiemannian manifold can be endowed in a canonical way with a smooth volume $\Vol$ that is called \emph{Popp measure}.  The construction can be found in \cite{cambridge}. Given an equiregular subRiemannian manifold $M$, it is clear that its topological dimension $\dim M=\sum_kn_k$. The \emph{homogeneous dimension} $Q$ of $M$ is defined to be $Q=\sum_kkn_k$, which is typically greater than the topological dimension of $M$.  Moreover, the Hausdorff measure $\mathcal{H}^Q$ (with respect to the subRiemannian distance) is equivalent to the Lebesgue measure on $M$ (see \cite{cambridge}). 

From now on, given an equiregular subRiemannian manifold $M$, we shall always use $d$ to denote the subRiemannian distance on $M$, use $Q$ to denote the homogeneous dimension, and use $r, s$ to denote the rank and step of $M$. 

Next, we introduce the concept of the non-holonomic order of a smooth function at a point in the subRiemannian manifold as in  \cite[Definition 2.12]{m02}.
\begin{definition}[Non-holonomic order]
	Let $f:M\to \R$ be a smooth function and $o\in M$.
	The \emph{non-holonomic order of $f$ at $o$} is defined as the maximum of $k\in\N$ such that for all $i<k$ and for any choice of horizontal vector fields $X_1,\dots,X_i\in\Dst^{(1)}$ it holds
	\[
	X_1X_2\cdots X_if(o) = 0 .
	\]
\end{definition}

For any $p\in M$, there exists a local basis $X_1$,...,$X_n$ for $TM$ satisfying that there are constants $m_1$,...,$m_s\in \N$ such that $X_1$,...,$X_{m_k}$ forms a frame for $\Dst^{(k)}$. In this case, we call them an \emph{equiregular basis} and define the degree $d_i$ of $X_i$ to be the maximal index of $\Dst^{(k)}$ to which $X_i$ belongs. 

\begin{definition}
	Let $o\in M$.
	A system of coordinates $(x_1,\dots,x_n):U\to \R^n$ centered at $o$ is a \emph{system of privileged coordinates} if the function $x_i$ has non-holonomic order $w_i$.
\end{definition}

Privileged coordinates exists at all points of $M$; see \cite{b96,cambridge} for more information on this topic.

%\subsection{Ball-Box Theorem}
For each $p\in M$ and $X\in \Gamma(TM)$, we denote by $\exp_p(X)$ the value of $\gamma(1)$ at time 1 of the integral curve of the vector field $X$ starting at $p$, \ie, the solution of
\begin{align*}
\dot{\gamma}(t)=X_{\gamma(t)}\quad \text{ and }\quad \gamma(0)=p.
\end{align*}
For $p\in M$, we define the \textit{exponential coordinates} as
\begin{align*}
\Phi:\bR^n\to M
\end{align*}
\begin{align*}
\Phi(t_1,\dots,t_n):=\exp_p\big(t_1X_1+\cdots+t_nX_n\big).
\end{align*}
Notice that such map might be defined only on a neighborhood of $0\in \bR^n$.

The \textit{box} with respect to $X_1,\dots,X_n$ is defined as
\begin{align*}
\Boxs(r):=\{(t_1,\dots,t_n)\in\bR^n:|t_j|\leq r^{d_j}\},
\end{align*}
where $d_j$, $j\in \{1,\dots,n\}$, is the degree of $X_j$.

The following well-known comparison theorem is due to Mitchell, Gershkovich, Nagel-Stein-Wainger (see for instance \cite[Theorem 10.67]{cambridge}) and is called the \textit{Ball-Box Theorem} since compare the boxes $\Boxs(r)$ in $\bR^n$ with the balls $B(p,r)$ with respect to the subRiemannian distance.

\begin{theorem}[Ball-Box Theorem]\label{thm:ball-box}
Let $M$ be an equiregular subRiemannian $n$-manifold of homogeneous  dimension $Q\geq 2$ and let $\Phi$ be some exponential coordinate map from a point $p\in M$ with respect to some equiregular basis $X_1,\cdots, X_n$. Then there exist a constant $C>0$ and a radius $r_p>0$ such that
\begin{align*}
\Phi\big(\Boxs(C^{-1}r)\big)\subset B(p,r)\subset \Phi(\Boxs(Cr))
\end{align*}
for all $r\in (0,r_p)$.
\end{theorem}

Recall that a metric space $X$ is said to be  \textit{linearly locally connected} (LLC) if there exists $\theta\geq 1$ such that for each $x\in X$ and all $0<r\leq \diam X$,

(i) ($\theta$-LLC-1) every two points $a,b\in B(x,r)$ can be joined in $B(x,\theta r)$, and

(ii) ($\theta$-LLC-2) every two points $a,b\in X\backslash \closure{B}(x,r)$ can be joined in $X\backslash \closure{B}(x,\theta^{-1}r)$.

Here, by joining $a$ and $b$ in $B$ we mean that there exists a path $\gamma:[0,1]\to B$ with $\gamma(0)=a,\gamma(1)=b$.

As a particular consequence of Theorem~\ref{thm:ball-box}, we point out that an equiregular subRiemannian manifold $M$ is locally LLC and locally Ahlfors $Q$-regular, \ie for each $x\in M$, there exists a radius $r_x>0$ such that the metric space $\big(B(x,r_x),d\big)$ ($d$ is the subRiemannian distance on $M$) is LLC and Ahlfors $Q$-regular (note that the constants associated to the LLC condition and the Ahlfors regularity condition~\eqref{eq:Ahlfors regular measure} depend on the point $x$).

By the results from~\cite{j86}, an equiregular subRiemannian manifolds locally supports a $(1,1)$-Poincar\'e inequality (with the constant associated to the Poincar\'e inequality depending on the locality).

\subsection{Tangent cone of an equiregular subRiemannian manifold}\label{subsec:tangent cone}

Let $(M,d)$ be an equiregular subRiemannian manifold with horizontal distribution $H\subset TM$.
Since the results of this section are local, we can assume that $H$ is generated by $r$ smooth vector fields $X_1,\dots,X_r\in\Gamma(TM)$ that are unit and orthogonal to each other at every point.

We assume that at each point $o\in M$ a system of privileged coordinates is chosen.
Let $\delta_\epsilon^o$ be the dilations with respect to the privileged coordinate. As the privileged coordinate is local, we restrict $\delta_\epsilon^o\colon U_1^o\to U_{\epsilon}^o$ when $\epsilon\geq 1$ and $\delta_\epsilon^o\colon U_{\epsilon}^o\to U_1^o$ when $\epsilon\leq 1$. When there is no confliction, we abbreviate the above symbol as $\delta_\epsilon$. Easily, it follows that  $U_{\epsilon}=U_{\frac{1}{\epsilon}}$.

It is a well-known result that these dilations $\delta_\epsilon^o$ permit to construct the metric tangent cone of $(M,d)$ at $o\in M$. We next briefly expose the procedure; see \cite{b96} for details.

For $\epsilon\in(0,1]$ and $j\in\{1,\dots,r\}$, define $X_j^{o,\epsilon}:=\epsilon\cdot d\delta_{\frac1\epsilon}\circ X_j\circ\delta_\epsilon\in\Gamma(TU^o_1)$.
Then there are $X_j^{o,0}\in\Gamma(TU^o_1)$ such that $X_j^{o,\epsilon}\to X_j^{o,0}$ uniformly on compact sets.
Up to shrinking the set $U^o_1$, we can assume the convergence to be uniform on $U^o_1$.
Notice that, for all $\epsilon\in(0,1]$, $(\delta_\epsilon^o)_*X_j^{o,0} = \dd\delta_\frac{1}{\epsilon}^o\circ X_j^{o,0}\circ\delta_\epsilon^o = \epsilon^{-1} X_j^{o,0}$. For the precise definition of $(\delta_\epsilon^o)_*$, see \cite[Section 5.2]{b96}.

For all $\epsilon\in[0,1]$, the vector fields $X_j^{o,\epsilon}$ define a subRiemannian metric $d_\epsilon^o$ on $U^o_1$.
For $\epsilon\neq0$, the metric space $(U^o_1,d_\epsilon^o)$ is isometric to a neighborhood of $o$ in $(M,\epsilon^{-1}d)$ via $\delta^o_\epsilon$. More precisely, \begin{equation}\label{bab}
	\epsilon^{-1}d(\delta^o_\epsilon y,\delta^o_\epsilon x)
	= d^{o}_{\epsilon}(y,x),\end{equation}
for $x$,$y$ in a neighborhood of $o$.
As $\epsilon\to 0$, $d_\epsilon^o$ converge uniformly on $U^o_1\times U^o_1$ to $d_0^o$.
This implies that $(U^o_1,d_0^o)$ is isometric to a neighborhood of the origin in the tangent cone of $(M,d)$ at $o$. We will always write $d^o_0$ or just $d^o$ for the tangent distance at the point $o$. 

More can be said about the tangent cone.
Let $\g_o\subset\Gamma(TU^o_1)$ be the Lie algebra generated by the vector fields $\{X_j^{o,0}\}_{j=1}^r$.
This is a finite dimensional, nilpotent, stratified Lie algebra, whose first layer is the span of $X_1^{o,0},\dots,X_r^{o,0}$.

Recall that a Lie algebra $\g$ is \emph{stratified} of step $s$ and rank $r$ if $\g=\bigoplus_{i=1}^s V_i$ with $\dim(V_1)=r$ and $[V_1,V_i]=V_{i+1}$ for all $i$.
When we speak of a stratified Lie algebra $\g$ we mean that the stratification $V_1,\dots,V_s$ is chosen.
A map $A:\g\to\g'=\bigoplus_{i=1}^{s'}V'_i$ is a \emph{morphism of stratified Lie algebras} if it commutes with Lie brackets and $A(V_i)\subset V'_i$.

Since a stratified Lie algebra is nilpotent, the Baker-Campbell-Hausdorff formula is a finite sum and it defines a map $*:\g\times\g\to\g$ that makes $(\g,*)$ into a Lie group.
More precisely, $(\g,*)$ is the unique simply connected Lie group whose Lie algebra is $\g$.
With this identification between Lie algebra and Lie group, any Lie algebra morphism is a Lie group morphism as well.

The group $\G=(\g,*)$ becomes a \emph{Carnot group} if $V_1$ is endowed by a scalar product and $\G$ is endowed with the induced left-invariant subRiemannian metric.

In the case of $\g_o\subset\Gamma(TU^o_1)$, the first layer is $V_1=\Span\{X_1^{o,0},\dots,X_r^{o,0}\}$, and the scalar product on $V_1$ is chosen by saying that $X_1^{o,0},\dots,X_r^{o,0}$ is a orthonormal basis.

The exponential map for vector fields $\exp:\Gamma(TU^o_1)\to U^o_1$ (which is not globally defined), restricted to $\g_o$ gives an isometry between an open neighborhood of $0\in\G_o=\g_o$ onto an open neighborhood of $o$ in $(U^o_1,d^o_0)$.
Moreover, by~\cite{cr15}, we may start with some special privileged coordinates on $U^o_1$ so that they correspond to exponential coordinates of the group $\G_o$.
Therefore, if $\bar o$ is another point on another subRiemannian manifolds, and $A:\g_o\to\g_{\bar o}$ is a Lie algebra morphism, then we can see $A$ as a map $U^o_1\to U^{\bar o}_1$ that is linear in these coordinates.

\subsection{Definition of a quasiregular mapping}
Recall that a map $f\colon M\to N$ between two equiregular subRiemannian manifolds of homogeneous dimension $Q\geq 2$ is said to be a \textit{branched cover} if it is continuous, discrete, open and sense-preserving. Here a mapping $f\colon X\to Y$ between topological manifolds is \textit{discrete} if each fiber is a discrete set in $X$, \ie for all $y\in Y$, $f^{-1}(y)$ is a discrete set in $X$, and is \textit{open} if it maps open set in $X$ onto open set in $Y$, and is \textit{sense-preserving} if the local index (or degree) at each image point is positive (see \cite{hr02} for a detailed introduction). 

Comparing with \eqref{eq:def for linear dilatation}, for a branched cover $f\colon M\to N$,  we also set
\begin{align*}
	H_f'(x)=\limsup_{r\to 0} H_f'(x,r)=\limsup_{r\to 0}\frac{L_f'(x,r)}{l_f(x,r)},
\end{align*}
where
\begin{align*}
	L_f'(x,r)=\sup\big\{d(f(x),f(y)):d(x,y)= r\big\}
\end{align*}
and
\begin{align*}
	l_f(x,r)=\inf\big\{d(f(x),f(y)):y\in\bdary B(x,r)\big\}.
\end{align*}
Similarly, the \emph{weak linear dilatation} of $f$ at $x\in M$ is defined to be
	$$h_f(x)=\liminf_{r\to 0}H_f(x,r).$$
	
We have the following two natural definitions for metrically quasiregular mappings. 
\begin{definition}[Metrically quasiregular mappings I]\label{def:metric qr}
	Let $M$ and $N$ be two equiregular subRiemannian manifolds. A branched cover $f:M\to N$ is said to be \textit{metrically $H$-quasiregular of type 1} if it satisfies
	\begin{itemize}
		\item[i).] $H_f(x)<\infty$ for all $x\in M$;
		\item[ii).] $H_f(x)\leq H$ for almost every $x\in M$.
	\end{itemize}
	We say that $f:M\to N$ is \textit{metrically quasiregular of type 1} if it is metrically $H$-quasiregular of type 1 for some $1\leq H<\infty$.
\end{definition}

\begin{definition}[Metrically quasiregular mappings II]\label{def:metric qr II}
	Let $M$ and $N$ be two equiregular subRiemannian manifolds. A branched cover $f:M\to N$ is said to be \textit{metrically $H$-quasiregular of type 2} if it satisfies
	\begin{itemize}
		\item[i).] $H_f'(x)<\infty$ for all $x\in M$;
		\item[ii).] $H_f'(x)\leq H$ for almost every $x\in M$.
	\end{itemize}
	We say that $f:M\to N$ is \textit{metrically quasiregular of type 2} if it is metrically $H$-quasiregular of type 2 for some $1\leq H<\infty$.
\end{definition}

\begin{remark}\label{rmk:on discreteness}
	We would like to point out that for a continuous and open mapping, the condition $H_f'(x)<\infty$ for all $x\in M$ implies that $f$ is discrete. Indeed, since $f$ is open and non-constant, $L_f'(x,r)>0$ for all $x\in M$ and $r>0$ such that $\closure{B}(x,r)$ is compact. On the other hand, if $f$ is not discrete, then for some $y\in N$, $f^{-1}(y)$ would contain an accumulation point $x$, \ie there is a sequence $\{x_i\}_{i\in \mathbb{N}}\subset f^{-1}(y)$ such that $x_i\to x$. By continuity of $f$, $f(x_i)\to f(x)$. Since $f(x_i)=y$ for all $i$, it is necessarily true that $y=f(x)$.
	
	Consider the point $x$ and denote by $r_i=d(x,x_i)$. Then $r_i\to 0$ as $i\to \infty$. Notice that $L_f'(x,r_i)>0$ and that
	\begin{align*}
		l_f(x,r_i)=\inf\big\{d(f(x),f(y)):y\in\bdary B(x,r_i)\big\}=0.
	\end{align*}
	We thus obtain
	\begin{align*}
		H_f'(x)=\limsup_{r\to 0}\frac{L_f'(x,r)}{l_f(x,r)}=\infty,
	\end{align*}
	contradicting the fact that $H_f'<\infty$ everywhere.
\end{remark}

It is clear from the definition that $f$ is metrically $H$-quasiregular of type 2 whenever it is metrically $H$-quasiregular of type 1. Indeed, the converse is also true and we will prove in Proposition~\ref{prop:type 1=type 2} below that a branched cover $f:M\to N$ is metrically $H$-quasiregular of type 1 if and only if it is metrically $H$-quasiregular of type 2.

We next introduce the so-called weakly metrically quasiregular mappings.
\begin{definition}[Weakly metrically quasiregular mappings]\label{def:weak metric qr}
	A branched cover $f:M\to N$ is said to be \textit{weakly metrically $H$-quasiregular} if it is constant or if it satisfies
	\begin{itemize}
		\item[i).] $h_f(x)<\infty$ for all $x\in M$;
		\item[ii).] $h_f(x)\leq H$ for almost every $x\in M$.
	\end{itemize}
	We say that $f:M\to N$ is \textit{weakly metrically quasiregular} if it is weakly metrically $H$-quasiregular for some $1\leq H<\infty$.
\end{definition}

We list some non-trivial examples of quasiregular mappings in the subRiemannian setting.
\begin{example}\label{exam:hh}
	i).(Example 6.23,~\cite{hh97}) Let $\mathbb{H}_1$ be the first Heisenberg group and let $X,Y,T$ be the associated vector fields. The mapping $f:\mathbb{H}_1\to \mathbb{H}_1$ defined as $(r,\varphi,t)\mapsto (r/2,2\varphi,t)$ in cylindrical coordinates is a quasiregular mapping with non-empty branch set. This mapping is a counterpart of the winding mapping in the Euclidean setting~\cite{r93} and the branch set of $f$ is the $t$-axis.
	
	ii).(Theorem 3.5,~\cite{flp14}) The sphere $\mathbb{S}^{2n+1}\subset\mathbb{C}^{n+1}$ admits a natural subRiemannian structure. Namely, one obtains the horizontal subbundle by taking a maximal complex subspace of $T\mathbb{S}^{2n+1}:H\mathbb{S}^{2n+1}=T\mathbb{S}^{2n+1}\cap iT\mathbb{S}^{2n+1}$. The Euclidean inner product on $\mathbb{C}^{n+1}$ then restricts to $H\mathbb{S}^{2n+1}$ as a subRiemannian metric tensor $g$ and the corresponding norm $|\cdot|$.
	
	For $a\in \mathbb{Z}$, consider the multi-twist mapping of $\mathbb{S}^{2n+1}$ given by
	\begin{align*}
		F_a(re^{i\theta_1},\cdots,r_{n+1}e^{i\theta_{n+1}})=(r_1e^{ia\theta_1},\cdots,r_{n+1}e^{ia\theta_{n+1}}).
	\end{align*}
	Then $F_a:\mathbb{S}^{2n+1}\to \mathbb{S}^{2n+1}$ is (metrically) $|a|$-quasiregular.
	
	iii). (Lemma 3.11,~\cite{flp14}) Let $p>1$ be an integer and $q_1,\dots,q_{n+1}\in \mathbb{N}$ relatively prime to $p$. Set $\textbf{q}=(q_1,\dots,q_{n+1})$, and define
	\begin{align*}
		R_{p,\textbf{q}}(z_1,z_2,\dots,z_{n+1})=(e^{2\pi iq_1/p},e^{2\pi iq_2/p}z_2,\cdots,e^{2\pi iq_{n+1}/p}z_{n+1}).
	\end{align*}
	The associated lens space $L_{p,\textbf{q}}$ is set to be the quotient space $\mathbb{S}^{2n+1}/\langle R_{p,\mathbf{q}}\rangle$. Moreover, $L_{p,\mathbf{q}}$ admits a natural subRiemannian structure; see for instance~\cite[Proposition 3.1]{flp14}.
	
	For each $a\in p\mathbb{Z}$ for some positive integer $p$, the multi-twist mapping $F_a$ from ii) induces a well-defined mapping on the lens space, namely
	\begin{align*}
		f_a:L_{p,\mathbf{q}}\to \mathbb{S}^{2n+1},\quad f_a({[z]}):=F_a(z)\ \text{ for } z\in \mathbb{S}^{2n+1}.
	\end{align*}
	Then $f_a$ is a quasiregular mapping. Moreover, if we denote by $\pi:\mathbb{S}^{2n+1}\to L_{p,\mathbf{q}}$ the usual projection, then the multi-twist mapping $\pi\circ f_a:L_{p,\mathbf{q}}\to L_{p,\mathbf{q}}$ of the lens spaces is quasiregular as well.
\end{example}

If $f:M\to N$ is a continuous map, $y\in N$ and $A\subset M$, we use the notation
\begin{equation}\label{def:multiplicity function}
	N(y,f,A)=\text{card}\{f^{-1}(y)\cap A\}
\end{equation}
for the multiplicity function.

\begin{definition}[Analytically quasiregular mappings]\label{def:analytic def for qr}
	A branched cover $f:M\to N$ is said to be \textit{analytically $K$-quasiregular} if $f\in N^{1,Q}_{loc}(M,N)$ and
	\begin{align*}
		g_f^Q(x)\leq KJ_f(x)
	\end{align*}
	for almost every $x\in M$.
\end{definition}

The ``Jacobian" above is given by
\begin{equation}\label{eq:analytic jacobian}
	J_f(x)=\frac{d\big(f^*\Vol_N\big)(x)}{d\Vol_M(x)},
\end{equation}
where the pull-back $f^*\Vol_N$ is defined as
\begin{equation}\label{eq:definition of pull-back}
	f^*\Vol_N(A)=\int_Y N(y,f,A)d\Vol_N(y).
\end{equation}
Note that a branched cover between manifolds necessarily has locally bounded multiplicity. It is not difficult to see that $f^*\Vol_N$ is a Radon measure.

We point out an useful observation by the third named author that the Jacobian $J_f(x)$ from~\eqref{eq:analytic jacobian} can be alternately described by
\begin{equation}\label{eq:another description of Jacobian}
	J_f(z)=\lim_{r\to 0}\frac{\Vol_N\big(f(B(z,r))\big)}{\Vol_M(B(z,r))}
\end{equation}
for almost every $z\in M$; see~\cite[Section 4.3]{g14} for a simple proof of this fact.

%%%%%%%%%%%%%%%%%%%%%%%%%%%%%%%%%%%%%%%%%%%%%%%%%%%%%%%%%%%%%%%
%%%%%%%%%%%%%%%%%%%%%%%%%%%%%%%%%%%%%%%%%%%%%%%%%%%%%%%%%%%%%%%
\section{Analytic foundations of quasiregular mappings}\label{sec:QR mapping}

The goal of this section is to prove the following theorem, which establishes the analytic foundations of metrically quasiregular mappings between equiregular subRiemannian manifolds. 

\begin{theorem}\label{thm:properties}
	Let $f\colon M\to N$ be a weakly metrically quasiregular mapping between two equiregular subRiemannian manifolds of homogeneous dimension $Q\geq 2$. Then
	
	1). $f$ satisfies Condition $N$, \ie $\Vol_N(f(E))=0$ if $\Vol_M(E)=0$;
	
	2). The area formula holds, namely, for all measurable function $h:N\to [0,\infty]$ and every measurable set $A\subset M$,
	\begin{align}\label{eq:area formula}
		\int_{A}h(f(x))J_f(x)d\Vol_M(x)=\int_{N}h(y)N(y,f,A)d\Vol_N(y),
	\end{align}
	where $N(y,f,A)=\card \big(f^{-1}(y)\cap A\big)$ is the multiplicity function of $f$ on $A$;
	
	3). $f$ satisfies Condition $N^{-1}$, \ie $\Vol_M(f^{-1}(E))=0$ if $\Vol_N(E)=0$, and is $P$-differentiable almost everywhere and the Jacobian $J_f> 0$ almost everywhere in $M$;
	
	4). $\Lip f$ is the minimal $Q$-weak upper gradient of $f$;
	
	5). $\Vol_M(\mathcal{B}_f)=0$, where $\mathcal{B}_f$ is the branch set of $f$, i.e. the set of all $x\in M$ such that $f$ fails to be a local homeomorphism at $x$.
\end{theorem}

Note in particular that Theorem \ref{thm:branch set} follows immediately from Theorem \ref{thm:properties} 1) and 5).

\subsection{Analytic properties of quasiregular mappings}
Recall that for a continuous mapping $f:X\to Y$ between two metric spaces, the \textit{upper Lipschitz constant} $\Lip f(x)$ of $f$ at $x\in X$ is defined as
\begin{align*}
	\Lip f(x)=\limsup_{r\to 0}\sup_{y\in B(x,r)}\frac{d(f(x),f(y))}{r}.
\end{align*}
Similarly, the \textit{lower Lipschitz constant} $\lip f(x)$ of $f$ at $x\in X$ is defined as
\begin{align*}
	\lip f(x)=\liminf_{r\to 0}\sup_{y\in B(x,r)}\frac{d(f(x),f(y))}{r}.
\end{align*}

The following simple lemma follows from the fact that an equiregular subRiemannian manifold $N$ is locally linearly locally connected.
\begin{lemma}\label{lemma:consequence of ball-box}
	Let $f:M\to N$ be a continuous, discrete and open mapping. Then for each $x\in M$, there exist a constant $C\geq 1$ and a radius $r_x>0$ such that for each $0<r\leq r_x$,
	\begin{align*}
		B\Big(f(x),\frac{l_f(x,r)}{C}\Big)\subset f\Big(B(x,r)\Big)\subset B\Big(f(x),CL_f'(x,r)\Big).
	\end{align*}
\end{lemma}
\begin{proof}
	We only prove the first inclusion since the proof of the second one is similar. Fix a point $x\in M$. Since $M$ is locally compact, there exists $r_x>0$ such that $\closure{B}(x,r_x)$ is compact. Since $N$ is equiregular, it is locally LLC and thus we may further assume that there exist a constant $C\geq 1$ and a radius $r_x>0$ such that each two points $a,b\in B\Big(f(x),\frac{l_f(x,r)}{C}\Big)$ can be joined in $B(f(x),l_f(x,r))$ whenever $0<r<r_x$. 
	
	For $0<r<r_x$, we claim that
	$$B\Big(f(x),\frac{l_f(x,r)}{C}\Big)\subset f\Big(B(x,r)\Big).$$
	%For if not,
	Assume the opposite, then there exists a point $b\in B\Big(f(x),\frac{l_f(x,r)}{C}\Big)\backslash f\Big(B(x,r)\Big)$. By the preceding assumption, we may find a path $\gamma$ in $B\big(f(x),l_f(x,r)\big)$ that joins $f(x)$ and $b$. Now $f\big(B(x,r)\big)$ is open, $\gamma\cap f(B(x,r))\neq \emptyset$ and $\gamma\cap \big(N\backslash f(B(x,r))\big)=\emptyset$, and so it follows that $\gamma\cap \partial f(B(x,r))\neq \emptyset$. Since $f$ is open and $\closure{B}(x,r)$ is compact, $\partial f(B(x,r))\subset f(\partial B(x,r))$, and so $\gamma\cap f\big(\partial B(x,r)\big)\neq \emptyset$. Choose $y\in \partial B(x,r)$ so that $f(y)\in \gamma\subset B(f(x),l_f(x,r))$. Then $d(f(x),f(y))<l_f(x,r)$ and $d(x,y)=r$, which contradicts the definition of $l_f(x,r)$ and so our claim holds.
\end{proof}

\begin{remark}\label{rmk:on equivalence of type 1 and type 2}
	i). Note that a direct consequence of Lemma~\ref{lemma:consequence of ball-box} is that  if $f:M\to N$ is a metrically quasiregular mapping of type 2, then it is a locally metrically quasiregular mapping of type 1.
	
	Indeed, by Lemma~\ref{lemma:consequence of ball-box}, for each $x\in M$, there exist a constant $C>0$ and a radius $r_x>0$ such that  $f(B(x,r))\subset B(f(x),CL_f'(x,r))$ for all $0<r<r_x$. This implies that for $0<r<r_x$, $L_f(x,r)\leq CL_f'(x,r)$ and so $$H_f(x)\leq CH_f'(x).$$
	
	ii). An equiregular subRiemannian manifold is locally geodesic and so in particular, it is 1-LLC-1 and thus the first inclusion in Lemma~\ref{lemma:consequence of ball-box} holds with $C=1$.
\end{remark}

The following lemma is a special case of the more general result from~\cite[Theorem A]{gw16}. For the convenience of the readers, we also include a sketch here.
\begin{lemma}\label{lemma:condition N}
	Let $f:M\to N$ be a weakly metrically quasiregular mapping. Then $f\in N^{1,Q}_{loc}(M,N)$ and satisfies Condition $N$.
\end{lemma}
\begin{proof}
	Note first that by~\cite[Lemma 3.3]{g14}, a continuous open mapping $f:M\to N$ is locally pseudomonotone and that by~\cite[Theorem 7.2]{hkst01}, a locally pseudomonotone map in $N^{1,Q}_{loc}(M,N)$ satisfies Condition $N$. Thus, it suffices to establish the Sobolev regularity.
	
	By~\cite[Proof of Theorem 1.1]{w14}, $f$ satisfies Condition $N$ on $Q$-almost every rectifiable curve. Note that in~\cite[Proof of Theorem 1.1]{w14}, $f$ is assumed to be a homeomorphism. However, the proof there also works for continuous, discrete and open mappings since the homeomorphism assumption was only used to deduce the local $L^1$-integrablity of the Jacobian.
	
	On such a curve, either $\int_\gamma \lip f ds=\infty$, in which case the upper gradient inequality~\eqref{ugdefeq} trivially holds on $\gamma$, or we may apply~\cite[Lemma 3.6]{z07} (more precisely the claim in the proof of ~\cite[Lemma 3.9]{z07}) to infer, as in the final proof of~\cite[Lemma 3.9]{z07}, that $\lip f$ satisfies~\eqref{ugdefeq} on $\gamma$. This means that $\lip f$ is a $Q$-weak upper gradient of $f$. To establish the desired Sobolev regularity for $f$, we only need to show that $\lip f\in L^Q_{loc}(M)$.

	By Remark~\ref{rmk:on equivalence of type 1 and type 2} ii) and Lemma~\ref{lemma:consequence of ball-box}, for each $B(x,r)\subset\subset M$, $B(x,l_f(x,r))\subset f(B(x,r))$ and so it follows from Lemma~\ref{lemma:consequence of ball-box} that
	\begin{align*}
		(\lip f)^Q(x)&\leq \liminf_{r\to 0}\frac{L_f(x,r)^Q}{r^Q}\leq c\liminf_{r\to 0}\frac{\Vol_N\big(B(f(x),l_f(x))\big)}{\Vol_M\big(B(x,r)\big)}\\
		&\leq c\liminf_{r\to 0}\frac{\Vol_N\big(f(B(x,r))\big)}{\Vol_M\big(B(x,r)\big)},
	\end{align*}
	where the constant $c$ depends on $x$, $Q$ and $K$ but it is bounded on compact subset of $M$ (due to the Ball-Box Theorem). On the other hand, for almost every $x\in M$,
	\begin{align*}
		\liminf_{r\to 0}\frac{\Vol_N\big(f(B(x,r))\big)}{\Vol_M\big(B(x,r)\big)}=J_f(x),
	\end{align*}
	from which we conclude that  $\lip f\in L^Q_{loc}(M)$. This completes our proof.
\end{proof}

As a corollary of the well-known area formula for Lipschtiz mappings~\cite[Theorem 1]{m11} and Lemma~\ref{lemma:condition N}, we thus obtain the area formula for quasiregular mapping between equiregular subRiemannian manifolds.

\begin{proposition}\label{prop:area formula Lip}
	Let $f:M\to N$ be a weakly metrically quasiregular mapping. Then the area formula
	\begin{equation}\label{eq:area formula Lip}
		\int_A u\circ f(x)J_f(x)d\Vol_M(x)=\int_N u(y)N(y,f,A)d\Vol_N(y)
	\end{equation}
	holds for each nonnegative measurable function $u$ on $N$ and every measurable set $A\subset M$.
\end{proposition}

\begin{proposition}\label{prop:Condition N inverse and positiveness of Jacobian}
	Let $f:M\to N$ be a weak metrically quasiregular mapping. Then $J_f(x)>0$ for almost every $x$ in $\Omega$ and $f$ satisfies Condition $N^{-1}$.
\end{proposition}
\begin{proof}
	Note that by the proof of Lemma~\ref{lemma:condition N}, for any $x\in M$, there exists a neighborhood $U_x$ of $x$ such that
	$$\lip f(p)^Q\leq CJ_f(p)$$
almost everywhere  in $U_x$, where $C$ allows to depend on $x$ as well. On the other hand, since by the proof of Lemma~\ref{lemma:condition N} $\lip f$ is an upper gradient of $f$, the preceding inequality implies that $f$ is locally analytically $C$-quasiregular according to the definition of~\cite{g14}. Since the issue is local, we may assume that $U_x$ and $f(U_x)$ each supports a $(1,1)$-Poincar\'e inequality. The result now follows from~\cite[Propositi\blue{on} 4.18]{g14}.
\end{proof}

\begin{remark}\label{rmk:on analytic def}
	It was asked in~\cite[Question 2.19]{flp14} \textit{whether condition i) in Definition~\ref{def:metric qr} is equivalent to the assumption that $H_f$ is locally bounded in $M$}. This is indeed the case, since locally, equiregular subRiemannian manifolds have bounded geometry. The proof of Theorem 7.1 from~\cite{or09} can be adapted to our setting by first noticing that, locally, metrically quasiregular mappings are analytically quasiregular (as seen in the proof of Proposition~\ref{prop:Condition N inverse and positiveness of Jacobian}) and hence  ~\cite[Theorem 6.2]{or09} holds in this setting. Moreover, by~\cite[Lemma 3.3]{g14}, the necessary topological Lemma 4.2 from~\cite{or09} holds in our setting as well. With these facts at hand, the proof of Theorem 7.1 from~\cite{or09} works directly in our situation; see also~\cite{gw16} for a proof of this fact in more general metric spaces.
\end{remark}

Recall that if $f:M\to N$ is $P$-differentiable at $x_0\in M$, then the Jacobian of the $P$-differential $Df(x_0)$ is defined as
\begin{align}\label{eq:Jacobian via P differential}
	J_{Df(x_0)}(0)=\frac{\Vol\big(Df(x_0)(B(0,1))\big)}{\Vol(B(0,1))}.
\end{align}
By \cite[Proposition 3.40]{clo14}, the Jacobian of a quasiregular mapping coincides with the Jacobian of the $P$-differential almost everywhere. 

\begin{proposition}[\cite{clo14}]\label{prop:Jacobian=det Df}
	Let $f:M\to N$ be a mapping that is $P$-differentiable at $x_0\in M$. Then
	\begin{align*}
		J_{Df(x_0)}(0)=J_f(x_0).
	\end{align*}
\end{proposition}

\subsection{Differentiability of quasiregular mappings}

In this section, we show that weak metrically quasiregular mappings are $P$-differentiable almost everywhere in $M$.

\begin{proposition}\label{prop:on differentiability}
	Let $f:M\to N$ be an open mapping between two equiregular subRiemannian manifolds of homogeneous dimension $Q\geq 2$. If $f\in N^{1,Q}_{loc}(M,N)$, then $f$ is $P$-differentiable almost everywhere in $M$.
\end{proposition}
\begin{proof}
	This is a simple consequence of Theorem~\ref{thm:Differentiability Stepanov} and~\cite[Theorem 7.2]{hkst01}. Indeed, the openness of $f$ implies that $f$ is locally \textit{pseudomonotone} (see for instance~\cite[Lemma 3.3]{g14}). A locally pseudomonotone mapping in $N^{1,Q}_{loc}(M,N)$ satisfies the so-called \textit{Rado--Reichelderfer condition}, namely, for each $x\in M$, there exists a radius $r_x>0$ such that
	\begin{align*}
		\diam f(B(x,r_x))^Q\leq C\int_{B(x,r_x)}g_f(x)^Qd\Vol_M(x).
	\end{align*}
	As a consequence of the Lebesgue differentiation theorem, we obtain that $\Lip f(x)<\infty$ for almost every $x\in M$ and thus the claim follows from Theorem~\ref{thm:Differentiability Stepanov}.
\end{proof}

\begin{theorem}\label{thm:differentiability of QR map}
	Let $f:M\to N$ be a weak metrically quasiregular mapping. Then $f$ is $P$-differentiable almost everywhere in $M$. Moreover, for almost every $x_0\in M$, the $P$-differential $Df(x_0)$ is a Carnot group isomorphism that commutes with group dilations.
\end{theorem}
\begin{proof}
	The almost every differentiability follows immediately from Lemma~\ref{lemma:condition N} and Proposition~\ref{prop:on differentiability}. Regarding the second claim, notice that by Proposition~\ref{prop:Condition N inverse and positiveness of Jacobian}, the Jacobian of $f$ is positive almost everywhere in $M$. On the other hand, by Proposition~\ref{prop:Jacobian=det Df}, for almost every $x_0\in M$, the Jacobian $J_f(x_0)$ coincides with the determinant of the $P$-differential $Df(x_0)$ and so $Df(x_0)$ is a group isomorphism.
\end{proof}

\subsection{The minimal upper gradient and Jacobian}
Let $f:M\to N$ be a mapping. If $f$ is $P$-differentiable at $x_0\in M$, then the maximal norm of the differential $Df(x_0)$ is defined as
\begin{align}\label{eq:def maximal norm}
	\|Df(x_0)\|:=\max\Big\{d_0^{f(x_0)}\big(0,Df(x_0)(v)\big):d_0^{x_0}(0,v)\leq 1\Big\}.
\end{align}
Similarly, the minimal norm of the differential $Df(x_0)$ is defined as
\begin{align}\label{eq:def for minimal norm}
	\|Df(x_0)\|_s:=\min\Big\{d_0^{f(x_0)}\big(0,Df(x_0)(v)\big):d_0^{x_0}(0,v)\geq 1\Big\}.
\end{align}
Note that by~\cite{cc06}, we have 
\begin{equation*}
	\|Df(x_0)\|=\max\Big\{|Df(x_0)v|:v\in V_1 \text{ and } |v|=1 \Big\}=\limsup_{r\to 0}\frac{L_f(x_0,r)}{r}
\end{equation*}
and that
\begin{equation*}
	\|Df(x_0)\|_s=\min\Big\{|Df(x_0)v|:v\in V_1 \text{ and } |v|=1 \Big\}=\liminf_{r\to 0}\frac{l_f(x_0,r)}{r}.
\end{equation*}

\begin{lemma}\label{lemma:simple}
	If $f:M\to N$ is a weak metrically $H$-quasiregular mapping,  then
	\begin{align}\label{eq:bound on norm in tangent cone}
		\frac{\|Df(x_0)\|}{\|Df(x_0)\|_s}\leq H
	\end{align}
	for almost every $x_0\in M$.
\end{lemma}
\begin{proof}
	This follows from~\cite[Proofs of Lemma~5.15 and Theorem 5.16]{flp14}. The only difference one has to notice is that the definition of a metrically $H$-quasiregular mapping in~\cite{flp14} is stronger than the one used in this paper, namely the branch set and its image are assumed to have zero measure. However, the proofs of Lemma~5.15 and Theorem 5.16 there only use the fact that at almost every $x_0\in M$, $f$ is $P$-differentiable and the $P$-differential $Df(x_0)$ is a Carnot group isomorphism, which is provided by Theorem~\ref{thm:differentiability of QR map}.
\end{proof}

It is a well-known fact that if $f:\bR^n\to \bR^n$ is differentiable at $x_0$, then
\begin{align*}
	\lip f(x_0)=\Lip f(x_0)=\|Df(x_0)\|.
\end{align*}
The next lemma can be viewed as a natural generalization of this fact to the subRiemannian setting.
\begin{lemma}\label{lemma:Lip=norm of differential}
	Let $f:M\to N$ be a mapping which is $P$-differentiable at $x_0\in M$. Then
	\begin{align}\label{Lip=norm of differential}
		\lip f(x_0)=\Lip f(x_0)=\|Df(x_0)\|.
	\end{align}
\end{lemma}
\begin{proof}
	Since the issue is local, we may assume that $M$, $N$ are open subsets of $\bR^n$ and $x_0=0=f(x_0)$. The $P$-differentiability of $f$ at $0$ implies that
	\begin{align*}
		\lim_{y\to 0}\frac{d(f(y),Df(0)y)}{d(0,y)}=0.
	\end{align*}
	Note that by the distance comparison estimate~\cite[Theorem 2.2]{j14}, we have
	\begin{align*}
		\sup_{y\in B(0,r)}\frac{d(f(y),f(0))}{r}&\leq \sup_{y\in B(0,r)}\frac{d(Df(0)y,f(0))}{r}+\sup_{y\in B(0,r)}\frac{d(f(y),Df(0)y)}{r}\\
		&= \sup_{y\in B(0,r)}\Big(\frac{d(Df(0)y,f(0))}{d^0(Df(0)y,0)}\cdot \frac{d^0(Df(0)y,0)}{r}\Big)+\sup_{y\in B(0,r)}\frac{o\big(d(0,y)\big)}{r}\\
		&\leq \sup_{y\in B^0\big(0,(1+O(r))r\big)}\big(1+O(r)\big)\frac{d^0(Df(0)y,0)}{r}+\frac{o(r)}{r}\\
		&=\sup_{y\in B^0\big(0,(1+O(r))r\big)}\big(1+O(r)\big)\frac{\|Df(0)\|d^0(0,y)}{r}+\frac{o(r)}{r},
	\end{align*}
	from which we conclude that $\Lip f(0)\leq \|Df(0)\|$.
	
	On the other hand, let $0<\varepsilon<1$ and a similar computation as above implies that
	\begin{align*}
		\sup_{y\in B(0,r)}\frac{d(f(y),f(0))}{r}\geq \sup_{y\in B^0\big(0,(1-O(r))r\big)}\big(1-O(r)\big)\frac{\bar{d}(Df(0)y,0)}{r}-\sup_{y\in B(0,r)}\frac{d(f(y),Df(0)y)}{r}.
	\end{align*}
	We may choose $y_0\in \bdary B\big(0,(1-\varepsilon)(1-O(r))r\big)$ realizing $\|Df(0)\|$ to deduce that
	\begin{align*}
		\sup_{y\in B^0\big(0,(1-O(r))r\big)}\frac{d(f(y),f(0))}{r}&\geq \big(1-O(r)\big)\frac{\|Df(0)\|d(y_0,0)}{r}-\sup_{y\in B(0,r)}\frac{d(f(y),Df(0)y)}{r}\\
		&=\big(1-O(r)\big)^2(1-\varepsilon)\|Df(0)\|-\sup_{y\in B(0,r)}\frac{o\big(d(0,y)\big)}{r}.
	\end{align*}
	This leads to $(1-\varepsilon)\|Df(0)\|\leq \lip f(0)$. Letting $\varepsilon\to 0$ gives
	\begin{align*}
		\Lip f(0)\leq \|Df(0)\|\leq \lip f(0)\leq \Lip f(0).
	\end{align*}
\end{proof}

%\begin{lemma}\label{lemma:linear algebra}
%Let $L:V\to W$ be a linear map between two Hilbert spaces. Then
%\begin{align*}
%\|L\|=\|L^*\|.
%\end{align*}
%\end{lemma}
%\begin{proof}
%By definition,
%\begin{align*}
%\|L^*\|&=\sup\big\{|L^*w^*|:w^*\in W^*,|w^*|=1\big\}\\
%&=\sup\big\{\langle L^*w^*,v\rangle:v\in V,|v|=1,w^*\in W^*,|w^*|=1\big\}\\
%&\leq \sup\big\{\langle w^*,Lv\rangle:v\in V,|v|=1,w^*\in W^*,|w^*|=1\big\}\\
%&\leq \sup\big\{|w^*||Lv|:v\in V,|v|=1, w^*\in W^*,|w^*|=1\big\}\\
%&=\sup\big\{|Lv|:v\in V,|v|=1\big\}=\|L\|.
%\end{align*}
%Notice that $\|L\|=\|L^{**}\|$, we may repeat the preceding argument to deduce that $\|L\|=\|L^{**}\|\leq \|L^*\|$. This completes the proof of the lemma.
%\end{proof}

We are ready to show that $\Lip f$ is the minimal $Q$-weak upper gradient of a quasiregular mapping $f$.
\begin{proposition}\label{prop:minimal upper gradient}
	Let $f:M\to N$ be a weak metrically quasiregular mapping. Then $\Lip f$ is the minimal $Q$-weak upper gradient of $f$, \ie
	\begin{align}\label{eq:minimal upper gradient}
		g_f=\Lip f.
	\end{align}
\end{proposition}

\begin{proof}
	Since the issue is local, we may assume that $M$ and $N$ are open subsets of $\bR^n$. We will follow the arguments used in~\cite[Proof of Theorem 5.2]{w12proc}. For a function $F\colon N\to \R$, we denote by $\nabla_H F=(X_1F,\cdots,X_rF)$ the horizontal gradient of the function $F$ (in~\cite[Section 11.1]{hk00}, the authors use $XF$ to represent the horizontal gradient of $F$).
	
	Note that $\Lip f$ is always a $Q$-weak upper gradient of $f$ (see for instance the proof of Lemma~\ref{lemma:condition N}) and we need to show the minimality, namely if $g$ is any other $Q$-weak upper gradient of $f$, then $g(x)\geq \Lip f(x)$ for almost every $x\in M$.
	
	To this end, let $\mathcal{D}$ be the set of $C^1$-functions $F:N\to \R$ such that $|\nabla_H F|=1$ on $N$ and $\mathcal{D}_0\subset \mathcal{D}$ a countable dense subset in the sense that for each $y\in N$, $\{\nabla_HF(y):F\in \mathcal{D}_0\}$ is dense in the horizontal unit sphere $\{v\in V_1: |v|=1\}$. Notice that if $g$ is a $Q$-weak upper gradient of $f$, then for each $F\in \mathcal{D}_0$ it is also a $Q$-weak upper gradient of $F\circ f$.
	
	Since $F\circ f:M\to \bR$ is a real-valued function, by the result of Cheeger~\cite{c99}, $\Lip (F\circ f)$ is the minimal $Q$-weak upper gradient of $F\circ f$. From this, we infer that
	$$g\geq \sup_{F\in \mathcal{D}_0}\Lip(F\circ f)$$
	almost everywhere in $M$. On the other hand, by~\cite[Theorem 11.7]{hk00}
	\begin{align*}
		\Lip(F\circ f)(x)&=|\nabla_H(F\circ f)(x)|=|Df(x)^*\nabla_HF(f(x))|
	\end{align*}
	for almost every $x\in M$, where $Df(p)^*:T_{f(p)}N\to T_pM$ is the formal adjoint of $Df(p)$. Note that if $x$ is a differentiable point of $f$, then the set $\{\nabla_HF(f(x)): F\in\mathcal{D}_0\}$ is dense in the horizontal unit sphere $\{v\in V_1: |v|=1\}$, and so at such a point $x\in M$,
	\begin{align*}
		\sup_{F\in \mathcal{D}_0}\Lip(F\circ f)(x)&=\sup_{F\in \mathcal{D}_0}|Df(x)^*\nabla_HF(f(x))|=\|Df(x)^*\|=\|Df(x)\|,
	\end{align*}
	where in the last equality we have used the standard fact that for a linear map $L:V\to W$ between two Hilbert spaces $V$ and $W$,  $\|L\|=\|L^*\|$. Therefore, it follows from the above estimate and Lemma~\ref{lemma:Lip=norm of differential} that
	\begin{align*}
		g\geq \sup_{F\in \mathcal{D}_0}\Lip(F\circ f)=\|Df\|=\Lip f
	\end{align*}
	almost everywhere in $M$. This completes the proof of Proposition~\ref{prop:minimal upper gradient}.
\end{proof}

\begin{remark}\label{rmk:on minimal upper gradient}
	It is already well-known that for a quasiregular mapping $f:M\to N$ between two equiregular subRiemannian manifolds of homogeneous dimension $Q\geq 2$, $\Lip f$ is always locally comparable with the minimal $Q$-upper gradient of $f$ (with a constant depending on the locality) (see for instance~\cite[Section 5]{w12proc}). However, the exact coincidence of these two functions is a highly non-trivial fact and indeed plays an important role in the identification of the different notions of 1-quasiconformality (or conformality) in the interesting work~\cite{clo14} and different $K$-quasiconformality in \cite{gl16}.
\end{remark}

\subsection{Equivalence of the two definitions of metric quasiregularity}
In this section, we show that the two definitions of metric quasiregularity are indeed equivalent. 
\begin{proposition}\label{prop:type 1=type 2}
	Let $f:M\to N$ be a branched cover. Then $f$ is metrically $H$-quasiregular of type 1 if and only if it is metrically $H$-quasiregular of type 2.
\end{proposition}
\begin{proof}
	Note first that by Remark~\ref{rmk:on equivalence of type 1 and type 2}, for each $x\in M$, there exists a constant $C>0$ such that, $H_f'(x)\leq H_f(x)\leq CH_f'(x)$. In particular, if $H_f'(x)<\infty$, then $H_f(x)<\infty$. We are thus left to show that $H_f(x)\leq H_f'(x)$ almost everywhere in $M$.
	
	By Theorem~\ref{thm:differentiability of QR map} and Proposition~\ref{prop:Condition N inverse and positiveness of Jacobian}, $f$ is $P$-differentiable almost everywhere and $J_f(x)>0$ for almost every such $x$ in $M$. Let $x$ be a differentiable point of $f$ with $J_f(x)>0$. Since the issue is local, we may assume that $M$ and $N$ are open subsets of $\bR^n$ and $x=0=f(x)$. By the $P$-differentiability of $f$ at $0$ and the distance comparison estimate~\cite[Theorem 2.2]{j14}, for $r=d(0,y)$ small enough,
	\begin{align*}
		L_f(0,r)&=\sup\big\{d(f(0),f(y)):d(0,y)\leq r\big\}\\
		&=\sup\big\{d(f(0),Df(0)y)+\big(d(f(0),f(y))-d(f(0),Df(0)y)\big):d(0,y)\leq r\big\}\\
		&\leq \sup\big\{d(f(0),Df(0)y):d(0,y)\leq r\big\}+o(r)\\
		&\leq \sup\Big\{\frac{d(f(0),Df(0)y)}{d^0(0,Df(0)y)}d^0(0,Df(0)y):d^0(0,y)\leq (1+O(r))r\Big\}+o(r)\\
		&\leq \big(1+O(r)\big)\sup\Big\{d^0(0,Df(0)y):d^0(0,y)=(1+O(r))r\Big\}+o(r)\\
		&\leq \big(1+O(r)\big)^2L_f'(0,r)+2o(r).
	\end{align*}
	Similarly, we deduce that
	\begin{align*}
		l_f(0,r)&=\inf\big\{d(f(0),f(y)):d(0,y)= r\big\}\\
		&\geq \big(1-O(r)\big)\inf\big\{d^0(0,Df(0)y):d(0,y)=(1-O(r))r\big\}-2o(r)\\
		&\geq \big(1-O(r)\big)^2\|Df(0)\|_sr-2o(r).
	\end{align*}
	As a consequence of the above estimates,
	\begin{align*}
		H_f(x)&=\limsup_{r\to 0}H_f(x,r)=\limsup_{r\to 0}\frac{L_f(x,r)}{l_f(x,r)}\\
		&\leq \limsup_{r\to 0}\frac{\big(1+O(r)\big)^2L_f'(x,r)}{l_f(x,r)}+\limsup_{r\to 0}\frac{o(r)}{l_f(x,r)}\\
		&\leq \limsup_{r\to 0}\big(1+O(r)\big)^2H_f'(x,r)+\limsup_{r\to 0}\frac{o(r)}{\big(1-O(r)\big)^2\|Df(x)\|_sr-o(r)}\\
		&=H_f'(x).
	\end{align*}
	This shows the coincidence of these two definitions as desired.
\end{proof}

\begin{remark}\label{rmk:on equivalence of weak metric quasiregular mappings}
	In principle, one could also define two types of weak metrically quasiregular mappings as the metrically quasiregular case. As an immediate consequence of the proof of Proposition~\ref{prop:type 1=type 2}, we conclude that the two types of definitions of weak metrically quasiregularity will be equivalent.
\end{remark}

\subsection{Size of the branch set}\label{subsec:branch set}

In this section, we prove that the branch set of a quasiregular mapping has measure zero, which was expected to be true in~\cite[Remark 1.2]{flp14}. This fact was indeed a corollary of the more general results obtained recently in~\cite{gw14}. We prefer presenting the alternative approach here since it is more elementary and the basic idea behind the proof is similar to the one used in the Euclidean case.

Recall that for each $p\in M$ and $r$ small enough (depending on $p$), the exponential mapping $\exp_p:U\to B(p,r)$ is a homeomorphism and satisfies
\begin{align}\label{eq:almost bi-Lipschitz exp}
	L^{-1}d(p,q)\leq d^p\big(\exp_p^{-1}(p),\exp_p^{-1}(q)\big)\leq Ld(p,q)
\end{align}
for all $q\in B(p,r)$, where $L=L(p)$ is a constant depending on $p$.

\begin{theorem}\label{thm:size of the Branch set}
	Let $f:M\to N$ be a weak metrically quasiregular mapping. Then
	\begin{align*}
		\Vol_M(\mathcal{B}_f)=\Vol_N(f(\mathcal{B}_f))=0.
	\end{align*}
\end{theorem}
\begin{proof}
	Since $f$ satisfies Condition $N$, we only need to prove that $\mathcal{B}_f$ has measure zero. We will prove that if $p\in M$ is a point such that $f$ is $P$-differentiable at $p$ and $J_f(p)>0$, then $p\notin \mathcal{B}_f$.
	
	Before turning to the detailed proof, let us briefly indicate the idea: for a point $p\in M$ as above, the $P$-differential $Df(p)$ of $f$ at $p$ is a group isomorphism of the corresponding tangent Carnot groups and we may approximate $f$ in a neighborhood of $p$ by the $P$-differential composed with the exponential mapping. Locally, the exponential mapping is a homeomorphism that satisfies~\eqref{eq:almost bi-Lipschitz exp}. This implies that $f$ is close to a homeomorphism in a neighborhood of $p$ and so we may use the standard homotopy argument to show that the local index of $f$ at $p$ is $\pm 1$. This means $p\notin \mathcal{B}_f$. Since all such points have full measure in $M$, $\mathcal{B}_f$ must have measure zero. See~\cite{hr02} for the definition of local index of a continuous mapping between manifolds and its homotopy invariance.
	
	Fix such a point $p$. Choosing $r$ sufficiently small if necessary, we may assume that the mappings $\exp_p^{-1}|_{B(p,r)}:U\to B(p,r)$ and $\exp_{f(p)}^{-1}|_{f(B(p,r))}:V\to f(B(p,r))$ are homeomorphisms that satisfy~\eqref{eq:almost bi-Lipschitz exp} with some constant $L=L(p)$. Since $f$ is $P$-differentiable at $p$ and $J_f(p)>0$, the $P$-differential $Df(p):(\G_p,d^p)\to (\G_{f(p)},d^{f(p)})$ is a group isomorphism. We define a homotopy $H:[0,1]\times \closure{B}(p,r)\to (\G_{f(p)},d^{f(p)})$ between $Df(p)\circ \exp_p^{-1}(x)$ and $\exp_{f(p)}^{-1}\circ f(x)$ by
	\begin{align*}
		H_t(x)=Df(p)\circ\exp_p^{-1}(x)\cdot\big(\delta_{t}Df(p)\circ  \exp_p^{-1}(x)\big)^{-1}\cdot \big(\delta_t\exp_{f(p)}^{-1}\circ f(x)\big).
	\end{align*}
	It is clear $H_0(x)=Df(p)\circ  \exp_p^{-1}(x)$ is a homeomorphism and thus we only need to verify that $H_t$ is a proper homotopy between $H_0$ and $H_1$, \ie $\exp_{f(p)}^{-1}\circ f(p)\notin H_t\big(\bdary B(p,r)\big)$.
	
	Note that the minimal norm of $Df(p)$ defined as in~\eqref{eq:def for minimal norm}
	\begin{align*}
		\lambda=\|Df(p)\|_s>0,
	\end{align*}
	since $J_f(p)>0$. By Triangle inequality,
	\begin{align*}
		&d^{f(p)}\big(H_0(x)\cdot(\delta_tH_0(x))^{-1}\cdot\delta_t H_1(x),\exp_{f(p)}^{-1}\circ f(p)\big)\\
		&\geq d^{f(p)}\big(H_0(x),\exp_{f(p)}^{-1}\circ f(p)\big)-d^{f(p)}\big(H_0(x)\cdot(\delta_tH_0(x))^{-1}\cdot\delta_t H_1(x),H_0(x)\big)\\
		&=:I_1-I_2.
	\end{align*}
	We first estimate $I_1$ from \blue{below}:
	\begin{align*}
		I_1&=d^{f(p)}\big(H_0(x),\exp_{f(p)}^{-1}\circ f(p)\big)=d^{f(p)}\big(H_0(x),Df(p)\circ  \exp_p^{-1}(p)\big)\\
		&=d^{f(p)}\big(Df(p)\circ \exp_p^{-1}(x),Df(p)\circ  \exp_p^{-1}(p)\big)\\
		&\geq \frac{\lambda}{2}d^{p}\big(\exp_p^{-1}(x),\exp_p^{-1}(p)\big)\geq \frac{\lambda}{2L}d(x,p).
	\end{align*}
	We next estimate $I_2$ from \blue{above} as follows:
	\begin{align*}
		I_2&=d^{f(p)}\big(H_0(x)\cdot(\delta_tH_0(x))^{-1}\cdot\delta_t H_1(x),H_0(x)\big)\\
		&=d^{f(p)}\big((\delta_{t}H_0(x))^{-1}\cdot\delta_t H_1(x),0\big)\\
		&=d^{f(p)}\big(\delta_t H_1(x),\delta_t H_0(x)\big)\\
		&=td^{f(p)}\big(H_1(x),H_0(x)\big).
		%%&\leq t\Big(\bar{d}\big(H_1(x),H_1(p)\big)+\bar{d}\big(H_1(p),H_0(x)\big)\Big)
	\end{align*}
	Note that the $P$-differentiability of $f$ at $p$ implies that
	\begin{align*}
		d^{f(p)}\big(H_1(x),H_0(x)\big)&=d^{f(p)}\big(\exp_{f(p)}^{-1}\circ f(x),Df(p)\circ  \exp_p^{-1}(x)\big)\\
		&=o(d(x,p))\quad \text{ as } x\to p.
	\end{align*}
	This implies that  $\exp_{f(p)}^{-1}\circ f(p)\notin H_t\big(\bdary B(p,r)\big)$ for all sufficiently small $r$ and all $t\in [0,1]$. Therefore, $i(p,f)=i(p,H_1)=i(p,H_0)=\pm 1$ and $p\notin \mathcal{B}_f$.
\end{proof}

%%%%%%%%%%%%%%%%%%%%%%%%%%%%%%%%%%%%%%%%%%%%%%%%%%%%%%%%%%%%%%%

\subsection{Open questions}\label{subsec:open questions}
Since the basic analytic theory of quasiregular mappings between equiregular subRiemannian manifolds was established, in this subsection, we list some natural open questions for further research.

In the smooth setting, there is another well-known approach to establish the theory of quasiregular mappings based on the non-linear potential theory, see for instance~\cite{hkm06}. This approach was used by Heinonen and Holopainen~\cite{hh97} in the study of quasiregular mappings between Carnot groups. It is natural to ask the following question.

\begin{question}\label{ques:non-linear potential approach}
	Is it possible to establish the theory of quasiregular mappings between equiregular subRiemannian manifolds via the non-linear potential theory?
\end{question}

Note that the study of (sub-)$Q$-harmonic equation is necessary in order to handle Question~\ref{ques:non-linear potential approach}, which might be of independent interest.

A deep theorem of Reshetnyak~\cite{re89} says that non-constant analytic quasiregular mappings between Euclidean spaces are indeed both discrete and open. This result was further generalized by Heinonen and Holopainen~\cite{hh97} to the setting of certain Carnot groups, namely, non-constant analytic quasiregular mappings between Carnot groups of Heisenberg type are both discrete and open. It is natural to inquire whether this is a general fact for quasiregular mappings between equiregular subRiemannian manifolds.

\begin{question}\label{ques:openness and discreteness}
	Is it true that a non-constant analytically quasiregular mapping $f:M\to N$ between two equiregular subRiemannian manifolds is both discrete and open?
\end{question}

A first attempt to Question~\ref{ques:openness and discreteness} would be the study of the case when $M$ and $N$ are Carnot groups.

Another interesting result from~\cite{hh97} says that quasiregular mappings between certain Carnot groups are in fact conformal. Thus, we could ask the following question.

\begin{question}\label{ques:on existence of QR}
	For which subRiemannian manifolds $M$ and $N$ does there exist a non-trivial quasiregular mapping $f:M\to N$?
\end{question}

The well-known global homeomorphism theorem of Gromov and Zorich states that a locally homeomorphic quasiregular mapping $f:M\to N$ from a $n$-dimensional Riemannian manifold $M$ into a simply connected $n$-dimensional Riemannian manifold $N$ is a homeomorphism onto its image, provided $n\geq 3$. Moreover, the exceptional set $N\backslash f(M)$ has zero $n$-capacity; see for instance~\cite{hp04,g13}.

\begin{question}\label{ques:global homeomorphism theorem}
	Does the global homeomorphism theorem hold for locally homeomorphic quasiregular mappings between equiregular subRiemannian manifolds?
\end{question}

A classical result of Martio, Rickman and V\"ais\"al\"a~\cite{mrv71} states that there exists a constant $\varepsilon(n)>0$ such that every non-constant $K$-quasiregular mapping with $K\leq 1+\varepsilon(n)$ in dimension $n\geq 3$ is locally homeomorphic when $n\geq 3$. 
%see also~\cite{r05} for a quantitative estimate of $\varepsilon(n)$.
\begin{question}\label{ques:empty branch set when K close to 1}
	Is there an $\varepsilon>0$, depending only on the data of the equiregular subRiemannian manifolds, such that every $(1+\varepsilon)$-quasiregular mapping between two equiregular subRiemannian manifolds of homogeneous dimension $Q\geq 3$ is locally homeomorphic?
\end{question}

\begin{question}\label{ques:Picard type}
	Do the Picard type results from~\cite{r93} and \cite{bh01} and Bloch's Theorem from~\cite{r07} hold for quasiregular mappings between equiregular subRiemannian manifolds?
\end{question}

%%%%%%%%%%%%%%%%%%%%%%%%%%%%%%%%%%%%%%%%%%%%%%%%%%%%%%%%%%%%%%%
%%%%%%%%%%%%%%%%%%%%%%%%%%%%%%%%%%%%%%%%%%%%%%%%%%%%%%%%%%%%%%%
\section{Differentiability of Lipschitz mappings}\label{sec:Pansu diff}
This section is devoted to prove the Stepanov's Theorem~\ref{thm:Differentiability Stepanov}.
%\begin{theorem}[Stepanov's Theorem]\label{thm:stepanov}
%	Let $f:(M,d)\to (\bar M,\bar d)$ be a Borel function between subRiemannian manifolds.
%	Then $f$ is $P$-differentiable for a.e.~$o$ in the set
%	\[
%	L(f) := \left\{o\in M: \limsup_{p\to o} \frac{\bar d(f(o),f(p))}{d(o,p)} < \infty \right\}.
%	\]
%\end{theorem}

%%%%%%%%%%%%%%%%%%%%%%%%%%%%%%%%%%%%%%%%%%%%%%%%%%%%%%%%%%%%%%%

%%%%%%%%%%%%%%%%%%%%%%%%%%%%%%%%%%%%%%%%%%%%%%%%%%%%%%%%%%%%%%%
\subsection{$P$-differentiability}
Let $\bar M$ be a smooth manifold.
We have on $\bar M$ all the same objects as on $M$, and we distinguish them by putting a bar on the ones for $\bar M$.
Let $f:M\to \bar M$ be a Borel mapping, $o\in M$ and $\bar o:=f(o)\in\bar M$.
\begin{definition}[$P$-Differential]\label{def1738}
	We say that $f:M\to \bar M$ is $P$-differentiable at $o$ if there exists a morphism of graded Lie algebras $A:\g_o\to\g_{\bar o}$ such that
	\[
	\lim_{\g_o\ni X\to 0} \frac{\bar d\left(\exp(A[X])(\bar o),f(\exp(X)(o))\right)}{\|X\|}=0
	\]
	where $\|\cdot\|$ is any homogeneous norm on $\g_o$. When $f$ is $P$-differentiable at $o$, we write $Df(o)$ instead of $A$ for the $P$-differential.
\end{definition}

Notice that this definition of $P$-differential depends on a choice of two systems of privileged coordinates, one centered at $o$ and the other at $\bar o$.
However, different choices of privileged coordinates commutes by isomorphisms.
More precisely, if $\g'_o$ and $\bar\g'_{\bar o}$ are the graded Lie algebras that arise from a different choice of privileged coordinates, then there are isomorphisms of graded Lie algebras $\phi:\g_o\to\g'_o$ and $\bar\phi:\bar\g_{\bar o}\to\bar\g'_{\bar o}$ with the following property:
for any map $f:M\to \bar M$ with $f(o)=\bar o$ and with a $P$-differential $A:\g_o\to\g_{\bar o}$, a morphism of Lie algebras $A':\bar\g'_{ o}\to\bar\g'_{\bar o}$ is the $P$-differential of $f$ at $o$ if and only if the following diagram commutes:
\[
\xymatrix{
	\bar\g_{ o}\ar[d]_{\phi}\ar[r]^A & \bar\g_{\bar o}\ar[d]^{\bar\phi} \\
	\bar\g'_{ o}\ar[r]_{A'} & \bar\g'_{\bar o}
}
\]

As the metric $d$ is comparable to the metric $d^o$ and, further, to the  homogeneous norm $\|\cdot\|$
in privileged coordinates, in these coordinates the $P$-differential is a linear map $Df(0)=Df(o):\R^n\to\R^{\bar n}$ such that
\begin{equation}
\lim_{y\to 0} \frac{\bar d(Df(0)y, f(y))}{d(0,y)} = 0.
\end{equation}
As \[|\bar d(f(o),f(y))-\bar d(f(o),Df(0)y)|\leq\bar d(Df(0)y, f(y)),\]
we have
\begin{equation}
\bar d(f(o),f(y)) = \bar d(f(o),Df(0)y) + o\left( d(o,y) \right).
\end{equation}

\begin{remark}
	If both $M$ and $\bar M$ are Carnot groups, then our Definition \ref{def1738} of $P$-differential is the same as the classical \emph{Pansu differential}.
	Indeed, suppose for sake of simplicity that $o$ and $\bar o$ are the neutral elements of $M$ and $\bar M$ respectively.
	Then we can identify $\g_o=\G_o=M$ and $\g_{\bar o}=\G_{\bar o}=\bar M$ as sets through the exponential maps.
	As a homogeneous norm on $\g_o$, we can choose $\|X\|:=d(o,X)$.
	Therefore, cleaning up the notation in Definition \ref{def1738}, we have
	\[
	\lim_{\g_o\ni X\to 0} \frac{\bar d\left(\exp(A[X])(\bar o),f(\exp(X)(o))\right)}{\|X\|}
	= \lim_{M\ni X\to 0} \frac{\bar d\left(A[X],f(X)\right)}{d(o,X)}=0 ,
	\]
	which is the usual definition of Pansu differential.
\end{remark}
\begin{remark}
The definition of the Pansu differential here is equivalent to the definition of the Pansu differential of Definition 5.11 in \cite{flp14}. 
In particular, as $\bar {d}^{\bar o}$ is locally equivalent to $\bar d$ by \cite[Theorem 2.2]{j14}, our definition of Pansu differential is equivalent to 
\[
\lim_{\g_o\ni X\to 0} \bar d^{\bar o}\left(\delta_{\frac{1}{\|X\|}}\exp(A[X])(\bar o),\delta_{\frac{1}{\|X\|}}f(\exp(X)(o))\right)=0,
\]
which is the same meaning as, by setting $Y=\delta_{\frac{1}{\|X\|}}{X}$ and $h=\frac{1}{\|X\|}$, 
\[\lim_{h\to 0} \bar d^{\bar o}\left(\exp(A[Y])(\bar o),\delta_hf(\exp(\delta_{\frac{1}{h}}Y)(o))\right)=0,\]
uniformly for all $Y\in S$, where $S$ is the unit sphere under the homogeneous norm. 
\end{remark}

\begin{remark}
In the non-equiregular case, the $P$-differential is not well-defined. As an example, let $M=\R$ with $X_1=\de_t$, and $\bar M=\R^3$ with 
	\[
	\bar X_1 := \de_x ,
	\qquad
	\bar X_2 := y\de_z ,
	\qquad
	\bar X_3 := \de_y .
	\]
	Notice that $\delta_\lambda(x,y,z) = (\lambda x, \lambda y, \lambda^2 z)$ are dilations of $\bar M$, i.e., the tangent cone of $\bar M$ at $(0,0,0)$ is $\bar M$ itself.
	In particular, 
	\[
	\bar g_{(0,0,0)}=\Span\{\bar X_1, \bar X_2, \bar X_3\}\oplus\Span\{\de_z\}\subset\Vec(\R^3) .
	\]
	
	Then it is easy to see that the mappings $A_a:\Span\{\de_t\}\to \bar g_{(0,0,0)}$ defined by
	\[
	A(\de_t) = \bar X_1 + a \bar X_2
	\] 
	are all $P$-differentials of the map $f:M\to\bar M$ according to Definition \ref{def1738}, $f(t) := (t,0,0)$ at $(0,0,0)$.

\end{remark}

%%%%%%%%%%%%%%%%%%%%%%%%%%%%%%%%%%%%%%%%%%%%%%%%%%%%%%%%%%%%%%%
\subsection{A variant of the Lebesgue differentiation theorem}
The aim of this section is to show the following differentiation theorem, which will be used in our later proof of the Stepanov's theorem. For every $p\in M$ let $B_p\subset U^1_p$ be a compact neighborhood of $p$.

\begin{proposition}\label{prop:Lebesg diff them}
Let $\Omega\subset M\times[0,1]$ be an open neighborhood of $M\times\{0\}$ and let $\phi:\Omega\to M$, $(p,t)\mapsto \phi_tp$, be the flow of a smooth vector field on $M$ which is nonzero everywhere. If $h:M\to\R$ be a locally bounded function, then, for almost all $o\in M$
	\[
	\lim_{\epsilon\to0}
	\dashint_{\delta_\epsilon^o(B_o)}\dashint_0^\epsilon |h(\phi_sp)-h(o)|\dd s\dd p = 0
	\]
\end{proposition}

We will use a version of the Lebesgue differentiation theorem due to Federer~\cite[Theorem 2.9.8, Page 156--165]{f69}.

%%%%%%%%%%%%%%%%%%%%%%%%%%%%%%%%%%%%%%%%%%%%%%%%%%%%%%%%%%%%%%%
\begin{lemma}\label{lem1356}
 	Define
	\begin{equation}\label{eq1343}
	 	\scr V:= \left\{\big(p,\delta^p_\epsilon B_p\big) : p\in M,\ \epsilon\in(0,1]\right\} .
	\end{equation}
	The family $\scr V$ is a Vitali relation,  in the sense of Federer~\cite[\S 2.8.16]{f69}.
\end{lemma}
\begin{proof}
We do this using~\cite[Theorem 2.8.17]{f69}.
Using Federer's notation, in our case we choose
\[
\tau=2\qquad\text{ and }\qquad\delta\big(p,\delta^p_\epsilon B_p\big) = \diam_d(\delta^p_\epsilon B_p).
\]
We need only to show that for almost all $o\in M$:
\begin{equation}
 	\limsup_{\epsilon\to0} \frac{|\widehat{\delta^o_\epsilon B_o}|}{|\delta^o_\epsilon B_o|} <+\infty,
\end{equation}
where
\[
\widehat{\delta^o_\epsilon B_o} = \bigcup \delta^p_\eta B_p,
\]
the union is taken on all $\delta^p_\eta B_p$ such that $\delta^p_\eta B_p\cap\delta^o_\epsilon B_o \neq\emptyset$ and $\diam_d(\delta^p_\eta B_p) \le 2 \diam_d(\delta^o_\epsilon B_o)$.
Hence
\[
\widehat{\delta^o_\epsilon B_o} \subset B_d(o,3\diam_d(\delta^o_\epsilon B_o))
\]
and we have to prove
\begin{equation}\label{eq1516}
 	\limsup_{\epsilon\to0} \frac{|B_d\big(o,3\diam_d(\delta^o_\epsilon B_o)\big)|}{|\delta^o_\epsilon B_o|} <+\infty .
\end{equation}

First we claim that
\[
\diam_d(\delta^o_\epsilon B_o) = O( \epsilon)
\]
Indeed, if $x\in B_o$, then by \eqref{bab}
\[
d(o,\delta^o_\epsilon x) = \epsilon\epsilon^{-1}d(\delta^o_\epsilon o,\delta^o_\epsilon x)
= \epsilon d^{o}_{\epsilon}(o,x)
\]
where $ d^{o}_{\epsilon}(o,x)\to d^{o}_{0}(o,x)$ uniformly in $x$ as $\epsilon\to 0$, and therefore $d(o,\delta^o_\epsilon x) = O(\epsilon)$, uniformly in $x$.

Thanks to Theorem~\ref{thm:ball-box} and the fact that the Popp measure is smooth, as $d$ and $d_0^o$ are comparable by \cite[Theorem 2.2]{j14}, we have
\[
|B_d(o,r)| \sim |B_{d_0^o}(o,r)| \sim r^Q
\]
where $ Q$ is the homogeneous dimension of $M$ at $o$.

Finally, since $\delta^o_\epsilon$ has determinant equal to $\epsilon^Q$ and the Popp measure is smooth in the coordinate, we have
\[
|\delta^o_\epsilon B_o| \sim \epsilon^Q,
\]
which leads to \eqref{eq1516}. 	
\end{proof}
%%%%%%%%%%%%%%%%%%%%%%%%%%%%%%%%%%%%%%%%%%%%%%%%%%%%%%%%%%%%%%%
\begin{lemma}\label{lemma:33}
 	Define
	\[
	R_\epsilon(p):=\dashint_0^\epsilon |h(\phi_sp)-h(p)|\dd s .
	\]
	Then for almost every $p\in M$ we have
	\begin{equation}\label{eq1257}
	\lim_{\epsilon\to0} R_\epsilon(p) = 0 .
	\end{equation}
\end{lemma}
\begin{proof}
	Applying the Lebesgue differentiation theorem to the function $t\mapsto h(\phi_tq)$ for any $q\in M$, we obtain that for almost all $(q,t)$% for a.e. $t\in[0,1]$
	$$\lim_{\epsilon\to0}\dashint_0^\epsilon |h(\phi_{t+s}q)-h(\phi_tq)| \dd s = 0,$$
which implies \eqref{eq1257} holds for $p=\phi_tq$. Since the map $(q,t)\mapsto\phi_t q$ is locally Lipschitz and surjective, it maps a set of full measure into a set of full measure, therefore for almost every $p\in M$, \eqref{eq1257} holds.
\end{proof}
%%%%%%%%%%%%%%%%%%%%%%%%%%%%%%%%%%%%%%%%%%%%%%%%%%%%%%%%%%%%%%%
\begin{lemma}\label{lem1321}
 	Let $F\subset M$ be a measurable subset.
	Then for almost all $o\in F$
	\begin{equation}\label{eq1321}
	\lim_{\epsilon\to0}\frac{|\delta^o_\epsilon B_o\setminus F|}{|\delta^o_\epsilon B_o|} = 0.
	\end{equation}
\end{lemma}
\begin{proof}
 	Since $\scr V$, defined in \eqref{eq1343}, is a Vitali relation, the claim follows by applying the Lebesgue differentiation theorem to  the characteristic function $\chi_F$ of $F$.
\end{proof}
%%%%%%%%%%%%%%%%%%%%%%%%%%%%%%%%%%%%%%%%%%%%%%%%%%%%%%%%%%%%%%%
\begin{proof}[Proof of Proposition~\ref{prop:Lebesg diff them}]
%\begin{lemma}\label{lemma:35}
We need to show that for almost all $o\in M$
 	\begin{equation}\label{eq1323}
	 	\lim_{\epsilon\to0}\dashint_{\delta^o_\epsilon B_o}\dashint_0^\epsilon |h(\phi_sp)-h(o)|\dd s\dd p = 0
	\end{equation}
%\end{lemma}

Thanks to Egorov Theorem and Lemma~\ref{lemma:33}, for every $\eta>0$ there is a measurable subset $F\subset M$ such that $|M\setminus F| \le \eta$ and $R_\epsilon$ converge uniformly to $0$ on $F$. Since $\eta$ is arbitrary, it suffices to show that \eqref{eq1323} holds for almost all $o\in F$.
	
	Since $\scr V$, defined in \eqref{eq1343}, is a Vitali family, by Lemma \ref{lem1321}, we deduce that for almost every $o\in F$  \eqref{eq1321} holds and
	\begin{equation}\label{eq1327}
	 	\lim_{\epsilon\to0}\dashint_{\delta^o_\epsilon B_o} |h(p)-h(o)| \dd p = 0 .
	\end{equation}
	
	For such an $o$ we have
 	\begin{multline*}
	\dashint_{\delta^o_\epsilon B_o}\dashint_0^\epsilon |h(\phi_sp)-h(o)|\dd s\dd p
	\\ \le
	\underbrace{\dashint_{\delta^o_\epsilon B_o}\dashint_0^\epsilon |h(\phi_sp)-h(p)|\dd s\dd p}_{A}
	+
	\underbrace{\dashint_{\delta^o_\epsilon B_o} |h(p)-h(o)| \dd p}_{B}\\
	\end{multline*}
	where part $B$ converges to $0$ as $\epsilon\to0$ because of \eqref{eq1327}.
	For part $A$, we have
	\begin{align*}
	 	A&=\dashint_{\delta^o_\epsilon B_o} R_\epsilon(p) \dd p
		= \frac1{|\delta^o_\epsilon B_o|}\int_{\delta^o_\epsilon B_o\cap F}R_\epsilon(p) \dd p
			+ \frac1{|\delta^o_\epsilon B_o|}\int_{\delta^o_\epsilon B_o\setminus F}R_\epsilon(p) \dd p \\
		&\le \dashint_{\delta^o_\epsilon B_o} R_\epsilon(p) \chi_F \dd p + C \frac{|\delta^o_\epsilon B_o\setminus F|}{|\delta^o_\epsilon B_o|},
	\end{align*}
	where $C>0$ is some constant that bounds $R_\epsilon$, which exists because $h$ is locally bounded. It is now clear that $A$ converges to $0$ as $\epsilon\to 0$.
\end{proof}

\subsection{Blow-up of Lipschitz functions}
Let $E\subset M$ be closed and let $f:E\to \bar M$ be an $L$-Lipschitz function. For every $p\in M$ let $B_p\subset M$ be a closed $d$-ball centered at $p$ such that the dilations $\delta^p_\epsilon$ are well defined for $\epsilon\in(0,1]$.
For $o\in E$, we define the functions
\[
f^{o,\epsilon}:=\bar{\delta}_{\frac1\epsilon}\circ f\circ\delta^o_\epsilon,
\]
where $\bar{\delta}_{\frac1\epsilon}$ is with respect to $f(o)$.

Define with these $B_p$
\[
 	\scr V:= \left\{\big(p,\delta^p_\epsilon B_p\big) : p\in M,\ \epsilon\in(0,1]\right\} .
\]
Lemma \ref{lem1356} implies that $\scr V$ is a Vitali relation, and hence almost every~$o\in E$ is a $\scr V$-density point of $E$.

Our main step is to show a sort of equicontinuity of $f^{o,\epsilon}$ at almost every point of $E$.

\begin{lemma}\label{lem1350}
 	For every $o\in E$, setting $\bar o=f(o)\in \bar M$, there is an open neighborhood of $o$ $U^o$ such that $U^o\subset U^o_1$ and $f^{o,\epsilon}$ is a well-defined map $U^o\cap \delta^{o}_{\frac1\epsilon}E\rightarrow U_1^{\bar o}$.
\end{lemma}
\begin{proof}
Fix a point $o\in E$ and choose $r>0$ such that \[B_{\bar{d}_0^{\bar{o}}}(\bar {o},2r)\subset U_1^{\bar {o}}.\]
Then, it is easy to deduce that
\[
B_{\bar{d}_0^{\bar{o}}}(\bar{o}, \epsilon r)=\bar{\delta}_\epsilon\left(B_{\bar{d}_0^{\bar{o}}}(\bar{o}, r)\right)  \subset \bar{\delta}_{\epsilon}U_{1}^{\bar{o}}.\]
Since $d$ and $d^o_0$ are comparable (by \cite[Theorem 2.2]{j14}) and $d^o_0$ is  homogeneous, there exist constants $A, B>0$ such that
\[
\epsilon^{-1} d\left(o, \delta_\epsilon p\right) \leq A d_0^o(o, p)\quad \text{and}\quad \bar{d}_0^{\bar{o}}(\bar{o}, \bar{p})\leq B \bar{d}(\bar{o}, \bar{p}).
\]
Thus,
\[
\bar{d}_0^{\bar{o}}\left(\bar{o}, f\left({\delta}_\epsilon p\right)\right) \leq B \bar{d}\left(\bar{o}, f\left({\delta}_\epsilon p\right)\right) \leq B L d\left(o, \delta_\epsilon p\right) \leq B L A \epsilon d_0^o(o, p),
\]
where $L$ is the Lipschitz constant of $f$.
Then, as long as $p\in U^o\cap {\delta}_{\frac{1}{\epsilon}}^oE$ satisfying that
$$
d_0^o(o, p) \leq \frac{r}{B L A} ,
$$
$\bar{\delta}_{\frac{1}{\epsilon}}\left(f\left(\delta_\epsilon p\right)\right)$ is well-defined and belongs to $U_{1}^{\bar{o}}$.
\end{proof}

Let $p\in E$ be a density point of $E$. We introduce the following type of convergence, adapted to the fact that $f^{p,\epsilon}$ is not defined in a neighborhood of $p$. Set $\bar p=f(p)$.
\begin{definition}\label{120}
Let $\epsilon_k\to 0$ be a sequence and $g$ be a continuous function on $B_p$. We say that $f^{p,\epsilon_k}\to g$ uniformly on $B_p$ if \[\sup\big\{\bar d^{\bar p}(f^{p,\epsilon_k}(q),g(q)):q\in\delta^p_{\frac1{\epsilon_k}}E\cap B_p\big\}\to 0.\]
\end{definition}
We denote the Hausdorff distance under the metric $d_0^o$ to be $d_H^o$ and the ball under the metric $d_0^o$ to be $B^o$.

\begin{lemma}\label{130}
 	If $o\in E$ is a $\scr V$-density point of $E$, then as $\epsilon_k\to 0$, $d_H^o(\delta^o_{\frac1{\epsilon_k}}E\cap B_o,B_o)\to 0$.
\end{lemma}
\begin{proof}
 	Assume that the thesis is false.
	Then, up to passing to a subsequence, there is $a>0$ such that for every $k\in\N$ there is $x_k\in B_o$ with
	\[
	B_o\cap B^o(x_k,a) \subset B_o\setminus \delta^o_{\frac1{\epsilon_k}}E
	\]
	i.e.,
	\[
	\lim_{k\to\infty} \frac{|\delta^o_{\epsilon_k}B_o\setminus E|}{|\delta^o_{\epsilon_k}B_o|}
	\sim \lim_{k\to\infty} \frac{ |B_o\setminus \delta^o_{\frac1{\epsilon_k}}E| }{|B_o|}
	\ge \lim_{k\to\infty} \frac{|B_o\cap B^o(x_k,a)|}{|B_o|} >0,
	\]
	i.e., $o$ is not a $\scr V$-density point of $E$.
\end{proof}
\begin{remark}\label{ttp}
Lemma~\ref{130} tells us that \[\lim_{j\to\infty}\sup_{y\in B_o}d^o(y,B_o\cap\delta_{\frac{1}{\epsilon_k}}E)=0.\] 
Thus, for every $y\in B_o$, there exists a sequence $y_k\in B_o\cap\delta_{\frac{1}{\epsilon_k}}E$ converges to $y$ uniformly with respect to $y$.
\end{remark}

%{\color{red} Change the proofs of Lemma~\ref{lem1403} and Proposition~\ref{prop2201}!}
%%%%%%%%%%%%%%%%%%%%%%%%%%%%%%%%%%%%%%%%%%%%%%%%%%%%%%%%%%%%%%%
	Set $\bar o:=f(o)$. The family of functions $\{f^{o,\epsilon}\}_{\epsilon\in(0,1]}$ is called \emph{eventually equicontinuous} with respect to $d^o$ and $\bar {d}^{\bar o}$ if
\begin{equation}\label{eq0047}
	\begin{array}{c}
	 	\forall \eta>0,\ \exists\delta>0,\ \forall\epsilon\in(0,\delta),\ \forall p,q\in \delta^o_{\frac1{\epsilon}}E\cap B_o \\
		d^o(p,q)\le \delta \quad\THEN\quad \bar {d}^{\bar o}(f_\epsilon(p),f_\epsilon(q)) \le \eta.
	\end{array}
\end{equation}

\begin{lemma}[Equicontinuity of $f^{o,\epsilon}$]\label{lem1403}
 	Let $o\in E$ be a $\scr V$-density point of $E$.
	Then the family $\{f^{o,\epsilon}\}_{\epsilon}$ is eventually equicontinuous with respect to $d^o$ and $\bar {d}^{\bar o}$.
\end{lemma}
\begin{proof}
	First, since $d^o_\epsilon\to d^o_0$ and $\bar d^{\bar o}_\epsilon\to \bar d^{\bar o}_0$ uniformly, then it holds
	\[
	\forall\eta>0,\ \exists\delta>0,\ \forall\epsilon\in (0,\delta),\ \forall p,q\in B_o,
	\quad |d_0^o(p,q)-d_\epsilon^o(p,q)|\le \eta
	\]
	and
	\[
	\forall\eta>0,\ \exists\delta>0,\ \forall\epsilon\in(0,\delta),\ \forall p,q\in B_o,
	\quad  |\bar d_\epsilon^{\bar o}(p,q)-\bar d_0^{\bar o}(p,q)|\le \eta
	\]
	
	On the other hand, since $f$ is $L$-Lipschitz, it follows easily that
	\begin{align*}
	\bar d^{\bar {o}}_\epsilon(f^{{ {o}},\epsilon}(p),f^{o,\epsilon}(q)) &=\bar{d}^{\bar{o}}_{\epsilon}(\bar{\delta}_{\frac{1}{\epsilon}}f\delta_{\epsilon}(p),\bar{\delta}_{\frac{1}{\epsilon}}f\delta_{\epsilon}(q)) \\ &=\epsilon^{-1}\bar{d}(f\delta_{\epsilon}(p),f\delta_{\epsilon}(q)) \\
&\le \epsilon^{-1}Ld(\delta_{\epsilon}(p),\delta_{\epsilon}(q))
\\
&= Ld^o_\epsilon(p,q).
	\end{align*}
%	Now we can prove that \eqref{eq0047} holds for $f^{o,\epsilon}$.
	Together, these two facts give \eqref{eq0047}.
\end{proof}
A set of finite points is called a $\frac{1}{n}$-net of $B_o$ if the balls with radius $\frac{1}{n}$, centered as these points cover $B_o$.
\begin{lemma}\label{110}
Let $o\in E$ be a $\scr V$-density point of $E$.
Then the family $\{f^{o,\epsilon}\}_{\epsilon}$ has a subsequence uniformly convergent to a continuous function $g$ on $B_o$ in the sense of Definition~\ref{120}.
\end{lemma}
\begin{proof}
Based on Lemma~\ref{130} and Lemma~\ref{lem1403}, we shall prove that for a sequence $\epsilon_k\to 0$, there exists a continuous function $g$ on $B_o$ and a subsequence of $\{f^{o,\epsilon_k}\}$ uniformly convergent to $g$. 

Set $f_n=f_{\epsilon_n}$. As $B_o$ is compact, for every $n$ we could choose a set of finite points such that the balls with radius $\frac{1}{n}$, centered as these points, cover $B_o$. Thus, by taking these sets in order, we obtain a sequence of points $\{z_n\}$ such that $z_1$, $z_2$, ..., $z_{k_n}$ are $\frac{1}{n}$-net. Thus, $\{z_k\}$ is dense subset of $B_p$. For $z_k$, denote the sequence as in Remark~\ref{ttp} by $x_k^n$. 
 
Then, we start to choose a subsequence. By Lemma~\ref{lem1403}, by choosing $B_o$ possibly smaller, we could assume that $f_n(B_o)$ is contained in a compact set. Thus, for $\{f_n(x_1^n)\}$, choose a subsequence $\{n^1_k\}$ such that $\{f_{n^1_k}(x^{n^1_k}_1)\}$ is a Cauchy sequence. Next, for $\{f_{n^1_k}(x^{n^1_k}_2)\}$, choose a subsequence $\{n^2_k\}$ such that $\{f_{n^2_k}(x^{n^2_k}_2)\}$ is a Cauchy sequence. Then, we constantly continue this process. Moreover, we fix $f_{n_k^k}$ as our subsequence. For every $i$, $\{f_{n_k^k}(x^{n^k_k}_i)\}$ is a Cauchy sequence. 

Next, we find the function $g$. For every $\eta$, there exists a constant $\delta$ as in the eventually equicontinuous property. There exists $N$ such that for every $k\geq N$, we have $\epsilon_{n_k^k}< \delta$. There exists $N^{\prime}$ such that $\frac{1}{N^{\prime}}<\frac{\delta}{3}$ and $\{z_1, ..., z_{k_{N^{\prime}}}\}$ is a $\frac{1}{n}$-net. There exists $N^{\prime\prime}$ such that for every $k_1,k_2>N^{\prime\prime}$ and every $1\leq i\leq k_{N^{\prime}}$ \begin{equation}\label{999}
\bar{d}^{\bar o}(f_{n_{k_2}^{k_2}}(x_i^{n_{k_2}^{k_2}}),f_{n_{k_1}^{k_1}}(x_i^{n_{k_1}^{k_1}}))<\eta
\end{equation} 
and for every $m>N^{\prime\prime}$ \[d^o(x^m_i,z_i)<\delta/3.\]
For every $a\in B_o$, there exists a sequence $a_k$ as in Remark~\ref{ttp}. There exists a constant $N^{\prime\prime\prime}$ which only depends on $\delta$, such that for every $k>N^{\prime\prime\prime}$ \[d^o(a_k,a)<\delta/3.\] The $\frac{1}{N^{\prime}}$-net property tells us that there exists $1\leq i\leq k_{N^{\prime}}$ such that \[d^o(z_i,a)<\frac{1}{N^{\prime}}<\frac{\delta}{3}.\]
Thus, for $k_1,k_2>\max(N,N^{\prime\prime},N^{\prime\prime\prime})$, by the above three inequalities, we have 
\[d^o(a_{n^{k_1}_{k_1}},x_i^{n^{k_1}_{k_1}})<\delta\quad\text{and}\quad d^o(a_{n^{k_2}_{k_2}},x_i^{n^{k_2}_{k_2}})<\delta.\]
Therefore, we use the eventually equicontinuous property to obtain 
\[\bar{d}^{\bar o}(f_{n_{k_1}^{k_1}}(a_{n^{k_1}_{k_1}}),f_{n_{k_1}^{k_1}}(x_i^{n^{k_1}_{k_1}}))<\eta\quad\text{and}\quad \bar{d}^{\bar o}(f_{n_{k_2}^{k_2}}(a_{n^{k_2}_{k_2}}),f_{n_{k_2}^{k_2}}(x_i^{n^{k_2}_{k_2}}))<\eta.\]
By \eqref{999} and the above two inequalities, we deduce
\begin{equation}\label{oppp}
\bar{d}^{\bar o}(f_{n_{k_2}^{k_2}}(a_{n^{k_2}_{k_2}}),f_{n_{k_1}^{k_1}}(a_{n^{k_1}_{k_1}}))<3\eta.
\end{equation}
Thus, we obtain a Cauchy sequence $\{f_{n_{k}^{k}}(a_{n^{k}_{k}})\}$ whose limit is defined to be $g(a)$.

Next, we prove that the function $g(a)$ is well-defined by showing $g(a)$ is independent of the choice of the sequence $a_k\to a$. Choose another sequence $a^{\prime}_k\to a$ as Remark~\ref{ttp}. Then, there exists $N^*$ such that for any $k>N^*$, $d^o(a_k,a)<\delta/2$ and $d^o(a^{\prime}_k,a)<\delta/2$. Thus, $d^o(a^{\prime}_k,a_k)<\delta$. Thus, for $k>\max(N^{\star},N)$, \[\bar{d}^{\bar o}(f_{n^k_k}(a_{n^k_k}),f_{n^k_k}(a^{\prime}_{n^k_k}))<\eta.\]
Sending $k\to\infty$, $\bar{d}^{\bar o}(g(a),g^{\prime}(a))\leq \eta$, for any $\eta>0$. Thus, $g(a)=g^{\prime}(a)$.

Next, we prove that the function $g$ is continuous. For any $a\in B_o$ and any $a^{\prime}\in B_0$ satisfying $d^o(a,a^{\prime})<\frac{\delta}{3}$, considering the sequence $a_k\to a$ and $a^{\prime}_k\to a$ as in Remark~\ref{ttp}, there exists $N$ only depending on $\delta$, for any $k>N$, $d^o(a_k,a)<\frac{\delta}{3}$ and $d^o(a^{\prime}_k,a^{\prime})<\frac{\delta}{3}$. Thus, $d^o(a^{\prime}_k,a_k)<\delta$, which implies \[\bar{d}^{\bar o}(f_{n^k_k}(a_{n^k_k}),f_{n^k_k}(a^{\prime}_{n^k_k}))<\eta.\]
Sending $k\to \infty$, $\bar{d}^{\bar o}(g(a^{\prime}),g(a))\leq\eta$. Thus, $g$ is continuous at $a$.

Finally, we estimate the uniform convergence as required. As the convergence of the Cauchy sequence $\{f_{n_{k}^{k}}(a_{n^{k}_{k}})\}$ is uniform, for any $\eta>0$, there exists $M$, for any $j\geq M$ and any $q\in B_o\cap\delta_{\frac{1}{\epsilon_{n^j_j}}}E$, considering the sequence $q_k$ as in Remark~\ref{ttp}, \[\bar{d}^{\bar o}(f_{n^j_j}(q_{n^j_j}),g(q))<\eta.\]
Moreover, for $\delta$ as above, there exists $M^{\prime}$ only depending on $\delta$ such that for any $j>M^{\prime}$, $d^o(q_{n_j^j},q)<\delta$. Thus, for any $j\geq \max(N,M^{\prime})$, where $N$ is as above, \[\bar{d}^{\bar o}(f_{n^j_j}(q_{n^j_j}),f_{n^j_j}(q))<\eta.\]   
The above two inequalities tell us that for any $j\geq \max(N,M,M^{\prime})$, \[\bar{d}^{\bar o}(g(q),f_{n^j_j}(q))<2\eta,\] 
from which we obtain the conclusion of this lemma.
\end{proof}
Finally, we introduce the following lemma that gives the existence of $P$-differential.
\begin{lemma}\label{lem2002}
Let $o$ be a density point of $E$. Suppose that $f^{o,\epsilon}\to f^{o,0}$ uniformly in the sense of Definition~\ref{120} and that there is a morphism of Lie algebras $A:\g_o\to\g_{\bar o}$ and a constant $\epsilon_0>0$ such that
\begin{equation}\label{140}
f^{o,0}(\exp(X)(o)) = \exp(A[X])(\bar o)
\end{equation}
holds for all $\|X\|\le \epsilon_0$. Then $f$ is $P$-differentiable at $o$ and $Df(o)=A$.
\end{lemma}
Note that the limit in $ P$-differentiability in the above lemma is with respect to only $X\in g_o$ such that $\exp (X)\in E$.
\begin{proof}
We use the exponential coordinate as our local privileged coordinate. As $f^{o,0}$ is the uniform limit of $f^{o,\epsilon}$ and the fact that $f^{o,\tau}\delta_{\epsilon}=\bar{\delta}_{\epsilon}f^{o,\tau\epsilon}$, we have \begin{equation}\label{131}
f^{o,0}(\exp(\delta_\epsilon X)(o))=\delta_\epsilon f^{o,0}(\exp(X)(o)).
\end{equation}
Then, by \eqref{140} and the above equation
\begin{equation*}
\exp(A[\delta_\epsilon X])(\bar o) =f^{o,0}(\exp(\delta_\epsilon X)(o))=\delta_\epsilon f^{o,0}(\exp(X)(o))=\delta_\epsilon \exp(A[X])(\bar o).
\end{equation*}
Thus, $A[\delta_\epsilon X]=\delta_\epsilon A[X]$ and $A$ is a morphism of graded Lie algebra. 

For $\epsilon\leq \epsilon_0$, any vector of norm $\epsilon$ can be written as $\delta_\epsilon X$ with $\|X\|=1$. Moreover, by \eqref{140} and \eqref{131}, we have
	\begin{align*}
	\frac{ \bar d(\exp(A[\delta_\epsilon X])(\bar o) , f(\exp(\delta_\epsilon X)(o)) }{\|\delta_\epsilon X\|} &=
	 \epsilon^{-1} \bar d(f^{o,0}(\exp(\delta_\epsilon X)(o)), f(\exp(\delta_\epsilon X)(o))) \\
	 &= \epsilon^{-1} \bar d(\bar\delta_\epsilon f^{o,0}(\exp(X)(o)), \bar\delta_\epsilon \bar\delta_{\frac1\epsilon} f(\delta_\epsilon \exp( X)(o)))  \\
	 &=  \bar d_\epsilon^o ( f^{o,0}(\exp(X)(o)), f^{o,\epsilon} (\exp(X)(o))).
	\end{align*}
Since both $d_\epsilon^o$ and $f^{o,\epsilon}$ converge uniformly, we get
	\[
	\lim_{\epsilon\to0} \bar d_\epsilon^o ( f^{o,0}(\exp(X)(o)), f^{o,\epsilon} (\exp(X)(o))) = 0
%	\bar d_0^o ( f^{o,0}(\exp(X)(o)), f^{o,0} (\exp(X)(o)))
	\]
	and the limit is uniform with respect to $X$.
\end{proof}

%%%%%%%%%%%%%%%%%%%%%%%%%%%%%%%%%%%%%%%%%%%%%%%%%%%%%%%%%%%%%%%
%%%%%%%%%%%%%%%%%%%%%%%%%%%%%%%%%%%%%%%%%%%%%%%%%%%%%%%%%%%%%%%
\subsection{Blow-up of horizontal vector fields}
Let $W=\sum_{j=1}^r w_j X_j$ be a horizontal vector field which is nonzero everywhere and let $(p,t)\mapsto \phi_tp$ be its flow, where $w_j$($1\le j \le r$) are smooth functions.

%We extend $f$ outside of $E$ in the following way.

Set \begin{equation}\label{babbb}
F(p,t) := f(\phi_tp).
\end{equation}
Then $F$ is well-defined on $\{(p,t):\phi_tp\in E\}$.
Since
$$\bar d(f(\phi_tp),f(\phi_sq)) \le L d(\phi_tp,\phi_sq),$$
$F$ is locally Lipschitz, where the Lipschitz constant depends only on the Lipschitz constant of $f$ and the Lipschitz constant of $\phi$ on a compact set.

\subsubsection{Extension of $F$ on $M\times\R$}
We seek an extension of $F$ (still denoted by $F$) on $M\times\R$ satisfying properties:
\begin{enumerate}
\item[(P1)] $F(p,t+s) = F(\phi_tp,s)$ for all $p\in M$ and all $s,t\in\R$ such that $\phi_tp$ exists.
\item[(P2)] For each $p\in M$ the curve $t\mapsto F(p,t)$ is locally Lipschitz.
\end{enumerate}

We first extend $F$ on $E\times\R$ in such a way that each curve $t\mapsto F(p,t)$ is locally Lipschitz.
More precisely: for $p\in E$ define $I_p:=\{t\in\R:\phi_tp\in E\}\subset\R$.
Since $E$ is closed, $I_p$ is closed as well.
Let $\hat t\in\R\setminus I_p$. Then there are two cases. In the first case, it happens that there is $t_1\in I_p$ such that  $\hat t\in(-\infty, t_1)\subset\R\setminus I_p$ or $t\in ( t_1,+\infty)\subset\R\setminus I_p$.
Then we set $F(p,t) = F(p,t_1)$.
In the second case, there are $t_1,t_2\in I_p$ with $\hat t\in (t_1,t_2)\subset\R\setminus I_p$.
Then
\[
\bar d(f(\phi_{t_1}p),f(\phi_{t_2}p)) \le Ld(\phi_{t_1}p,\phi_{t_2}p) \le  \tilde L |t_1-t_2|.
\]
Therefore there is a geodesic $\gamma:[t_1,t_2]\to\bar M$ joining $f(\phi_{t_1}p)$ to $f(\phi_{t_2}p)$ with constant velocity, i.e. $\tilde L$-Lipschitz.
In this case, we define $F(p,t)=\gamma(t)$ for $t\in(t_1,t_2)$, where the curve is chosen in such a manner that (P1) holds for all $p,\phi_t p\in E$.

Moreover, on the set $E^{\prime}=\{p:\text{there exists } t \text{ such that }\phi_tp\in E\}$, we may extend $F$ using the rule $F(p,t+s) = F(\phi_tp,s)$. For $p$ outside of this set, we simply define $F(p,t)=\bar p$ for some fixed point $\bar p\in\bar M$. 

Next, we verify the property (P1). In the case $p\in E$, if $\phi_t p\in E$, then according to our construction we already obtain the property (P1); if $\phi_t p\not\in E$, the construction on $E^{\prime}$ tells us (P1) holds as well. When $p\in E^{\prime}$, we have $\phi_tp\in E^{\prime}$. Then, $F(p,t+s) = F(\phi_{t^{\prime}}p,t+s-t^{\prime})$, where $\phi_{t^{\prime}}p\in E$. Thus, $F(\phi_tp,s) = F(\phi_{t^{\prime}}p,s-t^{\prime}+t)=F(p,t+s)$. Therefore, (P1) holds. In the case $p\not\in E^{\prime}$, (P1) is true as $F(p,t)$ is constant. 

Finally, we verify the property (P2). In the case $p\in E$, our choice of curves tells us the property holds. When $p\in E^{\prime}$, according to the rule $F(p,t+s) = F(\phi_tp,s)$, where $\phi_tp\in E$, we have the local Lipschitz property from the fact that the flow is local Lipschitz. For $p$ outside these two sets, as $F(p,t)$ is constant, (P2) holds. 
\begin{remark}\label{210}
The construction tells us that the Lipschitz constant in the condition (P2) is uniform with respect to $(p,t)$ in a compact set.
\end{remark}
%%%%%%%%%%%%%%%%%%%%%%%%%%%%%%%%%%%%%%%%%%%%%%%%%%%%%%%%%%%%%%%
\subsubsection{Blow-up of $F$}
According to (P1) and \cite[Proposition 3.50]{cambridge}, there are $h_j:M\times\R\to\R$, $j\in\{1,\dots,\bar r\}$ such that for each $p\in M$,
\[
\frac{\de F}{\de t} (p,t) = \sum_{j=1}^{\bar r} h_j(p,t) \bar X_j(F(p,t))
\]
holds for almost every $t$.
Notice that, by Remark~\ref{210}, $h_j$'s are locally bounded.
Furthermore, by (P1), we have 
\[\frac{\de F}{\de (t+s)}(p,t+s)=\frac{\de F}{\de s}(\phi_t p,s)\]
holds for almost every $r+s$. By setting $s=0$ and $t=t_0$, we obtain
\[\frac{\de F}{\de t}(p,t_0)=\frac{\de F}{\de t}(\phi_{t_0} p,0)\]
holds for almost every $t_0$. 

Using the fact that the flow $\phi_tp$ is locally Lipschitz and surjective, we infer that for almost every $p\in M$ the derivative $\frac{\de F}{\de t} (p,0)$ exists.
Thus, for almost all $p$, there holds
\[\frac{\de F}{\de t} (p,0)= \sum_{j=1}^{\bar r} h_j(p,0) \bar X_j(F(p,0)).\]

The aim of this section is to blow-up both manifolds $M$ and $\bar M$ keeping track of the map $F$.
The result is, in some sense, the flow of a left-invariant vector field on $\bar M$; See Proposition \ref{prop1044}.

For $o\in M$ and $\epsilon\in(0,1]$ we have the vector fields
\[
W^{o,\epsilon} := \epsilon\cdot\dd\delta^o_{\frac1\epsilon}\circ W\circ\delta^o_\epsilon
\]
whose flow is denoted by $\phi^{o,\epsilon}_tp := \exp(tW^{o,\epsilon})(p)$.
By \cite{b96}, $W^{o,\epsilon}\to W^{o,0} = \sum_{j=1}^rw_j(o) X_j^{o,0}$ uniformly on compact sets and so $\phi^{o,\epsilon}\to\phi^{o,0}$ uniformly for $(p,t)$ in compact sets.

For $o\in M$, $\bar o:=F(o,0)$ and $\epsilon>0$, set
\[
F^{o,\epsilon}(p,t) := \bar\delta^{\bar o}_{\frac1\epsilon}F(\delta^o_\epsilon p,\epsilon t) .
\]
Notice that
$
F^{o,\epsilon}(p,0) = f^{o,\epsilon}(p)
$
if $\delta^o_{\epsilon}p\in E$.

%{\color{red} The stetement of the next lemma is not clear!}
%{\color{blue} I changed it:}
\begin{lemma}\label{lemma:51}
	For all $o\in E$ and all $t\in\R$
	\begin{equation}\label{eq0036}
	F^{o,\epsilon}(p,t) = f^{o,\epsilon}(\phi^{o,\epsilon}_tp) ,
	\end{equation}
	if the right-hand side is well-defined.
	Moreover, for each $p$,
	\begin{equation}\label{eq2349}
	 	\frac{\de F^{o,\epsilon}}{\de t}(p,t) = \sum_{j=1}^{\bar r} h_j(\delta_\epsilon^op,\epsilon t) \bar X_j^{\bar o,\epsilon}(F^{o,\epsilon}(p,t))
	\end{equation}
holds for almost every $t$.
\end{lemma}

\begin{proof}
	Fix $p$ and set $\gamma(t)=\phi_tp$.
	Then $\gamma(0)=p$ and $\gamma'(t) = W(\gamma(t))$.
	Define $\eta_\epsilon(t)=\delta_{\frac1\epsilon}(\gamma(\epsilon t))$.
	Then $\eta_\epsilon(0)=\delta_{\frac1\epsilon}(p)$ and
	\begin{equation*}
	 	\eta_\epsilon' (t)
		= \dd\delta^o_{\frac1\epsilon}[\epsilon\gamma'(\epsilon t)]
		= \epsilon \dd\delta^o_{\frac1\epsilon}[W(\gamma(\epsilon t))] = \epsilon \dd\delta^o_{\frac1\epsilon}[W(\delta^o_\epsilon\delta^o_{\frac1\epsilon}\gamma(\epsilon t))]
		= W^{o,\epsilon}(\eta_\epsilon(t))
	\end{equation*}
	i.e. $\eta_\epsilon(t) = \exp(tW^{o,\epsilon})(\delta_{\frac1\epsilon}(p))$.
	In other words
	$%\begin{equation}
	\delta_{\frac1\epsilon}(\phi_{\epsilon t}p) = \phi^{o,\epsilon}_t(\delta_{\frac1\epsilon}p) .
	$
	Hence, if $\phi_{\epsilon t}\delta^o_\epsilon p\in E$, then by \eqref{babbb}, we have
	\begin{align*}
	F^{o,\epsilon}(p,t)= \bar {\delta}^{\bar o}_{\frac1\epsilon} F(\delta^o_\epsilon p,\epsilon t)
	=\bar {\delta}^{\bar o}_{\frac1\epsilon} \circ f\circ \phi_{\epsilon t}\delta^o_\epsilon p =\bar {\delta}^{\bar o}_{\frac1\epsilon} \circ f \circ \delta^o_\epsilon\circ\delta^o_{\frac1\epsilon} \circ \phi_{\epsilon t}\delta^o_\epsilon p
	= f^{o,\epsilon}(\phi^{o,\epsilon}_tp),
	\end{align*}
which gives \eqref{eq0036}. Regarding \eqref{eq2349}, we have
	\begin{align*}
	 	% \sum_{j=1}^{\bar r} h_j^{o,\epsilon}(p,t) \bar X_j^{o,\epsilon}(F^{o,\epsilon}(p,t))
		 \frac{\de F^{o,\epsilon}}{\de t}(p,t) %= \\
		 &= \frac{\de}{\de t} \bar \delta_{\frac1\epsilon}^{\bar o} F(\delta^o_\epsilon p, \epsilon t)
		 = \epsilon \dd\bar \delta_{\frac1\epsilon}^{\bar o} \frac{\de F}{\de t}\Big(\delta^o_\epsilon p, \epsilon t\Big)  \\
		 &= \epsilon \dd\bar \delta_{\frac1\epsilon}^{\bar o} \Big[\sum_{j=1}^{\bar r} h_j(\delta^o_\epsilon p,\epsilon t) \bar X_j(F(\delta^o_\epsilon p,\epsilon t))\Big] \\
		& = \sum_{j=1}^{\bar r} h_j(\delta^o_\epsilon p,\epsilon t) \epsilon \dd\bar \delta_{\frac1\epsilon}^{\bar o} \bar X_j(\bar\delta^{\bar o}_\epsilon \bar\delta^{\bar o}_{\frac1\epsilon}F(\delta^o_\epsilon p,\epsilon t)) \\
		 &= \sum_{j=1}^{\bar r} h_j(\delta^o_\epsilon p,\epsilon t) \bar X_j^{\bar o,\epsilon}(F^{o,\epsilon}(p,t)) .
	\end{align*}
\end{proof}

\begin{proposition}\label{prop1044}
 	Let $o\in E$ be a $\scr V$-density point of $E$. Let $\epsilon_k\to 0$ be a sequence such that $f^{o,\epsilon_k}$ converge uniformly to a continuous function $g:B_o\to \bar M$.
	Then there exists a neighborhood $K$ of the point $o$, such that for $t$ sufficiently small and $p\in K$,
	\[
	\exp\Big(t\sum_{j=1}^{\bar r} h_j(o,0)\bar X_j^{\bar o,0}\Big)\big(g(p)\big)
	=
	g\Big(\exp\big(t\sum_{j=1}^{r} w_j(o) X_j^{o,0}\big)(p)\Big),
	\]
holds for almost every density point $o\in E$.
\end{proposition}
%Here is the proof, which uses some estimates that are proven in the Lemmas below.
Since our notation is getting heavier and heavier, we will drop the subscript $k$ in $\epsilon_k$ and write just $\epsilon$.
\begin{proof}%[Proof of Proposition \ref{prop1044}]
	Define
	\[
	G(p,t) =g\Big(\exp\big(t\sum_{j=1}^{r} w_j(o) X_j^{o,0}\big)(p)\Big)
	\]
	where $w_j\in\Co^\infty(M)$ are the components of $W$ with respect to $X_1,\dots,X_r$.
	Now considering everything in the exponential coordinates, where the Euclidean distance in this coordinate is denoted by $|\cdot |$, the curve
	\[
	\gamma_p(t) = g(p) + \int_0^t \sum_{j=1}^{\bar r} h_j(o,0) \bar X_j^{\bar o,0}(G(p,s)) \dd s
	\]
	is well-defined.
	Notice that $\gamma_p(0)=g(p)$ and $\gamma_p'(t)=\sum_{j=1}^{\bar r} h_j(o,0) \bar X_j^{\bar o,0}(G(p,t))$.
	
In order to obtain the conclusion, it suffices to prove
	\[
	\gamma_p(t) = G(p,t),
	\]
for $t$ sufficiently small.
For this, as $g$ is continuous, we only need to prove that there exists a neighborhood $K$ of $o$ such that for $t$ sufficiently small,
	\begin{equation}\label{eq2345}
	\lim_{k\to\infty} \int_{K} |F^{o,\epsilon_k}(p,t)-G(p,t)|\dd p = 0
	\end{equation}
	and
	\begin{equation}\label{eq2346}
	\lim_{k\to\infty} \int_{K} |F^{o,\epsilon_k}(p,t)-\gamma_p(t)| \dd p = 0 .
	\end{equation}
The proof of \eqref{eq2345} is given in Lemma \ref{lem2204} below and we next prove \eqref{eq2346}.\\

Notice that thanks to \eqref{eq2349} we have
	\[
	F^{o,\epsilon}(p,t) = F^{o,\epsilon}(p,0) + \int_0^t \sum_{j=1}^{\bar r} h_j(\delta_{\epsilon}^op,\epsilon s) \bar X_j^{\bar o,\epsilon}(F^{o,\epsilon}(p,s))\dd s .
	\]
For $K=\delta^o_{\frac14}B_o$,
\begin{align*}
\int_{K}&\left| F^{o,\epsilon}(p,t) - g(p) - \int_0^t \sum_{j=1}^{\bar r} h_j(o,0) \bar X_j^{\bar o,0}(G(p,s)) \dd s \right| \dd p\\
&\leq \int_{K} \left| F^{o,\epsilon}(p,0) - g(p)\right| \dd p \\
&\qquad+\sum_{j=1}^{\bar r} \int_{K}\int_0^t \left|
		 h_j(\delta_{\epsilon}^op,\epsilon s) \bar X_j^{\bar o,\epsilon}(F^{o,\epsilon}(p,s)) -  h_j(o,0) \bar X_j^{\bar o,0}(G(p,s)) \right|\dd s \dd p\\
&\leq \underbrace{ \int_{B_o} \left| F^{o,\epsilon}(p,0) - g(p)\right| \dd p }_{(a)}\\
&\qquad +\sum_{j=1}^{\bar r}\underbrace{ \int_{B_o} \int_0^t |h_j(\delta^o_\epsilon p,\epsilon s) - h_j(o,0)| \cdot |\bar X_j^{\bar o,\epsilon}(F^{o,\epsilon}(p,s))| \dd s \dd p }_{(b)}\\
&\qquad\qquad\qquad +\sum_{j=1}^{\bar r}\underbrace{ \int_{K} \int_0^t |h_j(o,0)| \left| \bar X_j^{\bar o,\epsilon}(F^{o,\epsilon}(p,s)) - \bar X_j^{\bar o,0}(G(p,s)) \right| \dd s \dd p. }_{(c)}
\end{align*}
So we next estimate the three parts.

The proof of $(a)\to 0$ as $k\to \infty$ is given in Lemma \ref{lem2351} below.
	
To estimate part $(b)$, notice that $|\bar X_j^{\bar o,\epsilon}(F^{o,\epsilon}(p,s))|\le C$. We may use a change of variable and Proposition \ref{prop:Lebesg diff them} to infer that
	\begin{align*}
(b) \le  \frac C{\epsilon\cdot|\delta^o_\epsilon B_o|} \int_{\delta^o_\epsilon B_o} \int_0^\epsilon |h_j(p,s)-h_j(o,0)| \dd s \dd p
	\to 0
	\end{align*}
holds for almost every density point $o\in E$.	
Similarly, to estimate part $(c)$, we notice that $|h_j|\le C$ and hence
\begin{equation}\label{eq:for estimates}
	\begin{aligned}
	&\quad\int_{K} \int_0^t \left| \bar X_j^{o,\epsilon}(F^{o,\epsilon}(p,s)) - \bar X_j^{o,0}(G(p,s)) \right| \dd s \dd p\\
	&\le \int_0^t \int_{K}   \left| \bar X_j^{o,\epsilon}(F^{o,\epsilon}(p,s)) - \bar X_j^{o,0}(F^{o,\epsilon}(p,s)) \right|  \dd p \dd s\\ 
	&\qquad + \int_0^t \int_{K} \left| \bar X_j^{o,0}(F^{o,\epsilon}(p,s)) - \bar X_j^{o,0}(G(p,s)) \right| \dd p \dd s.
	\end{aligned}
\end{equation}
Since $\bar X_j^{o,\epsilon}\to \bar X_j^{o,0}$ uniformly,
	\[
	\left| \bar X_j^{o,\epsilon}(F^{o,\epsilon}(p,s)) - \bar X_j^{o,0}(F^{o,\epsilon}(p,s)) \right|\le \rho(\epsilon)
	\]
	for a function $\rho(\epsilon)$, independent on $F^{o,\epsilon}(p,s)$, with the property that $\rho(\epsilon_k)\to 0$ as $k\to\infty$. This implies that the first term in~\eqref{eq:for estimates} tends to zero as $k\to\infty$.

Regarding the second term in~\eqref{eq:for estimates}, observe that $\bar X_j^{o,0}$ is Lipschitz on compact sets. So we may use Lemma \ref{lem2204} below to conclude
	\begin{align*}
	\int_0^t \int_{K} &\left| \bar X_j^{o,0}(F^{o,\epsilon}(p,s)) - \bar X_j^{o,0}(G(p,s)) \right| \dd p \dd s \\
	&\le C\int_0^t \int_{K} \left|F^{o,\epsilon}(p,s) - G(p,s) \right| \dd p \dd s  \to 0,
	\end{align*}
from which \eqref{eq2346} follows.
\end{proof}
%%%%%%%%%%%%%%%%%%%%%%%%%%%%%%%%%%%%%%%%%%%%%%%%%%%%%%%%%%%%%%%
%%%%%%%%%%%%%%%%%%%%%%%%%%%%%%%%%%%%%%%%%%%%%%%%%%%%%%%%%%%%%%%
\begin{lemma}\label{lem2351}
 	\[
	\lim_{k\to\infty} \int_{B_o} |F^{o,\epsilon_k}(p,0)-g(p)|\dd p = 0
	\]
\end{lemma}
\begin{proof}
	Since $o$ is a $\scr V$-density point of $E$, we have
 	\begin{align*}
	 &\quad\int_{B_o} |F^{o,\epsilon}(p,0)-g(p)|\dd p\\
	 &\le  \int_{B_o\cap\delta^o_{\frac1\epsilon}E} |f^{o,\epsilon}(p)-g(p)|\dd p
	 	+  \int_{B_o\setminus\delta^o_{\frac1\epsilon}E} |F^{o,\epsilon}(p,0)-g(p)|\dd p \\
	&\le \int_{B_o\cap\delta^o_{\frac1\epsilon}E} |f^{o,\epsilon}(p)-g(p)|\dd p
		+ C |B_o\setminus\delta^o_{\frac1\epsilon}E|
	\to 0,
	\end{align*}
where we have used the fact that the Euclidean distance in the exponential coordinate is controlled by $\bar {d}_0^{\bar {o}}$.
\end{proof}
%%%%%%%%%%%%%%%%%%%%%%%%%%%%%%%%%%%%%%%%%%%%%%%%%%%%%%%%%%%%%%%
\begin{lemma}\label{lem2204}
	For $t$ sufficiently close to 0, it holds
	 	\[
	\lim_{k\to\infty} \int_{\delta^o_{\frac14}B_o} |F^{o,\epsilon_k}(p,t)-G(p,t)|\dd p = 0.
	\]
\end{lemma}
\begin{proof}
Choose $t$ small enough such that
	$\phi^{o,\epsilon}_t(\delta^o_{\frac14} B_o)\subset \delta^o_{\frac12} B_o$
	and
	$(\phi^{o,\epsilon}_t)^{-1}\delta^o_{\frac12} B_o\subset  B_o$.
Then,
\begin{align*}
&\int_{\delta^o_{\frac14}B_o}|F^{o,\epsilon}(p,t)-G(p,t)|\dd p\\
&=\int_{(\delta^o_{\frac14} B_o) \cap (\phi^{o,\epsilon}_t)^{-1}\delta^o_{\frac1\epsilon}E}
			|f^{o,\epsilon}(\phi^{o,\epsilon}_tp)-g(\phi^{o,0}_tp)|\dd p+\int_{(\delta^o_{\frac14} B_o) \setminus (\phi^{o,\epsilon}_t)^{-1}\delta^o_{\frac1\epsilon}E} |F^{o,\epsilon}(p,t)-G(p,t)|\dd p\\
&\le\int_{(\delta^o_{\frac14} B_o) \cap (\phi^{o,\epsilon}_t)^{-1}\delta^o_{\frac1\epsilon}E}
			|f^{o,\epsilon}(\phi^{o,\epsilon}_tp)-g(\phi^{o,\epsilon}_tp)|\dd p+\int_{(\delta^o_{\frac14} B_o) \cap (\phi^{o,\epsilon}_t)^{-1}\delta^o_{\frac1\epsilon}E}
			|g(\phi^{o,\epsilon}_tp)-g(\phi^{o,0}_tp)|\dd p	\\ &\qquad\qquad\qquad\qquad\qquad +\int_{(\delta^o_{\frac14} B_o) \setminus (\phi^{o,\epsilon}_t)^{-1}\delta^o_{\frac1\epsilon}E} |F^{o,\epsilon}(p,t)-G(p,t)|\dd p\\
\\\numberthis\label{eq:last equa}	
&\le C \int_{B_o\cap \delta^o_{\frac1\epsilon}E} |f^{o,\epsilon}(p)-g(p)|\dd p+
	\int_{\delta^o_{\frac14} B_o} |g(\phi^{o,\epsilon}_tp)-g(\phi^{o,0}_tp)|\dd p \\ &\qquad\qquad\qquad\qquad\qquad +C \int_{\delta^o_{\frac{1}{2}}B_o\setminus \delta^o_{\frac1\epsilon}E} |F^{o,\epsilon}({(\phi^{o,\epsilon}_t)}^{-1}p,t)-G({(\phi^{o,\epsilon}_t)}^{-1}p,t)|\dd p  \\
&\le C \int_{B_o\cap \delta^o_{\frac1\epsilon}E} |f^{o,\epsilon}(p)-g(p)|\dd p+
	\int_{\delta^o_{\frac14} B_o} |g(\phi^{o,\epsilon}_tp)-g(\phi^{o,0}_tp)|\dd p+ C|B_o\setminus \delta^o_{\frac1\epsilon}E|,
\end{align*}
which converges to $0$ as $\epsilon_k\to0$.
\end{proof}

\subsection{Stepanov's theorem}

\begin{proposition}[Rademacher's Theorem]
 	Let $E\subset M$ be a Borel set and $f:E\to\bar M$ a Lipschitz map.
	Then $f$ is $P$-differentiable almost everywhere in $E$.
\end{proposition}
Note that the limit in $ P$-differentiability in the above proposition is with respect to only $X$ belonging to the Lie algebra such that $\exp (X)\in E$.
\begin{proof}
Without loss of generality, assume that $E$ is a closed set. For $k\in\{1,\dots,r\}$ and $j\in\{1,\dots,\bar r\}$, let $h_{kj}(p,t)$ be such that, setting
    \[
    F_k(p,t) = f(\exp(tX_k)(p)).
    \]
By the extension of $F_k$ as in the above subsection, we have
    \[
    \frac{\de F_k}{\de t}(p,t) = \sum_{j=1}^{\bar r} h_{kj}(p,t) \bar X_j(F_k(p,t)) .
    \]
 	Let $o\in E$ be a $\scr V$-density point of $E$.
	Almost all $o\in E$ have this properties.
	
	Set $\bar o=f(o)$.
	We want to define a morphism of Lie algebras $A:\g_o\to\g_{\bar o}$.
	Remind that $A:\g_o\to\g_{\bar o}$ is a morphism of Lie algebras if and only if $A:\G_o\to\G_{\bar o}$ is a morphism of Lie groups, under the identification $\g=\G$ as in Section \ref{subsec:tangent cone}.

	Set $\X:=\{tX_i^{o,0}:i\in\{1,\dots,r\},t\in\R\}$.
	Remind that $\X$ generates $\G_o$ as a Lie group.

	Define $A:\X\to \g_{\bar o}$ as
	$$A[tX_k^{o,0}] := t\sum_{j=1}^{\bar r} h_{kj}(o,0) \bar X_j^{o,0}\in\g_{\bar o}.$$
	
	Let $\epsilon_k\to0$ be a sequence such that $f^{o,\epsilon_k}$ converge to some $g$ as in Lemma~\ref{110}.
    By Proposition \ref{prop1044}, for almost every point $o\in E$, we have a neighborhood $K$ of the point $o$ and a constant $t_0$ such that
    \begin{equation}\label{eq2140}
    \exp(A[tX_j^{o,0}])(g(p)) = g(\exp(tX_j^{o,0})(p)),
    \end{equation}
    for $p\in K$, $1\leq j\leq r$ and $|t|\leq t_0$.
      We identify $\g_o$ with $\G_o$ in the standard way. Let $V\in\g_o$. By Lemma 1.40 in \cite{SSS}, there are $\xi_1,\dots,\xi_\ell\in\X$ such that $V=\xi_1*\dots*\xi_\ell$, where $\ell$ is bounded by a constant and $*$ is the multiple in $\G_o$.
    Moreover, we have for all $p\in K$
    \[%\begin{multline*}
     	\exp(V)(p)
    	= \exp(\xi_1*\dots*\xi_\ell)(p)
    	= \exp(\xi_1)\circ\dots\circ\exp(\xi_\ell)(p).
    \]%\end{multline*}
    Set $\xi_i=t_iX_{j_i}^{o,0}.$
    Again, by Lemma 1.40 in \cite{SSS}, there is a constant $c_0$ such that if $\|V\|\leq c_0$, $|t_i|\leq t_0$ for any $1\leq i \leq\ell$.  
Furthermore, there is a smaller neighborhood $K^{\prime}$ of $o$ and a  constant $\epsilon_0^{\prime}$, such that for $p\in K^{\prime}$ and $\|V\|\leq \epsilon_0^{\prime}$ 
    \[\exp(\xi_j)\circ\dots\circ\exp(\xi_\ell)(p)\in K,\] for all $2\leq j\leq \ell$.   
In the following, fix $p\in K^{\prime}$ and $\|V\|\leq \min(\epsilon_0^{\prime},c_0,\epsilon_0)$.   
  Therefore, iterating \eqref{eq2140}, we have 
    \begin{equation}\label{eq2142}
    g\big(\exp(V)(p)\big) %= \exp(A\xi_1)\circ\dots\circ\exp(A\xi_\ell)\big(g(p)\big)
    = \exp(A\xi_1\bar*\dots\bar*A\xi_\ell)(g(p)),
    \end{equation}
    where $\bar *$ is the multiple in $\G_{\bar o}$.
	%The equation \eqref{eq2142} have many consequences.
If we take $p=o$, then \eqref{eq2142} becomes
    \begin{equation}\label{220}
    g(\exp(V)(o))
    = \exp(A\xi_1\bar*\dots\bar*A\xi_\ell)(\bar o)
    \end{equation}
    Since the right hand side does not depend on the sequence $\epsilon_k\to0$, $g$ is unique.

    On the other side, now $A$ can be extended to a map $\g_o\to\g_{\bar o}$.
    Indeed, if $V=\xi_1*\dots*\xi_\ell =\xi'_1*\dots*\xi'_{\ell'}$, then by \eqref{eq2142} the vector fields $A\xi_1\bar*\dots\bar*A\xi_\ell$ and $A\xi'_1\bar*\dots\bar*A\xi'_{\ell'}$ have the same action on the image of $g$, in particular they are same at $\bar o$, hence they are same everywhere.

    By construction, this extension of $A$ is a Lie algebra morphism and by Lemma \ref{lem2002} it is the $P$-differential of $f$ at $o$.

\end{proof}

We are now ready to prove Stepanov's Theorem:
\begin{proof}[Proof of Theorem \ref{thm:Differentiability Stepanov}]
It is easily to deduce that $L(f)=\cup_lF_l$, where \[F_l=\{x: \bar d(f(z),f(x))\leq ld(z,x),z\in B(x,1/l)\},\]
where $l\geq 1$. Next, we decompose the set $F_l$ into countable union $E_{l,i}$ such that $f|_{E_{l,i}}$ is $i$-Lipschitz, where integer $i\geq 1$. Moreover,  by the previous Rademacher type result, $f|_{E_{l,i}}$ is $ P$-differentiable at almost every point of $E_{l,i}$. Choose a point $x_0$ from them. According to the proof of the previous Rademacher type result, $x_0$ is a density point of $E_i$. Thus, as Lemma 3.15 in \cite{f69} and $f|_{E_{l,i}}$ is $ P$-differentiable at $x_0$, $f$ is also $ P$-differentiable at $x_0$. As $L(f)=\bigcup E_{l,i}$, we obtain that $f$ is $P$-differentiable at almost every point in $M$.
\end{proof}

\appendix
\section{Pansu differentiability in Carnot groups}
For the convenience of the reader, we prove Stepanov's theorem for mappings between Carnot groups, which serves as a special case of equiregular subRiemannian manifolds; see Section \ref{subsec:Carnot group} for a brief introduction of Carnot groups. 

Let $G$ and $\bar G$ be two Carnot groups. We have on $\bar G$ all the same objects as on $G$, and we distinguish them by putting a bar on the ones for $\bar G$.
Let $f:G\to \bar G$ be a Borel mapping, $o\in G$ and $\bar o:=f(o)\in\bar G$. Note that the exponential mapping on $G$ satisfies 
\[\exp(X_i)(g)=gX_i(0),\]
for any $1\leq i\leq r$ and any $g\in G$.

\begin{definition}[$P$-Differential on Carnot group]\label{aadef1738}
	We say that $f:G\to \bar G$ is $P$-differentiable at $o$ if there exists a morphism of graded Lie algebras $A:\g_o\to\g_{\bar o}$ such that
	\[
	\lim_{\g_o\ni X\to 0} \frac{\bar d\left(\exp(A[X])(\bar o),f(\exp(X)(o))\right)}{\|X\|}=0
	\]
	where $\|\cdot\|$ is any homogeneous norm on $\g_o$. When $f$ is $P$-differentiable at $o$, we write $Df(o)$ instead of $A$ for the $P$-differential.
\end{definition}

The aim of this appendix is to prove the following Stepanov's Theorem.

\begin{theorem}\label{thm:stepanov Carnot group}
	Let $f:(G,d)\to (\bar G,\bar d)$ be a Borel mapping between two Carnot groups.
	Then $f$ is $P$-differentiable for almost every~$o$ in the set
	\[
	L(f) := \left\{o\in M: \limsup_{p\to o} \frac{\bar d(f(o),f(p))}{d(o,p)} < \infty \right\}.
	\]
\end{theorem}

\subsection{A variant of the Lebesgue differentiation theorem}
The aim of this section is to show the following differentiation theorem, which will be used in our later proof of the Stepanov's theorem in Carnot groups. For every $p\in G$, let $B_p\subset U^1_p$ be a compact neighborhood of $p$. For a subset $A\subset G$, $|A|$ denotes the standard Lebesgue measure of the set $A$. 

\begin{proposition}\label{aaprop:Lebesg diff them}
Let $\Omega\subset G\times[0,1]$ be an open neighborhood of $G\times\{0\}$ and let $\phi:\Omega\to G$, $(p,t)\mapsto \phi_tp$, be the flow of a smooth vector field on $G$ which is nonzero everywhere. If $h:G\to\R$ be a locally bounded function, then, for almost all $o\in G$
	\[
	\lim_{\epsilon\to0}
	\dashint_{\delta_\epsilon^o(B_o)}\dashint_0^\epsilon |h(\phi_sp)-h(o)|\dd s\dd p = 0
	\]
\end{proposition}

We will use a version of the Lebesgue differentiation theorem due to Federer~\cite[Theorem 2.9.8, Page 156--165]{f69}. 

%%%%%%%%%%%%%%%%%%%%%%%%%%%%%%%%%%%%%%%%%%%%%%%%%%%%%%%%%%%%%%%
\begin{lemma}\label{aalem1356}
 	Define
	\begin{equation}\label{aaeq1343}
	 	\scr V:= \left\{\big(p,\delta^p_\epsilon B_p\big) : p\in G,\ \epsilon\in(0,1]\right\} .
	\end{equation}
	The family $\scr V$ is a Vitali relation,  in the sense of Federer~\cite[\S 2.8.16]{f69}.
\end{lemma}
\begin{proof}
We do this using~\cite[Theorem 2.8.17]{f69}.
Using Federer's notation, in our case we choose
\[
\tau=2\qquad\text{ and }\qquad\delta\big(p,\delta^p_\epsilon B_p\big) = \diam_d(\delta^p_\epsilon B_p).
\]
We need only to show that for almost all $o\in M$:
\begin{equation}
 	\limsup_{\epsilon\to0} \frac{|\widehat{\delta^o_\epsilon B_o}|}{|\delta^o_\epsilon B_o|} <+\infty,
\end{equation}
where
\[
\widehat{\delta^o_\epsilon B_o} = \bigcup \delta^p_\eta B_p,
\]
the union is taken on all $\delta^p_\eta B_p$ such that $\delta^p_\eta B_p\cap\delta^o_\epsilon B_o \neq\emptyset$ and $\diam_d(\delta^p_\eta B_p) \le 2 \diam_d(\delta^o_\epsilon B_o)$.
Hence
\[
\widehat{\delta^o_\epsilon B_o} \subset B_d(o,3\diam_d(\delta^o_\epsilon B_o))
\]
and we have to prove
\begin{equation}\label{aaeq1516}
 	\limsup_{\epsilon\to0} \frac{|B_d\big(o,3\diam_d(\delta^o_\epsilon B_o)\big)|}{|\delta^o_\epsilon B_o|} <+\infty .
\end{equation}

First we claim that
\[
\diam_d(\delta^o_\epsilon B_o) = O( \epsilon)
\]
Indeed, if $x\in B_o$, then
\[
d(o,\delta^o_\epsilon x) 
= \epsilon d(o,x)
\]
and therefore $d(o,\delta^o_\epsilon x) = O(\epsilon)$, uniformly in $x$.

Since the Legesgue measure is left-invariant and  $\delta^0_r$ has determinant equal to $r^Q$, we have
\[
|B_d(o,r)|=|\delta^o_r(B_d(o,1))|=|L_o\delta^0_rL_{o^{-1}}(B_d(o,1))| \sim r^Q
\]
where $ Q$ is the homogeneous dimension of $G$.

Similarly, we have
\[
|\delta^o_\epsilon B_o|=|L_o\delta^0_\epsilon L_{o^{-1}}(B_o)| \sim \epsilon^Q,
\]
which leads to \eqref{aaeq1516}. 	
\end{proof}

\begin{remark}\label{rmk:on variant of Leb in Carnot}
Compared with the case of equiregular subRiemannian manifolds, we only use the fundamental property that the Lebesgue measure is left-invariant. In the case of subRiemannian manifolds, we have to use the smoothness of Popp measure.
\end{remark}

%%%%%%%%%%%%%%%%%%%%%%%%%%%%%%%%%%%%%%%%%%%%%%%%%%%%%%%%%%%%%%%
\begin{lemma}\label{aalemma:33}
 	Define
	\[
	R_\epsilon(p):=\dashint_0^\epsilon |h(\phi_sp)-h(p)|\dd s .
	\]
	Then for almost every $p\in G$ we have
	\begin{equation}\label{aaeq1257}
	\lim_{\epsilon\to0} R_\epsilon(p) = 0 .
	\end{equation}
\end{lemma}
\begin{proof}
	Applying the Lebesgue differentiation theorem to the function $t\mapsto h(\phi_tq)$ for any $q\in G$, we obtain that for almost all $(q,t)$% for a.e. $t\in[0,1]$
	$$\lim_{\epsilon\to0}\dashint_0^\epsilon |h(\phi_{t+s}q)-h(\phi_tq)| \dd s = 0,$$
which implies \eqref{aaeq1257} holds for $p=\phi_tq$. Since the map $(q,t)\mapsto\phi_t q$ is locally Lipschitz and surjective, it maps a set of full measure into a set of full measure, therefore for almost every $p\in M$ \eqref{aaeq1257} holds.
\end{proof}
%%%%%%%%%%%%%%%%%%%%%%%%%%%%%%%%%%%%%%%%%%%%%%%%%%%%%%%%%%%%%%%
\begin{lemma}\label{aalem1321}
 	Let $F\subset G$ be a measurable subset.
	Then for almost all $o\in F$
	\begin{equation}\label{aaeq1321}
	\lim_{\epsilon\to0}\frac{|\delta^o_\epsilon B_o\setminus F|}{|\delta^o_\epsilon B_o|} = 0.
	\end{equation}
\end{lemma}
\begin{proof}
 	Since $\scr V$, defined in \eqref{aaeq1343}, is a Vitali relation, the claim follows by applying the Lebesgue differentiation theorem to  the characteristic function $\chi_F$ of $F$.
\end{proof}
%%%%%%%%%%%%%%%%%%%%%%%%%%%%%%%%%%%%%%%%%%%%%%%%%%%%%%%%%%%%%%%
\begin{proof}[Proof of Proposition~\ref{aaprop:Lebesg diff them}]
%\begin{lemma}\label{lemma:35}
We need to show that for almost all $o\in G$
 	\begin{equation}\label{aaeq1323}
	 	\lim_{\epsilon\to0}\dashint_{\delta^o_\epsilon B_o}\dashint_0^\epsilon |h(\phi_sp)-h(o)|\dd s\dd p = 0
	\end{equation}
%\end{lemma}

Thanks to Egorov Theorem and Lemma~\ref{aalemma:33}, for every $\eta>0$ there is a measurable subset $F\subset G$ such that $|G\setminus F| \le \eta$ and $R_\epsilon$ converge uniformly to $0$ on $F$. Since $\eta$ is arbitrary, it suffices to show that \eqref{aaeq1323} holds for almost all $o\in F$.
	
	Since $\scr V$, defined in \eqref{aaeq1343}, is a Vitali family, by Lemma \ref{aalem1321}, we deduce that for almost every $o\in F$  \eqref{aaeq1321} holds and
	\begin{equation}\label{aaeq1327}
	 	\lim_{\epsilon\to0}\dashint_{\delta^o_\epsilon B_o} |h(p)-h(o)| \dd p = 0 .
	\end{equation}
	
	For such an $o$ we have
 	\begin{multline*}
	\dashint_{\delta^o_\epsilon B_o}\dashint_0^\epsilon |h(\phi_sp)-h(o)|\dd s\dd p
	\\ \le
	\underbrace{\dashint_{\delta^o_\epsilon B_o}\dashint_0^\epsilon |h(\phi_sp)-h(p)|\dd s\dd p}_{A}
	+
	\underbrace{\dashint_{\delta^o_\epsilon B_o} |h(p)-h(o)| \dd p}_{B}\\
	\end{multline*}
	where part $B$ converges to $0$ as $\epsilon\to0$ because of \eqref{aaeq1327}.
	For part $A$, we have
	\begin{align*}
	 	A&=\dashint_{\delta^o_\epsilon B_o} R_\epsilon(p) \dd p
		= \frac1{|\delta^o_\epsilon B_o|}\int_{\delta^o_\epsilon B_o\cap F}R_\epsilon(p) \dd p
			+ \frac1{|\delta^o_\epsilon B_o|}\int_{\delta^o_\epsilon B_o\setminus F}R_\epsilon(p) \dd p \\
		&\le \dashint_{\delta^o_\epsilon B_o} R_\epsilon(p) \chi_F \dd p + C \frac{|\delta^o_\epsilon B_o\setminus F|}{|\delta^o_\epsilon B_o|},
	\end{align*}
	where $C>0$ is some constant that bounds $R_\epsilon$, which exists because $h$ is locally bounded. It is now clear that $A$ converges to $0$ as $\epsilon\to 0$.
\end{proof}

\subsection{Blow-up of Lipschitz functions}
Let $E\subset G$ be closed and let $f:E\to \bar G$ be an $L$-Lipschitz function. For every $p\in G$ let $B_p\subset G$ be a closed $d$-ball centered at $p$ such that the dilations $\delta^p_\epsilon$ are well defined for $\epsilon\in(0,1]$.
For $o\in E$, we define the functions
%\[
%f^{o,\epsilon}:=\delta^{f(o)}_{\frac1\epsilon}\circ f\circ\delta^o_\epsilon
%\]
\[
f^{o,\epsilon}:=\bar{\delta}_{\frac1\epsilon}\circ f\circ\delta^o_\epsilon,
\]
where $\bar{\delta}_{\frac1\epsilon}$ is with respect to $f(o)$.

Define with these $B_p$
\[
 	\scr V:= \left\{\big(p,\delta^p_\epsilon B_p\big) : p\in M,\ \epsilon\in(0,1]\right\} .
\]
Lemma \ref{aalem1356} implies that $\scr V$ is a Vitali relation, and hence almost every~$o\in E$ is a $\scr V$-density point of $E$.

Our main step is to show a sort of equicontinuity of $f^{o,\epsilon}$ at almost every point of $E$.

\begin{lemma}\label{aalem1350}
 	For every $o\in E$, setting $\bar o=f(o)\in \bar G$, there is an open neighborhood of $o$ $U^o$ such that $U^o\subset U^o_1$ and $f^{o,\epsilon}$ is a well-defined map $U^o\cap \delta^{o}_{\frac1\epsilon}E\rightarrow U_1^{\bar o}$.
\end{lemma}
\begin{proof}
Fix a point $o\in E$ and choose $r>0$ such that \[B_{\bar{d}}(\bar {o},2r)\subset U_1^{\bar {o}}.\]
Then, it is easy to deduce that
\[
B_{\bar{d}}(\bar{o}, \epsilon r)=\bar{\delta}_\epsilon\left(B_{\bar{d}}(\bar{o}, r)\right)  \subset \bar{\delta}_{\epsilon}U_{1}^{\bar{o}}.\]

Thus,
\[
\bar{d}\left(\bar{o}, f\left({\delta}_\epsilon p\right)\right) \leq L d\left(o, \delta_\epsilon p\right) \leq L \epsilon d(o, p),
\]
where $L$ is the Lipschitz constant of $f$.
Then, as long as $p\in U^o\cap {\delta}_{\frac{1}{\epsilon}}^oE$ satisfying that
$$
d(o, p) \leq \frac{r}{L} ,
$$
$\bar{\delta}_{\frac{1}{\epsilon}}\left(f\left(\delta_\epsilon p\right)\right)$ is well-defined and belongs to $U_1^{\bar o}$.
\end{proof}

Let $p\in E$ be a density point of $E$. We introduce the following type of convergence, adapted to the fact that $f^{p,\epsilon}$ is not defined in a neighborhood of $p$. Set $\bar p$ to be $f(p)$.
\begin{definition}\label{aa120}
Let $\epsilon_k\to 0$ be a sequence and $g$ be a continuous function on $B_p$. We say that $f^{p,\epsilon_k}\to g$ uniformly on $B_p$ if \[\sup\big\{\bar d(f^{p,\epsilon_k}(q),g(q)):q\in\delta^p_{\frac1{\epsilon_k}}E\cap B_p\big\}\to 0.\]
\end{definition}
We denote the Hausdorff distance under the metric $d$ to be $d_H$.

\begin{lemma}\label{aa130}
 	If $o\in E$ is a $\scr V$-density point of $E$, then as $\epsilon_k\to 0$, $d_H(\delta^o_{\frac1{\epsilon_k}}E\cap B_o,B_o)\to 0$.
\end{lemma}
\begin{proof}
 	Assume that the conclusion were false.
	Then, up to passing to a subsequence, there is $a>0$ such that for every $k\in\N$ there is $x_k\in B_o$ with
	\[
	B_o\cap B(x_k,a) \subset B_o\setminus \delta^o_{\frac1{\epsilon_k}}E.
	\]
Then, as the Lebesgue measure is left-invariant, 
	\[
	\lim_{k\to\infty} \frac{|\delta^o_{\epsilon_k}B_o\setminus E|}{|\delta^o_{\epsilon_k}B_o|}
	= \lim_{k\to\infty} \frac{ |B_o\setminus \delta^o_{\frac1{\epsilon_k}}E| }{|B_o|}
	\ge \lim_{k\to\infty} \frac{|B_o\cap B(x_k,a)|}{|B_o|} >0,
	\]
	i.e., $o$ is not a $\scr V$-density point of $E$.
\end{proof}

%%%%%%%%%%%%%%%%%%%%%%%%%%%%%%%%%%%%%%%%%%%%%%%%%%%%%%%%%%%%%%%
	Set $\bar o:=f(o)$. The family of functions $\{f^{o,\epsilon}\}_{\epsilon\in(0,1]}$ is called \emph{equicontinuous} if
\begin{equation}\label{aaeq0047}
	\begin{array}{c}
	 	\forall \eta>0,\ \exists\delta>0,\ \forall\epsilon,\ \forall p,q\in \delta^o_{\frac1{\epsilon}}E\cap B_o \\
		d(p,q)\le \delta \quad\THEN\quad \bar {d}(f_\epsilon(p),f_\epsilon(q)) \le \eta.
	\end{array}
\end{equation}

\begin{lemma}[Equicontinuity of $f^{o,\epsilon}$]\label{aalem1403}
 	Let $o\in E$ be a $\scr V$-density point of $E$.
	Then the family $\{f^{o,\epsilon}\}_{\epsilon}$ is equicontinuous.
\end{lemma}
\begin{proof}
Since $f$ is $L$-Lipschitz, it follows easily that
	\begin{align*}
	\bar d(f^{{ {o}},\epsilon}(p),f^{o,\epsilon}(q)) &=\bar{d}(\bar{\delta}_{\frac{1}{\epsilon}}f\delta_{\epsilon}(p),\bar{\delta}_{\frac{1}{\epsilon}}f\delta_{\epsilon}(q)) =\epsilon^{-1}\bar{d}(f\delta_{\epsilon}(p),f\delta_{\epsilon}(q)) \\
&\le \epsilon^{-1}Ld(\delta_{\epsilon}(p),\delta_{\epsilon}(q))= Ld(p,q).
	\end{align*}

\end{proof}

\begin{lemma}\label{aa110}
Let $o\in E$ be a $\scr V$-density point of $E$.
Then the family $\{f^{o,\epsilon}\}_{\epsilon}$ has a subsequence uniformly convergent to a continuous function $g$ on $B_o$ in the sense of Definition~\ref{aa120}.
\end{lemma}
\begin{proof}
Basing on Lemma~\ref{aa130} and Lemma~\ref{aalem1403}, using \cite[Lemma 8.6 and Lemma 8.7]{GGSS}, we conclude that for a sequence $\epsilon_k\to 0$, there exist a continuous function $g$ on $B_o$ and a subsequence of $\{f^{o,\epsilon_k}\}$ converging to $g$ uniformly. 
\end{proof}

\begin{remark}\label{rmk:equicontinuity Carnot gropu}
Compared with the case of subRimannian manifolds, the proof here is greatly simplified. In the subRiemannian case, we need to use the more subtle eventually equicontinuous property.
\end{remark}

Finally, we introduce the following lemma that will be used later.
\begin{lemma}\label{aalem2002}
Let $o$ be a density point of $E$. Suppose that $f^{o,\epsilon}\to f^{o,0}$ uniformly in the sense of Definition~\ref{aa120} and that there is a morphism of Lie algebras $A:\g_o\to\g_{\bar o}$ and a constant $\epsilon_0>0$ such that
\begin{equation}\label{aa140}
f^{o,0}(\exp(X)(o)) = \exp(A[X])(\bar o)
\end{equation}
holds for all $\|X\|\le \epsilon_0$. Then $f$ is $P$-differentiable at $o$ and $Df(o)=A$.
\end{lemma}
Note that the limit in $ P$-differentiability in the above lemma is with respect to only $X\in g_o$ such that $\exp (X)\in E$.
\begin{proof}
We use the standard exponential coordinates. As $f^{o,0}$ is the uniform limit of $f^{o,\epsilon}$ and $f^{o,\tau}\delta_{\epsilon}=\bar{\delta}_{\epsilon}f^{o,\tau\epsilon}$, we have \begin{equation}\label{aa131}
f^{o,0}(\exp(\delta_\epsilon X)(o))=\bar{\delta}_\epsilon f^{o,0}(\exp(X)(o)).
\end{equation}
Then, by \eqref{aa140} and the above equation
\begin{equation*}
\exp(A[\delta_\epsilon X])(\bar o) =f^{o,0}(\exp(\delta_\epsilon X)(o))=\bar{\delta}_\epsilon f^{o,0}(\exp(X)(o))=\bar{\delta}_\epsilon \exp(A[X])(\bar o).
\end{equation*}
Then $A[\delta_\epsilon X]=\bar{\delta}_\epsilon A[X]$ and thus, $A$ is a morphism of graded Lie algebra. 

For $\epsilon\leq \epsilon_0$, any vector of norm $\epsilon$ can be written as $\delta_\epsilon X$ with $\|X\|=1$.
	Moreover, by \eqref{aa140} and \eqref{aa131}, we have 
	\begin{align*}
	\frac{ \bar d(\exp(A[\delta_\epsilon X])(\bar o) , f(\exp(\delta_\epsilon X)(o)) }{\|\delta_\epsilon X\|} &=
	 \epsilon^{-1} \bar d(f^{o,0}(\exp(\delta_\epsilon X)(o)), f(\exp(\delta_\epsilon X)(o))) \\
	 &= \epsilon^{-1} \bar d(\bar\delta_\epsilon f^{o,0}(\exp(X)(o)), \bar\delta_\epsilon \bar\delta_{\frac1\epsilon} f(\delta_\epsilon \exp( X)(o)))  \\
	 &=  \bar d ( f^{o,0}(\exp(X)(o)), f^{o,\epsilon} (\exp(X)(o))).
	\end{align*}
Since $f^{o,\epsilon}$ converges uniformly, we get
	\[
	\lim_{\epsilon\to0} \bar d( f^{o,0}(\exp(X)(o)), f^{o,\epsilon} (\exp(X)(o))) = 0
%	\bar d_0^o ( f^{o,0}(\exp(X)(o)), f^{o,0} (\exp(X)(o)))
	\]
	and the limit is uniform with respect to $X$.
\end{proof}

%%%%%%%%%%%%%%%%%%%%%%%%%%%%%%%%%%%%%%%%%%%%%%%%%%%%%%%%%%%%%%%
%%%%%%%%%%%%%%%%%%%%%%%%%%%%%%%%%%%%%%%%%%%%%%%%%%%%%%%%%%%%%%%
\subsection{Blow-up of horizontal vector fields}
Fix $1\leq j\leq r$. Let $W=X_j$ and let $(p,t)\mapsto \phi_tp$ be its flow.

%We extend $f$ outside of $E$ in the following way.

Set \begin{equation}\label{aababbb}
F(p,t) := f(\phi_tp).
\end{equation}
Then $F$ is well-defined on $\{(p,t):\phi_tp\in E\}$.
Since
$$\bar d(f(\phi_tp),f(\phi_sq)) \le L d(\phi_tp,\phi_sq),$$
$F$ is locally Lipschitz, where the Lipschitz constant depends only on the Lipschitz constant of $f$ and the Lipschitz constant of $\phi$ on a compact set.

\subsubsection{Extension of $F$ on $G\times\R$}
We seek an extension of $F$ (still denoted by $F$) on $G\times\R$ satisfying properties:
\begin{enumerate}
\item[(P1)] $F(p,t+s) = F(\phi_tp,s)$ for all $p\in G$ and all $s,t\in\R$ such that $\phi_tp$ exists.
\item[(P2)] For each $p\in G$ the curve $t\mapsto F(p,t)$ is locally Lipschitz.
\end{enumerate}

We first extend $F$ on $E\times\R$ in such a way that each curve $t\mapsto F(p,t)$ is locally Lipschitz.
More precisely: for $p\in E$ define $I_p:=\{t\in\R:\phi_tp\in E\}\subset\R$.
Since $E$ is closed, $I_p$ is closed as well.
Let $\hat t\in\R\setminus I_p$. Then there are two cases. In the first case, it happens that there is $t_1\in I_p$ such that  $\hat t\in(-\infty, t_1)\subset\R\setminus I_p$ or $t\in ( t_1,+\infty)\subset\R\setminus I_p$.
Then we set $F(p,t) = F(p,t_1)$.
In the second case, there are $t_1,t_2\in I_p$ with $\hat t\in (t_1,t_2)\subset\R\setminus I_p$.
Then
\[
\bar d(f(\phi_{t_1}p),f(\phi_{t_2}p)) \le Ld(\phi_{t_1}p,\phi_{t_2}p) \le  \tilde L |t_1-t_2|.
\]
Therefore there is a geodesic $\gamma:[t_1,t_2]\to\bar G$ joining $f(\phi_{t_1}p)$ to $f(\phi_{t_2}p)$ with constant velocity, i.e. $\tilde L$-Lipschitz.
In this case, we define $F(p,t)=\gamma(t)$ for $t\in(t_1,t_2)$, where the curve is chosen in such a manner that (P1) holds for all $p,\phi_t p\in E$.

Moreover, on the set $E^{\prime}=\{p:\text{there exists } t \text{ such that }\phi_tp\in E\}$, we may extend $F$ using the rule $F(p,t+s) = F(\phi_tp,s)$. As (P1) holds for all $p,\phi_t p\in E$, the above definition is independent of the choice of $t$. For $p$ outside of this set, we simply define $F(p,t)=\bar p$ for some fixed point $\bar p\in\bar G$. 

Next, we verify the property (P1). In the case $p\in E$, if $\phi_t p\in E$, then according to our construction we already obtain the property (P1); if $\phi_t p\not\in E$, the construction on $E^{\prime}$ tells us (P1) holds as well. When $p\in E^{\prime}$, we have $\phi_tp\in E^{\prime}$. Then, $F(p,t+s) = F(\phi_{t^{\prime}}p,t+s-t^{\prime})$, where $\phi_{t^{\prime}}p\in E$. Thus, $F(\phi_tp,s) = F(\phi_{t^{\prime}}p,s-t^{\prime}+t)=F(p,t+s)$. Therefore, (P1) holds. In the case $p\not\in E^{\prime}$, (P1) is true as $F(p,t)$ is constant. 

Finally, we verify the property (P2). In the case $p\in E$, our choice of curves tells us the property holds. When $p\in E^{\prime}$, according to the rule $F(p,t+s) = F(\phi_tp,s)$, where $\phi_tp\in E$, we have the local Lipschitz property from the fact that the flow is local Lipschitz. For $p$ outside these two sets, as $F(p,t)$ is constant, (P2) holds.

\begin{remark}\label{aa210}
The construction tells us that the Lipschitz constant in the condition (P2) is uniform with respect to $(p,t)$ in a compact set.
\end{remark}
%%%%%%%%%%%%%%%%%%%%%%%%%%%%%%%%%%%%%%%%%%%%%%%%%%%%%%%%%%%%%%%
\subsubsection{Blow-up of $F$}
According to (P1) and Proposition 3.50 in \cite{cambridge}, there are $h_j:G\times\R\to\R$, $j\in\{1,\dots,\bar r\}$ such that for each $p\in G$,
\[
\frac{\de F}{\de t} (p,t) = \sum_{j=1}^{\bar r} h_j(p,t) \bar X_j(F(p,t))
\]
holds for almost every $t$.
Notice that, by Remark~\ref{aa210}, $h_j$ are locally bounded.
Furthermore, by (P1), we have that
\[\frac{\de F}{\de (t+s)}(p,t+s)=\frac{\de F}{\de s}(\phi_t p,s)\]
holds for almost every $r+s$. By setting $s=0$ and $t=t_0$,
\[\frac{\de F}{\de t}(p,t_0)=\frac{\de F}{\de t}(\phi_{t_0} p,0)\]
holds for almost every $t_0$. Then, using the fact that the flow $\phi_tp$ is locally Lipschitz and surjective, we have that for almost every $p\in G$ the derivative $\frac{\de F}{\de t} (p,0)$ exists.
Thus, \[\frac{\de F}{\de t} (p,0)= \sum_{j=1}^{\bar r} h_j(p,0) \bar X_j(F(p,0))\]
holds for almost all $p$.

The aim of this section is to blow-up both $G$ and $\bar G$ keeping track of the map $F$.
The result is, in some sense, the flow of a left-invariant vector field on $\bar G$; See Proposition \ref{aaprop1044}.

For $o\in G$ and $\epsilon\in(0,1]$ we have the vector fields
\[
W^{o,\epsilon} := \epsilon\cdot\dd\delta^o_{\frac1\epsilon}\circ W\circ\delta^o_\epsilon.
\]
By Lemma~\ref{ttt}, $W^{o,\epsilon}=W$.

For $o\in M$, $\bar o:=F(o,0)$ and $\epsilon>0$, set
\[
F^{o,\epsilon}(p,t) := \bar\delta^{\bar o}_{\frac1\epsilon}F(\delta^o_\epsilon p,\epsilon t) .
\]
Notice that
$
F^{o,\epsilon}(p,0) = f^{o,\epsilon}(p)
$
if $\delta^o_{\epsilon}p\in E$.

%{\color{red} The stetement of the next lemma is not clear!}
%{\color{blue} I changed it:}
\begin{lemma}\label{aalemma:51}
	For all $o\in E$ and all $t\in\R$
	\begin{equation}\label{aaeq0036}
	F^{o,\epsilon}(p,t) = f^{o,\epsilon}(\phi_tp) ,
	\end{equation}
	if the right-hand side is well-defined.
	Moreover, for each $p$,
	\begin{equation}\label{aaeq2349}
	 	\frac{\de F^{o,\epsilon}}{\de t}(p,t) = \sum_{j=1}^{\bar r} h_j(\delta_\epsilon^op,\epsilon t) \bar X_j(F^{o,\epsilon}(p,t))
	\end{equation}
holds for almost every $t$.
\end{lemma}

\begin{proof}
	Fix $p$ and set $\gamma(t)=\phi_tp$.
	Then $\gamma(0)=p$ and $\gamma'(t) = W(\gamma(t))$.
	Define $\eta_\epsilon(t)=\delta_{\frac1\epsilon}(\gamma(\epsilon t))$.
	Then $\eta_\epsilon(0)=\delta_{\frac1\epsilon}(p)$ and
	\begin{equation*}
	 	\eta_\epsilon' (t)
		= \dd\delta^o_{\frac1\epsilon}[\epsilon\gamma'(\epsilon t)]
		= \epsilon \dd\delta^o_{\frac1\epsilon}[W(\gamma(\epsilon t))] = \epsilon \dd\delta^o_{\frac1\epsilon}[W(\delta^o_\epsilon\delta^o_{\frac1\epsilon}\gamma(\epsilon t))]
		= W(\eta_\epsilon(t))
	\end{equation*}
	i.e. $\eta_\epsilon(t) = \exp(tW)(\delta_{\frac1\epsilon}(p))$.
	In other words
	$%\begin{equation}
	\delta_{\frac1\epsilon}(\phi_{\epsilon t}p) = \phi_t(\delta_{\frac1\epsilon}p) .
	$
	Hence, if $\phi_{\epsilon t}\delta^o_\epsilon p\in E$, then by \eqref{aababbb},
	\begin{align*}
	F^{o,\epsilon}(p,t)
	&=\bar{\delta}^{\bar o}_{\frac1\epsilon} F(\delta^o_\epsilon p,\epsilon t)
	=\bar{\delta}^{\bar o}_{\frac1\epsilon} \circ f\circ \phi_{\epsilon t}\delta^o_\epsilon p \\
	&=\bar{\delta}^{\bar o}_{\frac1\epsilon} \circ f\circ \delta^o_\epsilon\circ\delta^o_{\frac1\epsilon}\circ \phi_{\epsilon t}\delta^o_\epsilon p
	= f^{o,\epsilon}(\phi_tp),
	\end{align*}
which gives \eqref{aaeq0036}. Regarding \eqref{aaeq2349}, using Lemma~\ref{ttt}, we have
	\begin{align*}
	 	% \sum_{j=1}^{\bar r} h_j^{o,\epsilon}(p,t) \bar X_j^{o,\epsilon}(F^{o,\epsilon}(p,t))
		 \frac{\de F^{o,\epsilon}}{\de t}(p,t) %= \\
		 &= \frac{\de}{\de t} \bar \delta_{\frac1\epsilon}^{\bar o} F(\delta^o_\epsilon p, \epsilon t)
		 = \epsilon \dd\bar \delta_{\frac1\epsilon}^{\bar o} \frac{\de F}{\de t}\Big(\delta^o_\epsilon p, \epsilon t\Big)  \\
		 &= \epsilon \dd\bar \delta_{\frac1\epsilon}^{\bar o} \Big[\sum_{j=1}^{\bar r} h_j(\delta^o_\epsilon p,\epsilon t) \bar X_j(F(\delta^o_\epsilon p,\epsilon t))\Big] \\
		& = \sum_{j=1}^{\bar r} h_j(\delta^o_\epsilon p,\epsilon t) \epsilon \dd\bar \delta_{\frac1\epsilon}^{\bar o} \bar X_j(\bar\delta^{\bar o}_\epsilon \bar\delta^{\bar o}_{\frac1\epsilon}F(\delta^o_\epsilon p,\epsilon t)) \\
		 &= \sum_{j=1}^{\bar r} h_j(\delta^o_\epsilon p,\epsilon t) \bar X_j(F^{o,\epsilon}(p,t)) .
	\end{align*}
\end{proof}

\begin{proposition}\label{aaprop1044}
 	Let $o\in E$ be a $\scr V$-density point of $E$. Let $\epsilon_k\to 0$ be a sequence such that $f^{o,\epsilon_k}$ converge uniformly to a continuous function $g:B_o\to \bar M$.
	Then there exists a neighborhood $K$ of the point $o$, such that for $t$ sufficiently small and $p\in K$,
	\[
	\exp\Big(t\sum_{j=1}^{\bar r} h_j(o,0)\bar X_j\Big)\big(g(p)\big)
	=
	g\Big(\exp\big(tX_j\big)(p)\Big),
	\]
holds for almost every density point $o\in E$.
\end{proposition}

Since our notation is getting heavier and heavier, we will drop the subscript $k$ in $\epsilon_k$ and write just $\epsilon$.
\begin{proof}%[Proof of Proposition \ref{prop1044}]
	Define
	\[
	G(p,t) =g\Big(\exp\big(tX_j\big)(p)\Big).
	\]
	Now considering everything in the Euclidean coordinates, where the Euclidean distance in this coordinate is denoted by $|\cdot |$, the curve
	\[
	\gamma_p(t) = g(p) + \int_0^t \sum_{j=1}^{\bar r} h_j(o,0) \bar X_j(G(p,s)) \dd s
	\]
	is well-defined.
	Notice that $\gamma_p(0)=g(p)$ and $\gamma_p'(t)=\sum_{j=1}^{\bar r} h_j(o,0) \bar X_j(G(p,t))$.
	
In order to obtain the conclusion, it suffices to prove
	\[
	\gamma_p(t) = G(p,t),
	\]
for $t$ sufficiently small.
For this, as $g$ is continuous, we only need to prove that there exists a neighborhood $K$ of $o$ such that for $t$ sufficiently small,
	\begin{equation}\label{aaeq2345}
	\lim_{k\to\infty} \int_{K} |F^{o,\epsilon_k}(p,t)-G(p,t)|\dd p = 0
	\end{equation}
	and
	\begin{equation}\label{aaeq2346}
	\lim_{k\to\infty} \int_{K} |F^{o,\epsilon_k}(p,t)-\gamma_p(t)| \dd p = 0 .
	\end{equation}
The proof of \eqref{aaeq2345} is given in Lemma \ref{aalem2204} below and we next prove \eqref{aaeq2346}.\\

Notice that thanks to \eqref{aaeq2349} we have
	\[
	F^{o,\epsilon}(p,t) = F^{o,\epsilon}(p,0) + \int_0^t \sum_{j=1}^{\bar r} h_j(\delta_{\epsilon}^op,\epsilon s) \bar X_j(F^{o,\epsilon}(p,s))\dd s .
	\]
For $K=\delta^o_{\frac14}B_o$,
\begin{align*}
\int_{K}&\left| F^{o,\epsilon}(p,t) - g(p) - \int_0^t \sum_{j=1}^{\bar r} h_j(o,0) \bar X_j(G(p,s)) \dd s \right| \dd p\\
&\leq \int_{K} \left| F^{o,\epsilon}(p,0) - g(p)\right| \dd p \\
&\qquad+\sum_{j=1}^{\bar r} \int_{K}\int_0^t \left|
		 h_j(\delta_{\epsilon}^op,\epsilon s) \bar X_j(F^{o,\epsilon}(p,s)) -  h_j(o,0) \bar X_j(G(p,s)) \dd s \right| \dd p\\
&\leq \underbrace{ \int_{B_o} \left| F^{o,\epsilon}(p,0) - g(p)\right| \dd p }_{(a)}\\
&\qquad +\sum_{j=1}^{\bar r}\underbrace{ \int_{B_o} \int_0^t |h_j(\delta^o_\epsilon p,\epsilon s) - h_j(o,0)| \cdot |\bar X_j(F^{o,\epsilon}(p,s))| \dd s \dd p }_{(b)}\\
&\qquad\qquad\qquad +\sum_{j=1}^{\bar r}\underbrace{ \int_{K} \int_0^t |h_j(o,0)| \left| \bar X_j(F^{o,\epsilon}(p,s)) - \bar X_j(G(p,s)) \right| \dd s \dd p. }_{(c)}
\end{align*}
So we next estimate the three parts.

The proof of $(a)\to 0$ as $k\to \infty$ is given in Lemma \ref{aalem2351} below.
	
To estimate part $(b)$, notice that $|\bar X_j(F^{o,\epsilon}(p,s))|\le C$. We may use a change of variable and Proposition \ref{aaprop:Lebesg diff them} to infer that
	\begin{align*}
(b) \le  \frac C{\epsilon\cdot|\delta^o_\epsilon B_o|} \int_{\delta^o_\epsilon B_o} \int_0^\epsilon |h_j(p,s)-h_j(o,0)| \dd s \dd p
	\to 0
	\end{align*}
holds for almost every density point $o\in E$.	
Similarly, to estimate part $(c)$, we notice that $|h_j|\le C$ and hence
	\begin{align*}
	(c)\leq C\int_{K}& \int_0^t \left| \bar X_j(F^{o,\epsilon}(p,s)) - \bar X_j(G(p,s)) \right| \dd s \dd p.
	\end{align*}

Observe that $\bar X_j$ is Lipschitz on compact sets. So we may use Lemma \ref{aalem2204} below to conclude
	\begin{align*}
	\int_0^t \int_{K} &\left| \bar X_j(F^{o,\epsilon}(p,s)) - \bar X_j(G(p,s)) \right| \dd p \dd s \\
	&\le C\int_0^t \int_{K} \left|F^{o,\epsilon}(p,s) - G(p,s) \right| \dd p \dd s  \to 0,
	\end{align*}
from which \eqref{aaeq2346} follows.
\end{proof}
%%%%%%%%%%%%%%%%%%%%%%%%%%%%%%%%%%%%%%%%%%%%%%%%%%%%%%%%%%%%%%%
%%%%%%%%%%%%%%%%%%%%%%%%%%%%%%%%%%%%%%%%%%%%%%%%%%%%%%%%%%%%%%%
\begin{lemma}\label{aalem2351}
 	\[
	\lim_{k\to\infty} \int_{B_o} |F^{o,\epsilon_k}(p,0)-g(p)|\dd p = 0
	\]
\end{lemma}
\begin{proof}
	Since $o$ is a $\scr V$-density point of $E$, we have
 	\begin{align*}
	 \int_{B_o} &|F^{o,\epsilon}(p,0)-g(p)|\dd p\\
	 &\le  \int_{B_o\cap\delta^o_{\frac1\epsilon}E} |f^{o,\epsilon}(p)-g(p)|\dd p
	 	+  \int_{B_o\setminus\delta^o_{\frac1\epsilon}E} |F^{o,\epsilon}(p,0)-g(p)|\dd p \\
	&\le \int_{B_o\cap\delta^o_{\frac1\epsilon}E} |f^{o,\epsilon}(p)-g(p)|\dd p
		+ C |B_o\setminus\delta^o_{\frac1\epsilon}E|
	\to 0,
	\end{align*}
where we have used the fact that the Euclidean distance is controlled by $\bar {d}$.
\end{proof}
%%%%%%%%%%%%%%%%%%%%%%%%%%%%%%%%%%%%%%%%%%%%%%%%%%%%%%%%%%%%%%%
\begin{lemma}\label{aalem2204}
	For $t$ sufficiently close to 0, it holds
	 	\[
	\lim_{k\to\infty} \int_{\delta^o_{\frac14}B_o} |F^{o,\epsilon_k}(p,t)-G(p,t)|\dd p = 0.
	\]
\end{lemma}
\begin{proof}
Choose $t$ small enough such that
	$\phi_t(\delta^o_{\frac14} B_o)\subset \delta^o_{\frac12} B_o$
	and
	$(\phi_t)^{-1}\delta^o_{\frac12} B_o\subset  B_o$.
Then,
\begin{align*}
&\int_{\delta^o_{\frac14}B_o}|F^{o,\epsilon}(p,t)-G(p,t)|\dd p\\
&\le\int_{(\delta^o_{\frac14} B_o) \cap (\phi_t)^{-1}\delta^o_{\frac1\epsilon}E}
			|f^{o,\epsilon}(\phi_tp)-g(\phi_tp)|\dd p+\int_{(\delta^o_{\frac14} B_o) \setminus (\phi_t)^{-1}\delta^o_{\frac1\epsilon}E} |F^{o,\epsilon}(p,t)-G(p,t)|\dd p\\
\\\numberthis\label{aaeq:last equa}	
&\le C \int_{B_o\cap \delta^o_{\frac1\epsilon}E} |f^{o,\epsilon}(p)-g(p)|\dd p+C \int_{\delta^o_{\frac{1}{2}}B_o\setminus \delta^o_{\frac1\epsilon}E} |F^{o,\epsilon}({(\phi_t)}^{-1}p,t)-G({(\phi_t)}^{-1}p,t)|\dd p  \\
&\le C \int_{B_o\cap \delta^o_{\frac1\epsilon}E} |f^{o,\epsilon}(p)-g(p)|\dd p+C|B_o\setminus \delta^o_{\frac1\epsilon}E|,
\end{align*}
which converges to $0$ as $\epsilon_k\to0$.
\end{proof}

\subsection{Stepanov's theorem in Carnot groups}

\begin{proposition}[Rademacher's Theorem]
 	Let $E\subset G$ be a Borel set and $f:E\to\bar G$ a Lipschitz map.
	Then $f$ is $P$-differentiable almost everywhere in $E$.
\end{proposition}
Note that the limit in $P$-differentiability in the above proposition is with respect to only $X$ belonging to the Lie algebra such that $\exp (X)\in E$.
\begin{proof}
Without loss of generality, assume that $E$ is a closed set. For $k\in\{1,\dots,r\}$ and $j\in\{1,\dots,\bar r\}$, let $h_{kj}(p,t)$ be such that, setting
    \[
    F_k(p,t) = f(\exp(tX_k)(p)),
    \]
By extending $F_k$ as in the above subsection, we have
    \[
    \frac{\de F_k}{\de t}(p,t) = \sum_{j=1}^{\bar r} h_{kj}(p,t) \bar X_j(F_k(p,t)) .
    \]
 	Let $o\in E$ be a $\scr V$-density point of $E$.
	Almost all $o\in E$ have this properties.
	
	Set $\bar o=f(o)$.
	We want to define a morphism of Lie algebras $A:\g_o\to\g_{\bar o}$.
Set $\X:=\{tX_i:i\in\{1,\dots,r\},t\in\R\}$.
We view $\X$ as a subset of $\g_o$ by identifying $tX_i\in \X$ with $tX_i(o)\in \g_o=\G_o$. Remind that $\X$ generates $\G_o$ as a Lie group.

	Define $A$ on $\X$ as
	$$A[tX_k] := t\sum_{j=1}^{\bar r} h_{kj}(o,0) \bar X_j,$$
	where $t\sum_{j=1}^{\bar r} h_{kj}(o,0) \bar X_j$ can be viewed as an element of $\g_{\bar o}$ in the canonical way.

	Let $\epsilon_k\to0$ be a sequence such that $f^{o,\epsilon_k}$ converge to some $g$ as in Lemma~\ref{aa110}.
    By Proposition \ref{aaprop1044}, for almost every point $o\in E$, we have a neighborhood $K$ of the point $o$ and a constant $t_0$ such that
    \begin{equation}\label{aaeq2140}
    \exp(A[tX_j])(g(p)) = g(\exp(tX_j)(p)),
    \end{equation}
    for $p\in K$, $1\leq j\leq r$ and $|t|\leq t_0$.
   We identify $\g_o$ with $\G_o$ in the standard way. Let $V\in\g_o$. By Lemma 1.40 in \cite{SSS}, there are $\xi_1,\dots,\xi_\ell\in\X$ such that $V=\xi_1*\dots*\xi_\ell$, where $\ell$ is bounded by a constant and $*$ is the multiple in $\G_o$.
    Moreover, we have for all $p\in K$
    \[%\begin{multline*}
     	\exp(V)(p)
    	= \exp(\xi_1*\dots*\xi_\ell)(p)
    	= \exp(\xi_1)\circ\dots\circ\exp(\xi_\ell)(p).
    \]%\end{multline*}
    Set $\xi_i=t_iX_{j_i}.$
    Again, by Lemma 1.40 in \cite{SSS}, there is a constant $c_0$ such that if $\|V\|\leq c_0$, $|t_i|\leq t_0$ for any $1\leq i \leq\ell$.  
Furthermore, there is a smaller neighborhood $K^{\prime}$ of $o$ and a  constant $\epsilon_0^{\prime}$, such that for $p\in K^{\prime}$ and $\|V\|\leq \epsilon_0^{\prime}$ 
    \[\exp(\xi_j)\circ\dots\circ\exp(\xi_\ell)(p)\in K,\] for all $2\leq j\leq \ell$.   
In the following, fix $p\in K^{\prime}$ and $\|V\|\leq \min(\epsilon_0^{\prime},c_0,\epsilon_0)$. 
    Therefore, iterating \eqref{aaeq2140}, we have 
    \begin{equation}\label{aaeq2142}
    g\big(\exp(V)(p)\big) %= \exp(A\xi_1)\circ\dots\circ\exp(A\xi_\ell)\big(g(p)\big)
    = \exp(A\xi_1\bar*\dots\bar*A\xi_\ell)(g(p)),
    \end{equation}
    where $\bar*$ is the mlutiple in $\G_{\bar o}$.
	%The equation \eqref{eq2142} have many consequences.
If we take $p=o$, then \eqref{aaeq2142} becomes
    \begin{equation}\label{aa220}
    g(\exp(V)(o))
    = \exp(A\xi_1\bar*\dots\bar*A\xi_\ell)(\bar o)
    \end{equation}
    Since the right hand side does not depend on the sequence $\epsilon_k\to0$, $g$ is unique.

    On the other side, now $A$ can be extended to a map $\g_o\to\g_{\bar o}$.
    Indeed, if $V=\xi_1*\dots*\xi_\ell =\xi'_1*\dots*\xi'_{\ell'}$, then by \eqref{aaeq2142} the vector fields $A\xi_1\bar*\dots\bar*A\xi_\ell$ and $A\xi'_1\bar*\dots\bar*A\xi'_{\ell'}$ have the same action on the image of $g$, in particular they are same at $\bar o$, hence they are same everywhere.

    By construction, this extension of $A$ is a Lie algebra morphism and by Lemma \ref{aalem2002} it is the $P$-differential of $f$ at $o$.
\end{proof}

\begin{proof}[Proof of Theorem \ref{thm:stepanov Carnot group}]
It is easily to deduce that $L(f)=\cup_lF_l$, where \[F_l=\{x: \bar d(f(z),f(x))\leq ld(z,x),z\in B(x,1/l)\},\]
where $l\geq 1$. Next, we decompose the set $F_l$ into countable union $E_{l,i}$ such that $f|_{E_{l,i}}$ is $i$-Lipschitz, where integer $i\geq 1$. Moreover,  by the previous Rademacher type result, $f|_{E_{l,i}}$ is $ P$-differentiable at almost every point of $E_{l,i}$. Choose a point $x_0$ from them. According to the proof of the previous Rademacher type result, $x_0$ is a density point of $E_i$. Thus, as Lemma 3.15 in \cite{f69} and $f|_{E_{l,i}}$ is $ P$-differentiable at $x_0$, $f$ is also $ P$-differentiable at $x_0$. As $L(f)=\bigcup E_{l,i}$, we obtain that $f$ is $P$-differentiable at almost every point in $G$.
\end{proof}

\medskip 
\textbf{Acknowledgements}. The authors would like to thank Prof.~J. Tyson, Prof.~P. Koskela, Prof.~X. Xie, Prof.~E. Le Donne and Prof.~K. F\"assler for their very helpful discussions on related topics. 
%They are also grateful to Professor~S.K. Vodopyanov for kindly sharing the papers~\cite{v03,v07,v07diff}

\end{document}